\documentclass[aps,twocolumn,superscriptaddress,groupedaddress]{revtex4}  
\usepackage{graphicx}  
\usepackage{dcolumn}   
\usepackage{bm}        
\usepackage{overpic}   

\makeatletter
\renewcommand{\p@subsection}{}
\renewcommand{\p@subsubsection}{}
\makeatother

\usepackage{tikz}
\usetikzlibrary{matrix,arrows,positioning}
\usepgflibrary{fpu}

\usepackage{amsfonts,epsfig,amsmath,amsthm,amssymb,graphics,verbatim}

\newcommand\arrowonlineback[1]{node
[pos=#1,sloped,anchor=center,draw=none,fill=none,shape=rectangle,inner
   sep=0pt,outer sep=0pt] {\tikz\draw[<-,thick,dash pattern=on 0pt off 2pt]
(-1pt,0pt) -- (0pt,0pt);}} 

\usepackage{xcolor}
\usepackage[caption=false]{subfig}
\usepackage{booktabs}
\usepackage{tabularx}
\def\R{\mathbb{R}}
\def\C{\mathbb{C}}
\def\N{\mathbb{N}}

\def\P{\mathbb{P}}

\def\cA{\mathcal{A}}
\def\cB{\mathcal{B}}
\def\cC{\mathcal{C}}

\def\cK{\mathcal{K}}

\def\cO{\mathcal{O}}
\def\cP{\mathcal{P}}

\def\cX{\mathcal{X}}

\def\txta{{\textnormal{a}}}

\def\txtd{{\textnormal{d}}}
\def\txte{{\textnormal{e}}}

\def\txtr{{\textnormal{r}}}

\def\txtD{{\textnormal{D}}}

\newcommand{\rmd}{\mathrm{d}}

\newcommand{\fa}          {\quad \text{for all } \,}

\def\ra{\rightarrow}

\newcommand{\be}{\begin{equation}}
\newcommand{\ee}{\end{equation}}
\newcommand{\benn}{\begin{equation*}}
\newcommand{\eenn}{\end{equation*}}
\newcommand{\bea}{\begin{eqnarray}}
\newcommand{\eea}{\end{eqnarray}}
\newcommand{\beann}{\begin{eqnarray*}}
\newcommand{\eeann}{\end{eqnarray*}}

\newcommand{\eps}{\varepsilon}
\newcommand{\Ptrans}{P_{\text{trans}}}

\begin{document}

\title{A General View on Double Limits in Differential Equations}
\date{\today}

\author{Christian Kuehn}
\affiliation{Department of Mathematics, Technical University of Munich, 
Boltzmannstr.~3, 85748 Garching b.~M\"unchen, Germany}
\affiliation{Complexity Science Hub Vienna, Josefst\"adter Str.~39, 
1080 Vienna, Austria}

\author{Nils Berglund}
\affiliation{Institut Denis Poisson (IDP), Universit\'e d'Orl\'eans, 
Universit\'e de Tours, CNRS -- UMR 7013, B\^atiment de Math\'ematiques, B.P. 
6759, 45067~Orl\'eans Cedex 2, France}

\author{Christian Bick}
\affiliation{Department of Mathematics, Vrije Universiteit Amsterdam, 
De~Boelelaan 1111, Amsterdam, the Netherlands}
\affiliation{Department of Mathematics, University of Exeter, Exeter 
EX4 4QF, UK}
\affiliation{Institute for Advanced Study, Technical University of Munich, 
Lichtenbergstr.~2, 85748 Garching, Germany}

\author{Maximilian Engel}
\affiliation{Department of Mathematics, Freie Universtit{\"a}t Berlin, 
Arnimallee~6, 14195 Berlin, Germany}

\author{Tobias Hurth}
\affiliation{Institut de math\'ematiques, Universit\'e de Neuch\^atel, 
Rue Emile-Argand 11, CH-2000 Neuch\^atel} 

\author{Annalisa Iuorio}
\affiliation{Faculty of Mathematics, University of Vienna,
Oskar-Morgenstern-Platz 1, 1090, Vienna, Austria}

\author{Cinzia Soresina}
\affiliation{Institute for Mathematics and Scientific Computing, 
University of Graz, Heinrichstr.~36, 8010 Graz, Austria}

\begin{abstract}
In this paper, we review several results from singularly perturbed differential equations with multiple small parameters. In addition, we develop a general conceptual framework to compare and contrast the different results by proposing a three-step process. First, one specifies the setting and restrictions of the differential equation problem to be studied and identifies the relevant small parameters. Second, one defines a notion of equivalence via a property/observable for partitioning the parameter space into suitable regions near the singular limit. Third, one studies the possible asymptotic singular limit problems as well as perturbation results to complete the diagrammatic subdivision process. We illustrate this approach for two simple problems from algebra and analysis. Then we proceed to the review of several modern double-limit problems including multiple time scales, stochastic dynamics, spatial patterns, and network coupling. For each example, we illustrate the previously mentioned three-step process and show that already double-limit parametric diagrams provide an excellent unifying theme. After this review, we compare and contrast the common features among the different examples. We conclude with a brief outlook, how our methodology can help to systematize the field better, and how it can be transferred to a wide variety of other classes of differential equations.
\end{abstract}

\maketitle

\section{Introduction}

Effectively all problems arising from science and engineering are studied by only considering a suitably reduced model of reality. In particular, we would often like to reduce differential equations by assuming that certain physical effects or external influences do not play a major role for the scientific question of interest. Yet, this implicitly supposes we can also show that the terms we do neglect are in some sense ``small'' so that they do not change the answers to the relevant scientific questions. There is a vast number of differential equations where \textcolor{black}{direct approaches to remove small parameters fail and non-trivial correction terms appear when perturbing away from the limit. These differential equations are often called singularly perturbed.} \textcolor{black}{A single generally accepted definition of ``singularly perturbed'' does not exist as some definitions are too narrow, others are too broad. Here, we adopt a pragmatic approach and define a singularly-perturbed differential equation as one where taking the small parameter to zero yields a differential equation within a different structural class. Doubly-singular differential equations are then those, where two small parameters lead each in the singular limit to a different structural problem class.} \textcolor{black}{From a practical viewpoint, the first step is to identify the origins of small parameters} which control the strength of the terms we want to neglect. Some typical examples appearing in the context of differential equations are:  

\begin{itemize}
\item \textit{Time Scale Separation:} Two, or more, sets of variables evolve at different rates.
\item \textit{Noise Level:} Finite-size effects or external forces are modeled via noise.
\item \textit{Spatial Scale Separation:} Two, or more, sets of variables have differing spatial scales.
\item \textit{Network Coupling:} Operating a system within a network leads to new coupling dependencies. 
\end{itemize}

In this paper, we are going to focus on these areas to illustrate the types of results
one can obtain for (multiple) small parameters. Of course, there are many other areas in
differential equations, where small parameters appear, for example:

\begin{itemize}
\item \textit{Discretization Size:} Temporal and/or spatial discretization leads to small parameters.
\item \textit{Inverse Particle Number:} One wants to convert finite systems to a continuum model.
\item \textit{Interfaces:} Interfaces or boundary layers are often small.
\item \textit{Nonlocal Coupling:} Local derivatives are augmented by global integral terms.
\item \textit{Nonsmoothness:} Functions are taken smooth outside of small subsets of space.
\item \textit{Time Delay:} (Small) communication delay induces a time history dependence.
\item \textit{Near-Symmetry:} A system might be very close to a symmetric one.
\item \textit{Near-Integrability:} Perturbations of integrable and/or Hamiltonian systems are well-studied.
\end{itemize}

Even the combination of the two previous lists is just a restricted snapshot of all potential 
cases where small parameters may appear. From a historical perspective, small parameters in differential equations are a quite classical topic that can be traced back at least to the end of the 19th century. Among the first applications were celestial mechanics~\cite{Poincare} and fluid dynamics~\cite{Prandtl}. In celestial mechanics, since the two-body problem is solvable, the three-body problem lends itself to consider singular perturbations by assuming two large mass bodies and one very small mass. In fluid mechanics, assuming very large viscosity is helpful as this assumption usually precludes the existence of turbulent flow. In the limited space of this work, it is impossible to give proper credit to the very successful, long, and sometimes winding, history of singular perturbations in celestial mechanics and fluid dynamics, so we refer to \cite{Beutler,Holmes5,OMalley24,VanDyke} containing excellent historical accounts and references regarding the development of these areas. 

Within the 
20th century, the use of small parameters and perturbation techniques for differential equations \textcolor{black}{has} permeated effectively all areas of science and engineering, while more recently also quantitative modelling in the social sciences tends to rely on differential equation modelling. For some pointers to the vast literature, we refer to the books~\cite{BenderOrszag,BensoussanLionsPapanicolaou,DeJagerFuru,Holmes,Jones,KevorkianCole,KuehnBook,Nayfeh1,OMalley20,PavliotisStuart,SandersVerhulstMurdock,Verhulst,Wasow,Wechselberger4}, where classical cases of ordinary and partial differential equations (ODEs and PDEs) with one small parameter are considered from a number of different viewpoints. \textcolor{black}{These books also contain several variants and viewpoints on the definition of ``singular perturbation'' for ODEs and provide an outlook to the PDE case.} 

Although the literature is quite detailed, it has become apparent in recent years that several techniques have to be extended to deal with more complex 21st century challenges, where differential equations and small parameters still take center stage. First, one might wonder, why existing methods have to be developed further? 
The first key reason is that mathematical modeling of complex systems almost immediately dictates that the case of just one small parameter is very rare. For example, it would be very difficult to argue that global climate dynamics, socio-economic networked systems, or neuro-mechanical as well as systems biology problems, frequently contain just \emph{one} small parameter. Second, in complex systems we often deal with many instabilities. Each instability, even if it is localized in parameter and phase space, leads to a delicate balance between nonlinear terms. Hence, we cannot invoke simple principles that very stable leading-order linear terms dominate so that small contributions from external/internal model perturbations are irrelevant. This entails the need for larger phase and parameter spaces~\cite{KuehnCurse}. In summary, there is an imminent need to study the case of two or more small parameters carefully to obtain a good practical understanding of current important topics in differential equations. \textcolor{black}{More precisely, we will restrict here the focus on analyzing differential equations, where two small parameters $(\varepsilon,\delta)$ tend to zero from above, and we want to classify different scaling regimes for this double limit.}

As one might expect, this field also has \textcolor{black}{an intricate history within several sub-disciplines of differential equations being involved. This makes it often difficult to gain access and/or an overview, when studying double limits. The} most classical cases, where two small parameters have been analyzed first, were ODEs with a focus on direct asymptotic methods such as matching~\cite{Freund,Meyer,OMalley25,OMalley11}, although more recently also more geometric ODE approaches have gained popularity, see e.g.~\cite{DeMaesschalckDumortier3,DeMaesschalckWechselberger,KrupaPopovicKopell,KrupaPopovicKopellRotstein,KuehnSzmolyan,CardinTeixeira}. Although extensions of existing approaches are often key components for our understanding of multiple small parameters, the development is not nearly as systematic and detailed as for just one distinguished small parameter. One can view the situation in analogy with several other areas of differential equations, e.g., second-order scalar oscillators already show a lot of interesting behaviour, but eventually one has to go beyond a widely accepted standard class. Therefore, we believe it is now time to re-think and systematize double limits in 
differential equations. In fact, virtually within all areas of differential equations, multiple small parameters do appear. In this review, we try to reflect this broader perspective via several illustrating examples motivated by very different applications. We are going to describe many key challenges, where a naive direct approach of taking double limits fails.\medskip 

More precisely, a common, yet highly non-trivial, situation we want to understand are doubly-singularly perturbed differential equations, or more generally multiscale dynamics with multiple small parameters. As argued above, a unified framework to understand doubly-singular perturbations is still lacking, so this will be our starting point. Here we make a conceptual step towards improving this situation.\medskip

\begin{figure}[htbp]
	\centering
\begin{tikzpicture}
[>=stealth',point/.style={circle,inner sep=0.035cm,fill=white,draw},
encircle/.style={circle,inner sep=0.07cm,draw},
x=4cm,y=4cm,declare function={f(\x) = sqrt(\x);}]

\draw[->,semithick] (0,0) -> (1,0);
\draw[blue,very thick] (0,0) -> (0.9,0);
\draw[->,thick,dashed] (0,0) -> (0,1);

\draw[red,semithick,-,smooth,domain=0:0.9,samples=75,/pgf/fpu,/pgf/fpu/output
format=fixed] plot (\x, {0.5*f(\x)});

\draw[blue,very thick,-,smooth,domain=0:0.9,samples=75,/pgf/fpu,/pgf/fpu/output
format=fixed] plot (\x, {f(\x)});

\node[point] at (0,0) {};

\node[] at (0.92,-0.05) {$\eps$};
\node[] at (-0.05,0.9) {$\delta$};

\draw [black] (0.2,0.72) circle [radius=0.06];
\node at (0.2,0.72) {III};
\draw [black] (0.6,0.55) circle [radius=0.06];
\node at (0.6,0.55) {II};
\draw [black] (0.7,0.25) circle [radius=0.06];
\node at (0.7,0.25) {I};

\end{tikzpicture}
\vspace{-4mm}	
	\caption{\label{fig:01}\textcolor{black}{Partitioning} of the positive quadrant 
	$\cK$ near the doubly-singular limit $\varepsilon\ra 0$ and $\delta\ra 0$ 
	into three different regions (I)--(III), which are non-equivalent
	under a property $\cP$. The thick lines (in blue) indicate hard boundaries 
	between the different regions, e.g.~between (II) and (III) there is a precise 
	curve separating these regions. The thin line (in red) indicates that the boundary 
	is only asymptotic up to a constant between two regions. Dashed lines (in black) 
	indicate an unclassified axis (such as the vertical axis in this figure). The circle 
	at the origin also means that at this point a classification with respect to $\cP$ is 
	not known and/or may not even be possible.}
\end{figure}
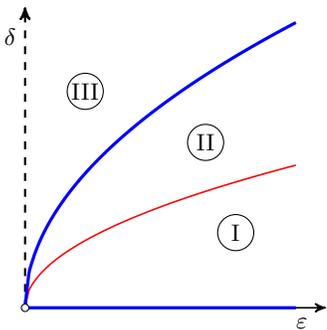 

Consider a doubly-singularly perturbed differential equation with two small 
non-negative parameters~$\varepsilon$ and~$\delta$. Often we are interested 
in the local behavior of the differential equation in the cone
\benn
\cK:= \{(\varepsilon,\delta)\in\R^2:\varepsilon\geq 0,\delta\geq 0\}
\eenn
intersected with a sufficiently small ball around the origin\textcolor{black}{, i.e., suppose we have tried already to neglect the small parameters but setting $\varepsilon=0=\delta$ does not provide a suitable description of the dynamics. Hence, }the natural step is to try to partition $\cK$ into different regions as shown in Figure~\ref{fig:01}. To make such a partitioning precise, we propose several steps:

\begin{enumerate}
 \item[(S1)] Specify the setting and restrictions of the problem $\cX$ to be studied.
 \item[(S2)] Define a notion of equivalence via a property/observable $\cP$ for the partitioning.
 \item[(S3)] Study the possible asymptotic limit problems $\cA$ to complete the diagram. 
\end{enumerate}

In the available literature, these steps can be found in various incarnations and various levels
of mathematical rigor. What tends to be missing in many problems is to \emph{recognize} (S1)--(S3) in 
a clear way to allow for a more \emph{comparative} and \emph{systematic} classification of possible 
behaviors. Already very simple classical examples, as discussed in Section~\ref{sec:classical}, show
that missing small details or slight changes in the setting $\cX$ or definition $\cP$ in the steps 
(S1)-(S2) can lead to completely different answers. We are going to show in this work that if the 
steps (S1)--(S3) are carried out carefully and within a uniform framework, a surprisingly coherent 
picture emerges, how doubly-singularly perturbed differential equations can be studied. The 
cross-connections between different classes of effects and methods thus become more visible. Universal 
classification diagrams emerge that concisely make the differences and similarities between different 
sub-fields of differential equations \textcolor{black}{much more prominent. Of course, we are still relying on well-established methods to carry out certain proofs or numerical explorations, particularly in step (S3), where the common viewpoint of singular perturbation theory to utilize the singular limit $\varepsilon=0=\delta$ takes center stage to understand scaling relations for $0<\varepsilon,\delta\ll 1$. From the viewpoint of singularity/bifurcation theory for ODEs, this often  means one is trying to unfold the dynamics in a suitable neighbourhood of a singular point. Yet, the key point is to always take into account, how $\cX$, $\cP$, $\cA$ are defined, which may depend crucially on the question and/or application. Indeed, this leads us beyond the notion of standard ODE classification via topological equivalence, which is not sufficient to fully understand double limits for different classes of differential equations.} \textcolor{black}{In summary, we contribute to provide a better starting point for }a systematic study of 
doubly-singular limits as another unifying scientific principle in the analysis of differential equations.\medskip

The remaining part of this paper is structured as follows: In Section~\ref{sec:classical}, we 
explain our approach via simple examples from analysis and algebra without a direct reference
to differential equations. The core part of this work is contained in Section~\ref{sec:doublysingular}, 
where numerous classes of known results for differential equation problems are re-cast precisely 
in the three steps (S1)--(S3) to provide a general framework, which highlights the unity of area. This
includes problems from fast-slow ODE dynamics, small noise stochastic differential 
equations (SDEs) and piecewise deterministic Markov processes (PDMPs), spatial problems arising from the
bifurcation analysis of partial differential equations (PDEs), and a problem in network
dynamics. In Section~\ref{sec:comparison}, we then contrast and compare the results. 
Section~\ref{sec:outlook} provides an outlook towards a more systematic study of multi-parameter 
singular limits for differential equations.

\section{Classical Examples}
\label{sec:classical}

Before starting with the development of a singular limit analysis of various classes of differential
equations, we illustrate some basic principles that occur in the steps (S1)--(S3) in simpler settings.

\subsection{Elementary Algebra}
\label{ssec:algebra}

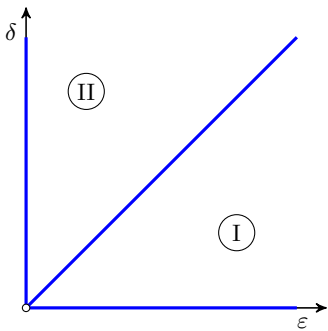
\begin{figure}[htbp]
	\centering
\begin{tikzpicture}
[>=stealth',point/.style={circle,inner sep=0.035cm,fill=white,draw},
encircle/.style={circle,inner sep=0.07cm,draw},
x=4cm,y=4cm,declare function={f(\x) = sqrt(\x);}]

\draw[->,semithick] (0,0) -> (1,0);
\draw[blue,very thick] (0,0) -> (0.9,0);
\draw[->,semithick] (0,0) -> (0,1);
\draw[blue,very thick] (0,0) -> (0,0.9);
\draw[blue,very thick] (0,0) -> (0.9,0.9);

\node[point] at (0,0) {};

\node[] at (0.92,-0.05) {$\eps$};
\node[] at (-0.05,0.92) {$\delta$};

\draw [black] (0.2,0.72) circle [radius=0.06];
\node at (0.2,0.72) {II};
\draw [black] (0.7,0.25) circle [radius=0.06];
\node at (0.7,0.25) {I};

\end{tikzpicture}
\vspace{-4mm}	
	\caption{\label{fig:02}Classification diagram with respect to the property $\cP_{[-1,1]}$. 
	In region II we have no zeros while in region I we have two zeros (counting multiplicity).}
\end{figure}

Consider the root-finding problem of a very simple quadratic polynomial
\benn
\label{eq:roots} \tag{{$\cX_{\textnormal{rts}}$}}
f(x;\varepsilon,\delta):=\varepsilon x^2 - \delta \stackrel{!}{=} 0. 
\eenn
For the problem~\eqref{eq:roots}, we assume that we do not allow any coordinate changes and/or
preliminary algebraic scaling operations for the problem, i.e., we want to find the roots as is. 
For any $\varepsilon,\delta>0$, we have the roots $x_\pm=\pm\sqrt{\delta/\varepsilon}$. Now it
crucially depends on the choice of the property $\cP$ what a classification diagram in a form
similar to Figure~\ref{fig:01} would look like. Suppose we take as a definition that two problems of the
form~\eqref{eq:roots} are equivalent if they have the same property
\benn
\cP_{[-1,1]}:=\text{cardinality}\{x\in[-1,1]:f(x;\varepsilon,\delta)=0\},
\eenn  
where we count roots according to multiplicity. Then one just calculates 
$|x_\pm|^2=\delta/\varepsilon\leq 1$ which yields $\delta\leq \varepsilon$. 
Hence, there are just two regions in the $(\varepsilon,\delta)$-diagram separated by the 
diagonal $\{\delta=\varepsilon\}\cap \cK$. Above the diagonal, we have $\delta>\varepsilon$
so $\cP_{[-1,1]}|_{\delta>\varepsilon}=0$, while on or below the diagonal we have 
$\cP_{[-1,1]}|_{\delta\leq\varepsilon}=2$. Of course, the point at the origin is special
leading to a solution set which is uncountable so we decide to leave it out in our classification;
see Figure~\ref{fig:02}. 

The splitting into two main regions is also visible via considering the two singular limit 
problems of~\eqref{eq:roots}, namely
\be
\lim_{\varepsilon\ra 0}f(x;\varepsilon,\delta)=- \delta \stackrel{!}{=} 0,
 \tag{{$\cA_{\textnormal{rts}}^{\varepsilon=0}$}}
\ee
and
\be
\lim_{\delta\ra 0}f(x;\varepsilon,\delta)= \varepsilon x^2 \stackrel{!}{=} 0, 
\tag{{$\cA_{\textnormal{rts}}^{\delta=0}$}}
\ee  
where we get no roots and a double-root respectively. In summary, there is also an inherent 
non-commutativity in the limits. Yet, the precise setting of~\eqref{eq:roots} and the 
specification $\cP_{[-1,1]}$ are crucial. For example, if we use $\cP_{\R}$ instead, looking
for all the real roots, then there is only one singular line remaining in parameter
space given by $\{\varepsilon=0,\delta>0\}$ with no roots and the usual singular situation 
at the point $(\varepsilon,\delta)=(0,0)$. Also, given the function 
\be
f(x;\varepsilon,\delta):=\varepsilon x^2-\delta 
\ee
we could have used a completely different property $\cP$ to check for equivalence. For 
example, we could ask for a binary classification and set
\be
\cP_{\textnormal{cvx}}=\left\{\begin{array}{ll}
1\qquad & \textnormal{if $f$ is convex in $x$},\\
0\qquad & \textnormal{if $f$ is not convex in $x$.}
\end{array}\right.
\ee

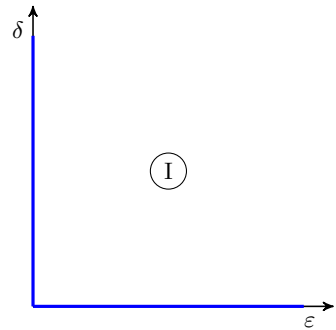
\begin{figure}[htbp]
	\centering
\begin{tikzpicture}
[>=stealth',point/.style={circle,inner sep=0.035cm,fill=white,draw},
encircle/.style={circle,inner sep=0.07cm,draw},
x=4cm,y=4cm,declare function={f(\x) = sqrt(\x);}]

\draw[->,semithick] (0,0) -> (1,0);
\draw[blue,very thick] (0,0) -> (0.9,0);
\draw[->,semithick] (0,0) -> (0,1);
\draw[blue,very thick] (0,0) -> (0,0.9);


\node[] at (0.92,-0.05) {$\eps$};
\node[] at (-0.05,0.92) {$\delta$};

\draw [black] (0.45,0.45) circle [radius=0.06];
\node at (0.45,0.45) {I};

\end{tikzpicture}
\vspace{-4mm}	
	\caption{\label{fig:03}Classification diagram with respect to the property $\cP_{\textnormal{cvx}}$. 
	We just have a single region as $f(x;\varepsilon,\delta)=\varepsilon x^2-\delta$ is always convex
	on $\cK$.}
\end{figure} 

Then we always have $\cP_{\textnormal{cvx}}|_{(\varepsilon,\delta)\in\cK}=1$ so the singular
limit classification is somewhat trivial as shown in Figure~\ref{fig:03}. This demonstrates 
that, although many a-priori natural-looking
mathematical properties could be used for double limits, it is vital to have a good 
motivation from applications and modeling to select the most important ones.

\subsection{Elementary Analysis}
\label{ssec:analysis}

The issues illustrated in the last section are evidently not limited to just purely algebraic 
problems. For example, let us consider the classical function
\be
\label{eq:Clairaut}\tag{{$\cX_{\textnormal{partials}}$}}
\tilde{f}(x,y):=\left\{
\begin{array}{ll}
\frac{xy(x^2-y^2)}{x^2+y^2}\qquad &\textnormal{if $(x,y)\neq (0,0)$,}\\
0 \qquad &\textnormal{if $(x,y)= (0,0)$,}
\end{array}
\right.
\ee
which is known to be a simple counter-example in the context of Schwarz's Theorem since the partial 
derivatives do not commute at zero
\be
\label{eq:failClairaut}
-1=\partial_{xy}\tilde{f}(0,0)\neq\partial_{yx}\tilde{f}(0,0)=1.
\ee
Evidently, we can also just understand this via double limits in defining
\benn
\begin{array}{lcr}
f(x,y;\varepsilon,\delta)&:=&\frac{\tilde{f}(x+\varepsilon,y+\delta)
-\tilde{f}(x,y+\delta)}{\varepsilon\delta}\\
&&+ \frac{\tilde{f}(x,y)-\tilde{f}(x+\varepsilon,y)}{\varepsilon\delta},
\end{array}
\eenn
and then~\eqref{eq:failClairaut} just means that 
\benn
\lim_{\delta\ra 0}\lim_{\varepsilon\ra 0}f(x,y;\varepsilon,\delta)\neq 
\lim_{\varepsilon\ra 0}\lim_{\delta\ra 0}f(x,y;\varepsilon,\delta)
\eenn
if we evaluate the two limits at $(x,y)=(0,0)$.
Evidently the subdivision of the cone $\cK$ again depends crucially on the choice of the property
$\cP$. However, here we shall fix the relevant property via second partial derivatives below. 

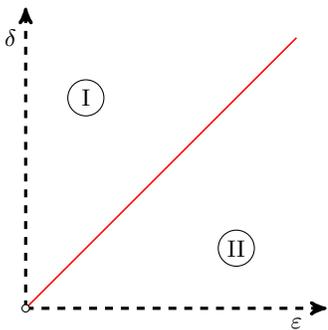
\begin{figure}[htbp]
	\centering
\begin{tikzpicture}
[>=stealth',point/.style={circle,inner sep=0.035cm,fill=white,draw},
encircle/.style={circle,inner sep=0.07cm,draw},
x=4cm,y=4cm,declare function={f(\x) = sqrt(\x);}]

\draw[->,very thick,dashed] (0,0) -> (1,0);
\draw[->,very thick,dashed] (0,0) -> (0,1);
\draw[red,semithick] (0,0) -> (0.9,0.9);

\node[point] at (0,0) {};

\node[] at (0.9,-0.05) {$\eps$};
\node[] at (-0.05,0.9) {$\delta$};

\draw [black] (0.2,0.7) circle [radius=0.06];
\node at (0.2,0.7) {I};
\draw [black] (0.7,0.2) circle [radius=0.06];
\node at (0.7,0.2) {II};

\end{tikzpicture}
\vspace{-4mm}	
	\caption{\label{fig:04}Classification diagram with respect to the property $\cP_{\partial\partial}$.
	The two regions correspond to the two possible partial derivative values at the origin of
	the function \eqref{eq:Clairaut} given for our elementary analysis problem. The thin line (red) could 
	have been chosen at any fixed slope as it is an asymptotic subdividing line of the form $\{\delta=
	\kappa \varepsilon,\varepsilon>0\}$ for some fixed constant $\kappa>0$.}
\end{figure}  

Since we are in an analytic setting, and not in an algebraic one, it often makes 
sense not to aim for a point-wise subdivision of the cone $\cK$. Instead, we are 
going to use an asymptotic subdivision by assuming that $\delta=\delta(\varepsilon)$ with
$\delta\in\C^0(\R^+_0,\R^+_0)$ and $\delta(0)=0$, which just means that $\delta$ is a 
continuous function of $\varepsilon$ vanishing simultaneously. If we define 
\benn
\cP_{\partial\partial}:=\lim_{\varepsilon\ra 0}f(x,y;\varepsilon,\delta(\varepsilon))|_{(x,y)=(0,0)}
\eenn 
then there are two main regions in $\cK$. Either we have $\delta(\varepsilon)=o(\varepsilon)$ as 
$\varepsilon\ra 0$ (alternatively written $\delta\ll \varepsilon$), which yields $\cP_{\partial\partial}=-1$. 
Or we have $\delta\gg \varepsilon$ leading to $\cP_{\partial\partial}=+1$. Hence, it is natural to divide 
$\cK$ into two regions via a line $\delta=\kappa \varepsilon$ for a fixed constant $\kappa>0$. The 
constant $\kappa$ is somewhat arbitrary as long it is independent of $\varepsilon$ and $\delta$ so we 
just write for the codimension-one subdivision line $\cK\cap \{\varepsilon\simeq \delta\}$; see
Figure~\ref{fig:04}. 

\section{Doubly-Singular Systems}
\label{sec:doublysingular}

As a next step, it is important to demonstrate that different classes of doubly-singularly 
perturbed differential equations fit within and benefit from the more unified view described
so far. We shall illustrate this aspect with several very recent examples, where one cannot
only re-cast the problem within our framework but where the main strategy and effects become
very transparent as a result.
 
\subsection{Multiple Time Scale Systems}
\label{ssec:mts}

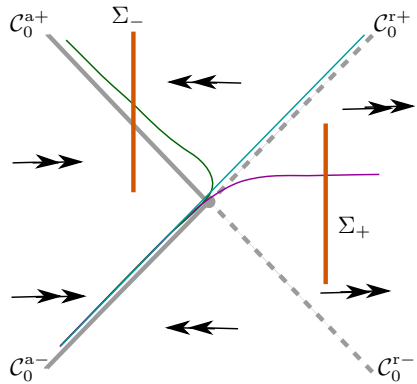
\begin{figure}[htbp]
	\centering
	\begin{overpic}[width=0.3\textwidth]{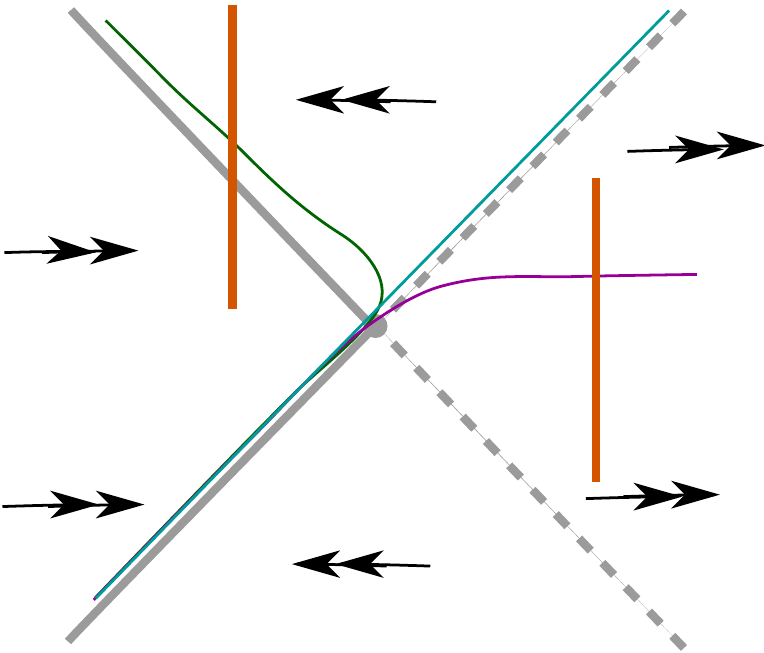}
	\put(0,85){$\cC_0^{\txta+}$}
	\put(0,0){$\cC_0^{\txta-}$}
  \put(90,85){$\cC_0^{\txtr+}$}
	\put(91,0){$\cC_0^{\txtr-}$}
	\put(25,87){$\Sigma_-$}
	\put(81,35){$\Sigma_+$}
	\end{overpic}
	\caption{\label{fig:05}Sketch of the possible dynamics of~\eqref{eq:tc} in $(x,y)$-coordinates. 
	The critical manifold $\cC_0$ is shown in gray (repelling parts with dashed lines and
	attracting parts with solid lines). Three trajectories (green, cyan, magenta) for 
	$0<\varepsilon\ll1$ are indicated for three different choices of $\delta$ (corresponding 
	to the exchange-of-stability, canard, and critical transition cases respectively). Double
	arrows show the direction of the fast subsystem flow for orientation; the slow subsystem 
	dynamics on $\cC_0$ is always directed upwards at unit speed.}
\end{figure}

We start with arguably one of the most classical~\cite{BenderOrszag,DeJagerFuru,OMalley20} cases of singular perturbation problems~\cite{OMalley24,Verhulst}, namely ordinary differential equations (ODEs) with two time scales, so-called fast-slow systems~\cite{Jones,Kaper,KuehnBook,DeMaesschalckDumortierRoussarie}. A good illustration within this context is to consider the transcritical fast-slow bifurcation normal form
\be
\label{eq:tc}\tag{{$\cX_{\textnormal{tc}}$}}
\begin{array}{rcrcl}
\frac{\txtd x}{\txtd t} &=& x' &=&(x-y)(x+y) + \frac{\varepsilon^2}{\delta},\\
\frac{\txtd y}{\txtd t} &=& y' &=& \varepsilon,
\end{array}
\ee
which is a well-studied system~\cite{KruSzm4}. As before, we shall assume that $\varepsilon\geq 0$ 
is a small parameter and then consider the case when $\delta\geq 0$ is a second
small parameter. Taking the fast subsystem limit~\eqref{eq:tc} given by $\varepsilon\ra 0$ 
yields
\be
\label{eq:t1}\tag{{$\cA_{\textnormal{tc,f}}^{\varepsilon=0}$}}
\begin{array}{rcl}
x' &=&x^2-y^2,\\
y' &=&0,
\end{array}
\ee
which is just a standard transcritical bifurcation with the slow variable $y$ acting as 
a bifurcation parameter. If we re-scale time as $s:=\varepsilon t$ and take the singular limit 
$\varepsilon\ra0$ again, then one obtains the slow subsystem
\be
\label{eq:t2}\tag{{$\cA_{\textnormal{tc,s}}^{\varepsilon=0}$}}
\begin{array}{rcl}
0 &=&(x-y)(x+y),\\
\frac{\txtd y}{\txtd s} &=& 1.
\end{array}
\ee
The fast and slow subsystems~\eqref{eq:t1}-\eqref{eq:t2} already show a singular structure
as the systems are differential equations of different types\textcolor{black}{, i.e., we go from a differential equation to a parameterized differential equation and differential algebraic equation respectively}. The algebraic constraint within the
slow subsystem is given by the critical manifold
\benn
\cC_0:=\{(0,0)\}\cup \cC_0^{\txta-}\cup \cC_0^{\txta+} \cup \cC_0^{\txtr-} \cup \cC_0^{\txtr+}
\eenn
where $\cC_0^{\txta-}:=\{|x|=|y|,x<0,y<0\}$, $\cC_0^{\txta+}:=\{|x|=|y|,x<0,y>0\}$,
$\cC_0^{\txtr-}:=\{|x|=|y|,x>0,y<0\}$, and $\cC_0^{\txtr+}:=\{|x|=|y|,x>0,y>0\}$ are normally hyperbolic
since the linearization with respect to the fast variables yields $\txtD_x(x^2-y^2)=2x$, which
is nonzero on $\cC_0\setminus\{(0,0)\}$. The critical manifold $\cC_0$ consists of equilibrium
points for the fast subsystem; see also Figure~\ref{fig:05}. Fenichel Theory~\cite{Fenichel4,Jones,KuehnBook} 
implies that there exist associated slow manifolds $\cC_\varepsilon^{\txta\pm}$ and 
$\cC_\varepsilon^{\txtr\pm}$. 

\vspace{0.3cm}

A generally very important question in many applications is
how trajectories of fast-slow systems pass through the region of a transcritical bifurcation
of the fast subsystem; for example, there are applications in 
ecology~\cite{BoudjellabaSari,KooiPoggialeAugerKooijman}, chemistry~\cite{KuehnUM}, numerical 
analysis~\cite{EngelKuehn}, epidemiology~\cite{Jardonetal} and network science~\cite{Jardon-Kojakhmetov2020}. 
Suppose we start with a trajectory $\gamma=\gamma(t)$ at 
a typical point on the attracting critical manifold $\cC_0^{\txta-}$, say $\gamma(0)=(x(0),y(0))
=(-3,-3)$ for concreteness as the following arguments do not change up to scaling by fixed constants. 
By Fenichel Theory, we have that $\gamma(0)$ is $\cO(\varepsilon)$-close to the slow manifold 
$\cC^{\txta-}_\varepsilon$. We are going to define two one-dimensional sections:
\benn
\Sigma_-:=\{x=-2,y\in[1,3]\},\quad \Sigma_+:=\{x=2,y\in[-1,1]\}. 
\eenn
One may easily prove that the trajectory $\gamma$ will first get attracted to 
$\cC^{\txta-}_\varepsilon$ exponentially fast and then track this manifold up towards the origin
due to the slow dynamics. Then there are three cases~\cite{KruSzm4}:

\begin{itemize}
 \item[(I)] If $\delta(\varepsilon)=\varepsilon(1+\cO(|\varepsilon|^p))$, 
for some $p>0$, then the trajectory will intersect $\Sigma_-$.
 \item[(II)] If $\delta(\varepsilon)=\varepsilon(1-\cO(|\varepsilon|^p))$, 
for some $p>0$, then the trajectory will intersect $\Sigma_+$.
 \item[(III)] If $\delta(\varepsilon)=\varepsilon(1\pm\cO(\exp(-K/\varepsilon)))$, 
then the trajectory will never intersect $\Sigma_\pm$.
\end{itemize} 

\begin{figure}[htbp]
	\centering
\begin{tikzpicture}
[>=stealth',point/.style={circle,inner sep=0.035cm,fill=white,draw},
encircle/.style={circle,inner sep=0.07cm,draw},
x=4cm,y=4cm,declare function={f(\x) = 1.5*exp(-1.75/\x);}]

\draw[->,semithick] (0,0) -> (1,0);
\draw[blue,very thick] (0,0) -> (0.9,0);
\draw[->,semithick] (0,0) -> (0,1);
\draw[blue,very thick] (0,0) -> (0,0.9);

\draw[blue,very thick] (0,0) -> (0.7,0.7);
\draw[red,semithick,-,smooth,domain=0.01:0.75,samples=75,/pgf/fpu,
/pgf/fpu/output format=fixed] plot (\x, {\x + f(\x)});
\draw[red,semithick,-,smooth,domain=0.01:0.75,samples=75,/pgf/fpu,
/pgf/fpu/output format=fixed] plot (\x, {\x - f(\x)});

\node[point] at (0,0) {};

\node[] at (0.92,-0.05) {$\eps$};
\node[] at (-0.05,0.92) {$\delta$};

\draw [black] (0.8,0.8) circle [radius=0.06];
\node at (0.8,0.8) {III};
\draw [black] (0.2,0.72) circle [radius=0.06];
\node at (0.2,0.72) {II};
\draw [black] (0.7,0.25) circle [radius=0.06];
\node at (0.7,0.25) {I};

\end{tikzpicture}
\vspace{-4mm}	
	\caption{\label{fig:06}Classification diagram with respect to the property $\cP_{\textnormal{tcd}}$
	for the problem~\eqref{eq:tc}.
	The three regions correspond to the cases (I)--(III) above yielding exchange-of-stability, critical
	transition, and canard cases respectively.}
\end{figure}
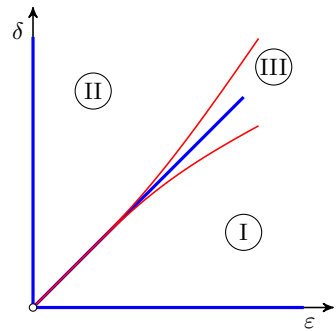

This classification is important as in case (I) we have an exchange-of-stability as $\gamma$ 
starts to track the attracting slow manifold $\cC^{\txta+}_\varepsilon$, while in case (II), there is
a critical transition leading to a jump near the fast subsystem bifurcation point.
In case (III), we have that $\gamma$ starts to track the repelling branch $\cC^{\txtr+}_\varepsilon$
for a slow time of order $\cO(1)$, which means that we have a canard 
orbit~\cite{BenoitCallotDienerDiener,DeMaesschalckDumortierRoussarie,KuehnBook}. Hence,
since these three cases differ crucially for application purposes, it makes sense to define a property 
$\cP_{\textnormal{tcd}}$ by a variable having just three possible values corresponding to the cases 
(I)-(III) respectively. This provides us with the double singular limit in the cone $\cK$ shown in 
Figure~\ref{fig:06}. 

In particular, the line $\delta(\varepsilon)=\varepsilon$ becomes a dividing line
around which we find an asymptotically exponentially small wedge. Outside this wedge, we have 
two completely different dynamical behaviors (I) and (II) as described above. Note that it also
makes sense to formally continue the classification of (I) and (II) onto the two lines 
$\{\varepsilon=0,\delta>0\}$ and $\{\varepsilon>0,\delta=0\}$ by using so-called candidate
trajectories obtained by concatenating orbits of the suitable fast and slow 
subsystem singular limit problems. Yet, we evidently cannot make a meaningful classification 
at the origin $(\varepsilon,\delta)=(0,0)$ itself regarding our property due to the undefined 
expression $\varepsilon^2/\delta$ in the last term of the fast variable dynamics in this case.\medskip

\begin{figure}[htbp]
	\centering
\begin{tikzpicture}
[>=stealth',point/.style={circle,inner sep=0.035cm,fill=white,draw},
encircle/.style={circle,inner sep=0.07cm,draw},
x=4cm,y=4cm,declare function={f(\x) = \x^2;}]

\draw[->,thick,dashed] (0,0) -> (1,0);
\draw[->,thick,dashed] (0,0) -> (0,1);

\draw[red,semithick,-,smooth,domain=0:0.65,samples=75,/pgf/fpu,/pgf/fpu/output
format=fixed] plot (\x, {2*f(\x)});

\draw[red,semithick,-,smooth,domain=0:0.9,samples=75,/pgf/fpu,/pgf/fpu/output
format=fixed] plot (\x, {0.7*f(\x)});

\node[point] at (0,0) {};

\node[] at (0.92,-0.05) {$\eps$};
\node[] at (-0.05,0.92) {$\delta$};

\draw [black] (0.2,0.72) circle [radius=0.06];
\node at (0.2,0.72) {III};
\draw [black] (0.7,0.55) circle [radius=0.06];
\node at (0.7,0.55) {II};
\draw [black] (0.7,0.15) circle [radius=0.06];
\node at (0.7,0.15) {I};

\end{tikzpicture}
\vspace{-4mm}	
	\caption{\label{fig:07}Classification diagram with respect to the property $\cP_{\textnormal{osc}}$
	for the problem~\eqref{eq:Olsen}. Note that the two parabolic thin curves (red) have the same
	function form $\delta(\varepsilon)=\kappa\varepsilon^2$ just with two different constants
	$\kappa>0$.
	The three regions correspond to the cases (I)--(III) above yielding non-standard relaxation
	oscillations, MMOs, and chaotic oscillations respectively.}
\end{figure}
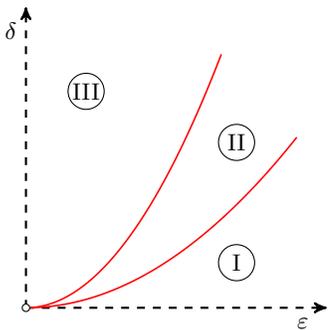

Obviously, the fast-slow normal form transcritical bifurcation case we have discussed here is 
just one of many cases in multiple time scale dynamics where several small parameters 
appear~\cite{KuehnBook}. Another important system directly motivated by a particular application to 
the peroxidase-oxidase reaction is the Olsen model~\cite{DegnOlsenPerram,Olsen}. It is given by
\be
\label{eq:Olsen} \tag{{$\cX_{\textnormal{Ol}}$}}
\begin{array}{rcl} 
\frac{\txtd a}{\txtd s}&=& \delta^2(p_1-\alpha a) -aby,\\
\frac{\txtd b}{\txtd s}&=& \varepsilon (\delta\varepsilon-\delta bx)- \delta aby,\\
\varepsilon \frac{\txtd x}{\txtd s}&=& -x^2 +\varepsilon (b-p_2) x +3aby +\varepsilon^2 p_4,\\
\frac{\txtd y}{\txtd s}&=&p_3(x^2-y-aby),
\end{array}
\ee
where $(a,b,x,y)\in(\R^+)^4$, we fix the parameters $p_1=0.97$, $p_2=0.98$, $p_3=3.93$, 
$p_4=1.2\cdot 10^{-5}$ to the classical values considered by Olsen~\cite{DegnOlsenPerram,Olsen} 
and take $\varepsilon,\delta$ as the small parameters. Then one can prove~\cite{KuehnSzmolyan} 
that for 
\benn
\varepsilon^2\ll \delta,
\eenn 
the system~\eqref{eq:Olsen} exhibits non-standard, but regularly periodic, relaxation oscillations~\cite{MisRoz}.
A singular limit geometric phase space description~\cite{KuehnSzmolyan}, as well as numerical
simulations~\cite{Olsen,DegnOlsenPerram} and numerical continuation 
calculations~\cite{DesrochesKrauskopfOsinga1,MusokeKrauskopfOsinga}, strongly suggest that there are at least two further important asymptotic regimes namely
\benn
\varepsilon^2\gg\delta \quad \text{ and }\quad \kappa\varepsilon^2=\delta=\delta(\varepsilon),~\kappa=\cO(1), 
\eenn
as $\varepsilon\ra 0$. In these cases one observes mixed-mode oscillations 
(MMOs)~\cite{Desrochesetal} and chaotic oscillations respectively, i.e., we have
for the Olsen model
\begin{itemize}
 \item[(I)] $\varepsilon^2\ll \delta$: non-standard relaxation oscillations,
 \item[(II)] $\cO(\varepsilon^2)=\delta$: chaotic oscillations,
 \item[(III)] $\varepsilon^2\gg\delta$: mixed-mode oscillations,
\end{itemize} 
which is illustrated in Figure~\ref{fig:07}. 

If we want to distinguish the three different classes of oscillation patterns (relaxation, MMO, 
chaos), then it does not suffice to rely on distinguishing properties individually such
as number of maxima for one variable within a time interval $\cP_{\textnormal{max}}$, the sign of the top/leading Lyapunov exponent $\cP_{\textnormal{Lya}}$, or topological equivalence of the phase portraits 
$\cP_{\textnormal{top}}$. For example, one expects that stable relaxation oscillations and MMOs 
may have topologically equivalent phase portraits and negative Lyapunov exponents for certain 
parameters. Furthermore, the number of maxima is also not a good indicator alone as for a given 
initial condition and a fixed time interval it is easily conceivable that an MMO and a chaotic 
oscillation have the same $\cP_{\textnormal{max}}$. Yet, suppose we fix a generic initial condition 
in the positive quadrant and a positive sufficiently large fixed time $T=K\varepsilon$ with a 
constant $K>0$ such that for $\cP_{\textnormal{max}}$ we have a fixed number 
$\cP_{\textnormal{max}}=K_0>0$ for the case of non-standard relaxation oscillations. Let us define
\benn
\cP_{\textnormal{osc}}:=\cP_{\textnormal{max}}\cP_{\textnormal{Lya}};
\eenn
then we expect that all three cases are different. Indeed, we conjecture that 
\begin{itemize}
 \item $\cP_{\textnormal{osc}}=-K_0$: stable non-standard relaxation oscillations; 
 \item $\cP_{\textnormal{osc}}<-K_0$: stable mixed-mode oscillations;
 \item $\cP_{\textnormal{osc}}>0$: chaotic oscillations.
\end{itemize}
Evidently, this is not a full classification, nor yet rigorously proven beyond the non-standard
relaxation case. However, it is very helpful to conceptually understand the Olsen
model and its analysis; see Figure~\ref{fig:07}. The difficulties of the problem are
now made precise and much more apparent. Already defining the property $\cP$ can be crucial to
make the mathematical analysis tractable as proving a precise shape of a trajectory as well as
an estimate of the Lyapunov exponent are highly non-trivial for global orbits of non-linear 
systems. Although methods from geometric singular perturbation theory exist to try to deal with
this situation~\cite{GuckenheimerWechselbergerYoung}, we expect that for the Olsen model these may 
have to be augmented by computer-assisted proof techniques~\cite{Haiduc1} to actually deal with 
tracking the dynamics in certain two- and three-dimensional reduced systems. As another question,
Figure~\ref{fig:07} points us immediately to the transition regimes, i.e., one should ask how 
trajectories are deformed near the separating asymptotic boundary curves and what happens near/on the non-negative cone $\partial\cK$ in $(\varepsilon,\delta)$-parameter space. Such a discussion is beyond the scope of this work.\medskip

For our examples so far, the second small parameter arose due to the need to study a bifurcation
problem, and the bifurcation parameter produced a double singular limit. Yet, in many applications,
there are additional ``physical'' modeling constraints, which lead to two small parameters.
A typical case is the effect of small noise, which is going to be discussed in the next subsection. 

\subsection{Stochastic Fast-Slow Systems}
\label{ssec:stochfss}

Among the most popular models for random noise acting on a dynamical system are 
stochastic differential equations (SDEs) driven by a Wiener process. There is a 
broad literature on such equations, based on different approaches such as 
analysing the Fokker--Planck equation \cite{HorsthemkeLefever}, the theory of 
large deviations \cite{FreidlinWentzell_book}, and random dynamical systems 
\cite{Arnold98}. Stochastic systems with multiple timescales have been more 
particularly analysed in 
\cite{BGbook,KabanovPergamenshchikov_2003,PavliotisStuart}. The stochastic 
dynamics near bifurcation points has been studied, for instance, in 
\cite{Stocks_Manella_McClintock_89,Swift_Hohenberg_Ahlers91,
Crauel_Flandoli_98,Jansons_Lythe98,Kuske99}. A particularly important field of 
application is neuroscience. In this respect, we refer to \cite{Tuckwell} for 
an overview, and to \cite{Lindner_Schimansky_1999,Lindneretal_04, 
KosmidisPakdaman_03,MuratovVandeneijnden2007,HitczenkoMedvedev_09,
Borowski_Kuske_etal_2011,Baxendale_Greenwood_11,Ditlevsen_Greenwood_12,
SimpsonKuske_2011} for examples of specific problems involving bifurcations.

Consider a stochastic differential equation of the form 
\begin{equation}
\label{eq:SDE} 
 \6x_t = f(x_t,\eps t)\6t + \sigma\6W_t\;,
\end{equation} 
where $f: \R^2\to\R$ is sufficiently smooth, and $W_t$ is a Wiener process 
describing white noise. The small parameters are $\eps$, which measures the 
slow drift of the \lq\lq parameter\rq\rq\ $y = \eps t$, and $\sigma$, which 
measures the noise intensity. 

In order to understand the influence of the noise on time scales, let us start 
by considering the case where $f = f(x)$ does not depend on the second 
variable, and let $V$ be a potential such that $f(x) = -V'(x)$. Assume that $V$ 
has a minimum at $x=0$. Then the theory of large 
deviations~\cite{FreidlinWentzell_book} implies that the probability of a 
solution of the SDE starting from $x_0=0$ to reach a point $x$ in a time of 
order $1$ is of order $\txte^{-V(x)/(2\sigma^2)}$, assuming $V$ is monotonous 
between $0$ and $x$.  This implies the so-called Arrhenius 
law~\cite{Arrhenius}, which states that the expected time for the solution to 
reach $x$ has order $\txte^{V(x)/(2\sigma^2)}$. Solutions of the SDE thus tend 
to spend exponentially long time spans near stable stationary points of $f$. 

When considering the slowly time-dependent system~\eqref{eq:SDE}, it is 
convenient to scale time by a factor $\eps$, so that $f$ changes by order $1$ 
in times of order $1$. The rescaled system reads 
\begin{equation}
\label{eq:SDE_slow} \tag{{$\cX_{\textnormal{sfs}}$}}
 \6x_t = \frac{1}{\eps} f(x_t,t)\6t 
 + \frac{\sigma}{\sqrt{\eps}}\6W_t\;,
\end{equation} 
where the factor $\sqrt{\eps}$ is due to the scaling property of the Wiener 
process. \textcolor{black}{We remark that from \eqref{eq:SDE_slow} it is clear that the problem is singularly perturbed as it is fast-slow in $\eps$ and degenerates from an SDE to an ODE for $\sigma=0$.} \textcolor{black}{This is also apparent in the infinitesimal generator of the SDE \eqref{eq:SDE_slow}, which is given by $\frac{\sigma^2}{2\eps}\Delta + \frac{1}{\eps}f\cdot\nabla$. The evolution of the probability density of the SDE, as well as its exit distribution from a domain, are thus described by parabolic or elliptic PDEs having a small parameter multiplying the highest derivative.}

Assume $f$ has a smooth stable equilibrium branch $x^*(t)$ acting as a 
critical manifold for~\eqref{eq:SDE_slow}. This means that $f(x^*(t),t) = 0$ 
for all $t$ in some interval $I$, and that $a^*(t) = \partial_x f(x^*(t),t)$ is 
negative, bounded away from $0$ in $I$. In the deterministic case $\sigma=0$, it 
is well known~\cite{Tihonov,Fenichel4} that  for small $\eps$, 
\eqref{eq:SDE_slow} admits a so-called slow solution $\bar x(t)$ staying 
at a distance of order $\eps$ from $x^*(t)$.  

\begin{figure}
\begin{center}
\begin{tikzpicture}
[>=stealth',point/.style={circle,inner sep=0.035cm,fill=white,draw},
encircle/.style={circle,inner sep=0.07cm,draw},
x=4cm,y=4cm,declare function={f(\x) = 3*exp(-1/\x^2);}]

\draw[->,semithick] (0,0) -> (1,0);
\draw[->,semithick] (0,0) -> (0,1);

\node[point] at (0,0) {};

\draw[blue,thick,-,smooth,domain=0.01:0.9,samples=75,/pgf/fpu,/pgf/fpu/output
format=fixed] plot ({f(\x)}, \x);

\node[] at (0.9,-0.05) {$\eps$};
\node[] at (-0.05,0.9) {$\sigma$};

\node[] at (0.5,0.1) {$P(\sigma,\eps) \simeq 0$};
\node[] at (0.23,0.87) {$P(\sigma,\eps) \simeq 1$};

\draw [black] (0.2,0.72) circle [radius=0.06];
\node at (0.2,0.72) {II};
\draw [black] (0.5,0.3) circle [radius=0.06];
\node at (0.5,0.3) {I};

\end{tikzpicture}
\vspace{-4mm}
\end{center}
\caption{Probability $P(\sigma,\eps)$ that a solution of the 
SDE~\eqref{eq:SDE_slow} leaves the neighborhood of a stable critical manifold 
in slow time of order $1$, in the parameter space $(\eps,\sigma)$. The 
probability is close to $0$ or $1$, except near the curve 
$\eps = \exp[-H/(2\sigma^2)]$.}
\label{fig:NB_sigma_epsilon}
\end{figure}
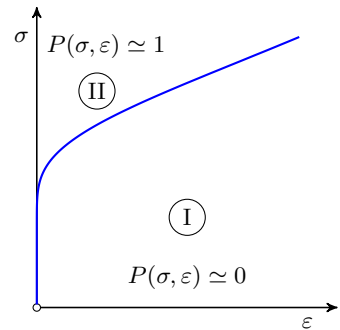

Let us now fix, say, $I=[0,1]$, and consider the solution 
of~\eqref{eq:SDE_slow} starting at time $0$ in $\bar x(0)$. Denote by 
$P(\sigma,\eps)$ the probability that the solution leaves a neighborhood of 
$\bar x(t)$ at or before time $1$. Then 
\begin{itemize}
\item[(I)] 	on one hand, the large-deviation results just mentioned imply 
that when $\sigma$ decreases to $0$ for fixed $\eps>0$, $P(\sigma,\eps)$ 
converges to $0$; 
\item[(II)] 	on the other hand, irreducibility of the Markov process 
$(x_t)_{t\ge0}$ implies that when $\eps$ decreases to $0$ for fixed 
$\sigma>0$, $P(\sigma,\eps)$ converges to $1$. 
\end{itemize}
Hence, the regimes (I)-(II) induce a property $\cP_{\textnormal{sfs}}$,
which divides the $(\varepsilon,\sigma)$-space via the escape probability.
The transition between $P(\sigma,\eps)$ close to $0$ and close to $1$ occurs 
when $\eps$ is of order $\txte^{-H/(2\sigma^2)}$ for an $H>0$ depending on the 
considered neighborhood (Figure~\ref{fig:NB_sigma_epsilon}). 

A more precise formulation of the regime $\sigma\searrow0$ has been given 
in~\cite{BG1,BGbook}. Let $\cB(h)$ be a family of strips centered in 
$x^* (t)$, of width $h/\sqrt{2|a(t)|}$, where $a(t) = \partial_x f(\bar 
x(t),t)$ is the linearization of $f$ at the slow solution. These strips act as 
a kind of \lq\lq confidence intervals\rq\rq, in the sense that the probability 
$P_t(h,\sigma,\eps)$ of leaving $\cB(h)$ before time $t$ satisfies  
\begin{equation}
 P_t(h,\sigma,\eps) \simeq 
 \biggl[\frac{1}{\eps} \int_0^t |a(s)|\6s\biggr] 
 \,\frac{h}{\sigma} \txte^{-h^2/(2\sigma^2)}
\end{equation} 
as long as $t \ll \eps\txte^{ch^2/\sigma^2}$ for some constant $c>0$
(see~\cite[Theorem~3.1.10]{BGbook} for a precise formulation). Choosing 
$h$ of order $\sigma\sqrt{2\log(t/(\eps p))}$, we obtain $P_t(h,\sigma,\eps) 
\simeq p$, so that $\cB(h)$ is indeed a confidence strip at level $p$.  

This first example of a two-scale behavior for an SDE is somewhat atypical 
compared to other examples given in this review, in the sense that the 
transition between qualitatively different regimes occurs when $\eps$ is 
exponentially small in $\sigma$. Of course, one can \lq\lq regularize\rq\rq\ 
things by writing $\eps = \txte^{-\lambda/\sigma^2}$ and describing the 
behavior in terms of $\lambda$ and $\sigma$. This is the approach adopted 
in~\cite{Freidlin1} for instance.  

\begin{figure}
\begin{center}
\begin{tikzpicture}
[>=stealth',point/.style={circle,inner sep=0.035cm,fill=white,draw},
encircle/.style={circle,inner sep=0.07cm,draw},
x=2.5cm,y=2.5cm,declare function={f(\x) = sqrt((\x-0.1)^2+0.15)-0.1;
g(\x) = 0.2 -\x - sqrt(\x^2 + 0.1);}]

\draw[->,semithick] (-1,0) -> (1.2,0);
\draw[->,semithick] (0,-1) -> (0,1.2);

\draw[thick,gray] (-1,1) -- (0,0) -- (1,1);
\draw[thick,gray,dashed] (-1,-1) -- (0,0) -- (1,-1);

\draw[semithick,blue] plot[smooth,tension=.6]
  coordinates{(-1,1.05) (-0.5,0.6) (0,0.35) (0.5,0.4) (1,0.95)};

\draw[red,semithick,-,smooth,domain=-0.9:1,samples=75,/pgf/fpu,
/pgf/fpu/output format=fixed] plot ({\x}, {f(\x) +0.03*rand});

\draw[red,semithick,-,smooth,domain=-0.95:0.9,samples=75,/pgf/fpu,
/pgf/fpu/output format=fixed] plot ({0.5*g(\x)}, {\x +0.1*rand});

\node[] at (1,-0.1) {$t$};
\node[] at (-0.1,1) {$x$};

\node[gray] at (0.65,0.93) {$x^*_+(t)$};
\node[gray] at (0.65,-0.93) {$x^*_-(t)$};

\node[blue] at (0.9,0.6) {$\bar x(t)$};
\node[red] at (1.05,0.8) {$x_t$};
\node[red] at (0.17,-0.8) {$x_t$};

\draw[semithick] (0.4,0.04) -- (0.4,-0.04);
\draw[semithick] (-0.4,0.04) -- (-0.4,-0.04);
\draw[semithick] (-0.04,0.4) -- (0.04,0.4);

\node[] at (0.4,-0.13) {$\sqrt{\eps}$};
\node[] at (-0.5,-0.13) {$-\sqrt{\eps}$};
\node[] at (0.15,0.4) {$\sqrt{\eps}$};

\end{tikzpicture}
\vspace{-4mm}
\end{center}
\caption{Slow passage through a transcritical bifurcation. The blue curve is a 
deterministic solution of~\eqref{eq:SDE_transcritical} with $\sigma=0$, which 
stays at distance at most of order $\sqrt{\eps}$ from the stable critical 
curve $x^*_+(t) = |t|$. Red paths sketch the behavior of typical stochastic 
solutions $x_t$, in parameter regimes (I) (upper path) and (II) (lower path).}
\label{fig:NB_transcritical}
\end{figure}
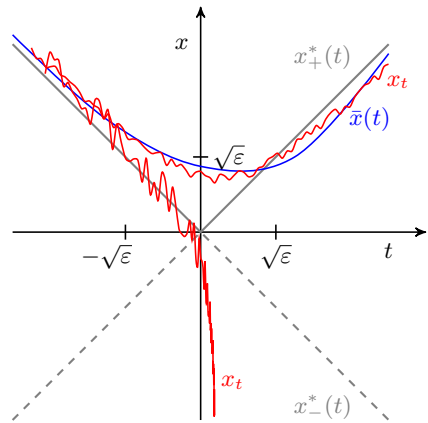

More standard examples of double limits can however be found in the vicinity of 
bifurcation points. Consider for instance the fast-slow 
SDE~\eqref{eq:SDE_slow} 
\begin{equation}
\label{eq:SDE_transcritical} \tag{{$\cX_{\textnormal{tcs}}$}}
 \6x_t = \frac{1}{\eps} (t^2 - x_t^2)\6t 
 + \frac{\sigma}{\sqrt{\eps}}\6W_t\;,
\end{equation} 
which is a stochastic version of~\eqref{eq:tc}. The critical manifold of the 
deterministic equation $\eps\dot x = t^2 - x^2$ is composed of a stable branch 
$\{x = x^*_+(t) = |t| \colon t\neq 0\}$ and an unstable branch $\{x = x^*_-(t) 
= -|t| \colon t\neq 0\}$. It is well-known (see for instance~\cite{Haberman}) 
that when $\sigma=0$, the equation~\eqref{eq:SDE_transcritical} admits a slow 
solution $\bar x(t)$ of order $\max\{|t|, \sqrt{\eps}\}$. This solution tracks 
the stable branch $x^*_+(t)$ at a distance of order 
$\eps/\max\{|t|,\sqrt{\eps}\}$ (Figure~\ref{fig:NB_transcritical}). 

In the case $\sigma>0$, we can define as above a strip $\cB(h)$ centered in the 
slow solution $\bar x(t)$, and of width $h/\sqrt{2|a(t)|}$. Note that this 
time, the linearization $|a(t)|$ has order  $\max\{|t|, \sqrt{\eps}\}$. The 
width of $\cB(h)$ is maximal near $t=0$, where it has order $h/\eps^{1/4}$. It 
turns out that one then has two qualitatively different situations~\cite{BG2}: 

\begin{itemize}
\item[(I)] 	If $\sigma \ll \eps^{3/4}$, we can take $h$ of order 
$\eps^{3/4}$ and still have a strip $\cB(h)$ staying away from the origin. One 
can then show that the probability of a solution of~\eqref{eq:SDE_transcritical} 
leaving $\cB(h)$ before, say, time $1$, has order $\exp[-h^2/(2\sigma^2)] = 
\exp[-\eps^{3/2}/\sigma^2]$, which is exponentially small in this regime. 

\item[(II)] 	If $\sigma \gg \eps^{3/4}$, on the other hand, any strip 
$\cB(h)$ with $h\ge\sigma$ intersects the $t$-axis already at or before a time 
of order $-\sigma^{2/3}$. One can then show that it is very likely that the 
solution $x_t$ becomes negative, of order $1$, shortly after time 
$-\sigma^{2/3}$. The probability that $x_t$ remains positive up to time $1$ has 
order $\txte^{-\sigma^{4/3}/(\eps\log(\sigma^{-1}))}$. 
\end{itemize}

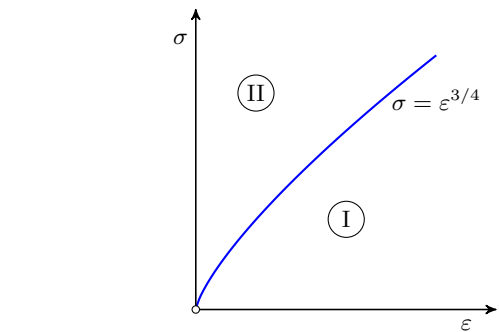
\begin{figure}
\begin{center}
\begin{tikzpicture}
[>=stealth',point/.style={circle,inner sep=0.035cm,fill=white,draw},
encircle/.style={circle,inner sep=0.07cm,draw},
x=4cm,y=4cm,declare function={f(\x) = \x^(3/4);}]
\draw[blue,thick,-,smooth,domain=0.001:0.8,samples=75,/pgf/fpu,/pgf/fpu/output
format=fixed] plot (\x, {f(\x)});

\draw[->,semithick] (0,0) -> (1,0);
\draw[->,semithick] (0,0) -> (0,1);

\node[point] at (0,0) {};

\node[] at (0.9,-0.05) {$\eps$};
\node[] at (-0.05,0.9) {$\sigma$};

\node[] at (0.8,0.7) {$\sigma = \eps^{3/4}$};

\draw [black] (0.5,0.3) circle [radius=0.06];
\node at (0.5,0.3) {I};
\draw [black] (0.2,0.72) circle [radius=0.06];
\node at (0.2,0.72) {II};

\end{tikzpicture}
\vspace{-4mm}
\end{center}
\caption{The probability $\Ptrans$ that the solutions $x_t$ of 
equation~\eqref{eq:SDE_transcritical} starting on the stable slow solution 
$\bar x(t)$ becomes negative behaves differently in the two shown parameter 
regions. In Region (I), $\Ptrans$ has order $\exp[-\eps^{3/2}/\sigma^2]$, while 
in Region (II), $1-\Ptrans$ has order 
$\exp[-\sigma^{4/3}/(\eps\log(\sigma^{-1}))]$.}
\label{fig:NB_transcritical_parameters}
\end{figure}

One can summarize the difference between the two regimes by considering the 
transition probability 
\begin{equation}
 \Ptrans(\sigma,\eps) = \P^{(\bar x(t_0),t_0)} \bigl\{ \exists t \le 1 \colon 
x_t = -1 \bigr\}\;,
\end{equation} 
where the superscript $(\bar x(t_0),t_0)$ indicates the initial condition. For 
negative $t_0$ of order $1$, we have 
\begin{equation}
 \Ptrans(\sigma,\eps) 
 \begin{cases}
  \le \txte^{-\kappa\eps^{3/2}/\sigma^2} & \text{in Regime (I)}\;, \\
  \ge 1 - \txte^{-\kappa\sigma^{4/3}/(\eps\log(\sigma^{-1}))} & \text{in Regime 
(II)}\;,
 \end{cases}
\end{equation} 
for a constant $\kappa>0$  (Figure~\ref{fig:NB_transcritical_parameters}). Hence,
we can again use a suitable transition probability to define a property 
$\cP_{\textnormal{tcs}}$, which provides at least two clearly distinct asymptotic
regimes (I)-(II) in the double limit. See~\cite[Theorems 3.5.1 and 3.5.2]{BGbook} for 
precise formulations of these results. 

An interesting generalization of Example~\eqref{eq:SDE_transcritical} is the 
SDE 
\begin{equation}
\label{eq:SDE_avoided_transcritical} \tag{{$\cX_{\textnormal{tcd}}$}}
 \6x_t = \frac{1}{\eps} (t^2 - x_t^2 + \delta)\6t 
 + \frac{\sigma}{\sqrt{\eps}}\6W_t\;, 
\end{equation} 
where the parameter $\delta > 0$ plays the same role as $\eps/\delta$ 
in~\eqref{eq:tc}. Note that we are now dealing with three small parameters 
$\eps$, $\sigma$, and $\delta$. The critical manifolds are given here by 
$x^*_\pm(t) = \pm\sqrt{t^2+\delta}$, so that they do not quite touch: their 
minimal distance is $2\sqrt{\delta}$. 

A similar analysis as for the transcritical 
bifurcation~\eqref{eq:SDE_transcritical} can be made, and results in the 
following case distinction 
(Figure~\ref{fig:NB_avoided_transcritical_parameters}, see~\cite[Theorems 2.6 
and 2.7]{BG2} for precise formulations):

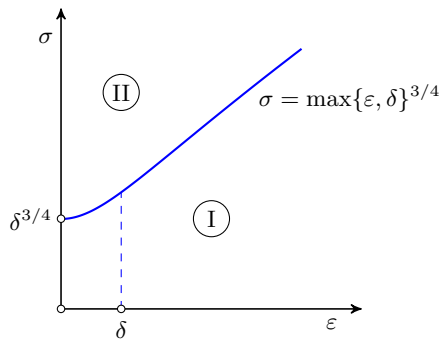
\begin{figure}[h]
\begin{center}
\begin{tikzpicture}
[>=stealth',point/.style={circle,inner sep=0.035cm,fill=white,draw},
encircle/.style={circle,inner sep=0.07cm,draw},
x=4cm,y=4cm,declare function={f(\x,\delta) = (\x^2 + \delta^2)^(3/8);}]

\pgfmathsetmacro{\mydelta}{0.2}

\draw[->,semithick] (0,0) -> (1,0);
\draw[->,semithick] (0,0) -> (0,1);

\node[point] at (0,0) {};

\draw[blue,thick,-,smooth,domain=0.001:0.8,samples=75,/pgf/fpu,/pgf/fpu/output
format=fixed] plot (\x, {f(\x,\mydelta)});

\node[] at (0.9,-0.05) {$\eps$};
\node[] at (-0.05,0.9) {$\sigma$};

\node[] at (0.96,0.7) {$\sigma = \max\{\eps,\delta\}^{3/4}$};

\draw [black] (0.5,0.3) circle [radius=0.06];
\node at (0.5,0.3) {I};
\draw [black] (0.2,0.72) circle [radius=0.06];
\node at (0.2,0.72) {II};

\node[point] at (0, {f(0,\mydelta)}) {};
\node[] at (-0.1, {f(0,\mydelta)}) {$\delta^{3/4}$};

\draw[blue,dashed] (\mydelta,0) -- (\mydelta, {f(\mydelta,\mydelta)});

\node[point] at (\mydelta,0) {};
\node[] at (\mydelta, -0.07) {$\delta$};
\end{tikzpicture}
\vspace{-4mm}
\end{center}
\caption{$(\eps,\sigma)$-parameter plane for the 
SDE~\eqref{eq:SDE_avoided_transcritical} describing an avoided transcritical 
bifurcation, for a fixed $\delta>0$.  In Region (I), the transition probability 
$\Ptrans$ has order $\exp[-\max\{\eps,\delta\}^{3/2}/\sigma^2]$, while in Region 
(II), $1-\Ptrans$ has order $\exp[-\sigma^{4/3}/(\eps\log(\sigma^{-1}))]$.}
\label{fig:NB_avoided_transcritical_parameters}
\end{figure}

\begin{itemize}
\item[(I)] 	If $\sigma \ll \max\{\eps,\delta\}^{3/4}$, solutions tend to 
stay near the slow solution $\bar x(t)$ tracking $x^*_+(t)$, and the transition 
probability $\Ptrans$ is exponentially small. 

\item[(II)] 	If $\sigma \gg \max\{\eps,\delta\}^{3/4}$, solutions are likely 
to escape to negative values of $x$ as soon as $t$ is slightly larger than 
$-\sigma^{2/3}$. 
\end{itemize}

This results in a transition probability behaving as 
\begin{equation}
 \Ptrans(\sigma,\eps) 
 \begin{cases}
  \le \txte^{-\kappa\max\{\eps,\delta\}^{3/2}/\sigma^2} & \text{in Regime 
(I)}\;, \\
  \ge 1 - \txte^{-\kappa\sigma^{4/3}/(\eps\log(\sigma^{-1}))} & \text{in Regime 
(II)}\;. 
 \end{cases}
\end{equation} 
The parameter $\delta$ thus causes a saturation effect at small values of 
$\eps$. 

The examples considered so far were all particular cases of the slowly 
time-dependent SDE~\eqref{eq:SDE_slow}. Other types of bifurcations, such as 
the saddle-node bifurcation, which results in similar regimes with different 
exponents, are described in~\cite[Chapter~3]{BGbook}. One can however also 
consider fully coupled fast-slow systems of the form 
\begin{align}
\nonumber
 \6x_t &= \frac{1}{\eps} f(x_t,y_t)\6t + \frac{\sigma}{\sqrt{\eps}} 
F(x_t,y_t)\6W_t\;, \\
 \6y_t &= g(x_t,y_t)\6t + \sigma'G(x_t,y_t) \6W_t\;,
\end{align}
where $x\in\R^m$, $y\in\R^n$, and $W_t$ is a $k$-dimensional Wiener process.
In a similar way as for~\eqref{eq:SDE}, one can obtain concentration results 
for solutions near stable normally hyperbolic critical manifolds, 
see~\cite{BG6}.

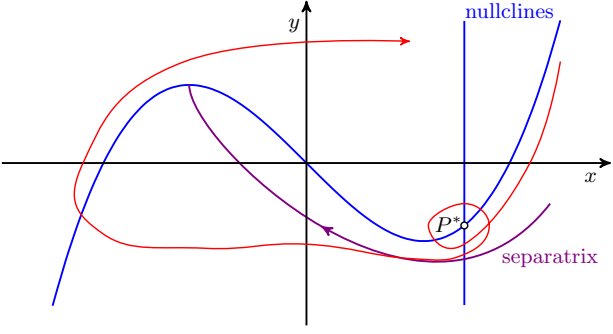
\begin{figure}
\begin{center}
\scalebox{0.9}{
\begin{tikzpicture}[>=stealth',main node/.style={circle,minimum
size=0.1cm,inner sep=0pt,fill=white,draw},x=3cm,y=3cm,
declare function={
yy(\x) = \x^3 - \x;}]

\pgfmathsetmacro{\xfhn}{1/sqrt(3)}
\pgfmathsetmacro{\afhn}{\xfhn +0.2}
\pgfmathsetmacro{\yfhn}{\xfhn *2/3}
\pgfmathsetmacro{\xminus}{-\xfhn}
\pgfmathsetmacro{\bfhn}{\afhn^3 - \afhn}




\draw[blue,thick,-,smooth,domain=-1.25:1.25,samples=40,/pgf/fpu,/pgf/fpu/output
format=fixed] plot (\x, {yy(\x)});
\draw[blue,thick] (\afhn,-0.7) -- (\afhn,0.7);

\draw[violet,thick] (\xminus,\yfhn) .. controls (\xminus,0.1) and
(\xfhn,-1) .. (1.2,-0.2)
\arrowonlineback{0.5};


\draw[->,thick] (-1.5,0) -- (1.5,0);
\draw[->,thick] (0,-0.8) -- (0,0.8);

\draw[red,semithick] plot[smooth,tension=0.8]
  coordinates{(1.25,0.5) (1.1,0) (\afhn,-0.4) (0.6,\bfhn) (\afhn,-0.2)
(0.9,\bfhn) (\afhn,-0.45) (0.5,-0.47) (0,-0.4) (-0.5,-0.42) (-1,-0.35) (-1.1,0)
(-0.6,0.5) (0.5,0.6)};
\draw[red,semithick,->] (0.5,0.6) -- (0.51,0.6);



\node[main node,semithick] at (\afhn,\bfhn) {};

\node[] at (0.7,-0.3) {{\small $P^*$}};

\node[] at (1.4,-0.07) {$x$};
\node[] at (-0.06,0.68) {$y$};

\node[violet] at (1.2,-0.47) {separatrix};
\node[blue] at (1,0.75) {nullclines};

\end{tikzpicture}
}
\end{center}
\vspace{-8mm}
\caption{Phase space of the stochastic FitzHugh--Nagumo system~\eqref{eq:FHN}. 
The separatrix is defined as the deterministic negative orbit of the  
local maximum of the $x$-nullcline $\{y=x^3-x\}$. When $P^*$ is a focus, 
stochastic solutions tend to perform small oscillations around $P^*$ before 
crossing the separatix, and making a large excursion (or spike) before 
returning near $P^*$.}
\label{fig:NB_FHN_phase_space}
\end{figure}

A particularly interesting case is the stochastic FitzHugh--Nagumo system 
modelling action potential dynamics of individual neurons, investigated 
in~\cite{MuratovVandeneijnden2007,BerglundLandon}. We consider here the 
particular case 
\begin{align}
\nonumber 
 \6x_t &= \frac{1}{\eps} \bigl[x_t - x_t^3 + y_t\bigr]\6t + 
\frac{\sigma}{\sqrt{\eps}} 
\6W_t^{(1)}\;, \\
 \6y_t &= \bigl[a - x_t\bigr] \6t + \sigma \6W_t^{(2)}\;, 
\label{eq:FHN} \tag{{$\cX_{\textnormal{FHs}}$}}
\end{align}
where $W_t^{(1)}$ and $W_t^{(2)}$ are independent Wiener processes. In the 
deterministic case $\sigma=0$, the system~\eqref{eq:FHN} has a unique 
equilibrium point $P^* = (a,a^3-a)$. The eigenvalues of the linearisation 
at $P^*$ are given by 
\begin{equation}
 \lambda_\pm = \frac{-\delta\pm\sqrt{\delta^2-\eps}}{\eps}\;, 
 \qquad 
 \delta = \frac{3a^2-1}{2}\;.
\end{equation} 
Hence $P^*$ is a stable node for $\delta > \sqrt{\eps}$, 
a stable focus for $0 < \delta < \sqrt{\eps}$, 
an unstable focus for $-\sqrt{\eps} < \delta < 0$, 
and an unstable node for $\delta < -\sqrt{\eps}$. 

We are interested here in the excitable regime $0 < \delta \ll 1$, 
$\sigma > 0$. In that situation, though $P^*$ is stable in the deterministic 
case, it lies close to a (pseudo-)separatrix 
(Figure~\ref{fig:NB_FHN_phase_space}). Whenever the noise kicks it over the 
separatrix, the system makes a large excursion before returning to its rest 
state, producing a so-called \emph{spike} of the neuron's membrane potential. 

\begin{figure}
\begin{center}
\begin{tikzpicture}[>=stealth',main node/.style={circle,minimum
size=0.08cm,inner sep=0pt,fill=white,draw},encircle/.style={circle,inner 
sep=0.07cm,draw},
x=16cm,y=16cm]

\pgfmathsetmacro{\epsi}{0.1}
\pgfmathsetmacro{\xmax}{0.4}
\pgfmathsetmacro{\ymax}{0.26}
\pgfmathsetmacro{\xx}{sqrt(\epsi)}
\pgfmathsetmacro{\yy}{sqrt(\epsi *\xx)}





\draw[thick,blue] (0,0) -- (\xx,\yy);
\draw[thick,blue,smooth,domain=\xx:\xmax,samples=40,/pgf/fpu,
/pgf/fpu/output format=fixed] plot (\x, {\x^(1.5)});
\draw[thick,blue,dashed,smooth,domain=0:\xx,samples=60,/pgf/fpu,
/pgf/fpu/output format=fixed] plot (\x, {sqrt(\epsi *\x)});
\draw[thick,blue] (0,\yy) -- (\xx,\yy);
\draw[thick,dashed,blue] (\xx,\yy) -- (\xx,0);

\draw[->,thick] (0,0) -- (0.43,0);
\draw[->,thick] (0,0) -- (0,0.29);


\node[main node] at (\xx,\yy) {};

\draw[thick] (-0.007,\yy) -- node[left=0.08cm] {{$\eps^{3/4}$}}
(0.007,\yy);

\draw[thick] (\xx,-0.007) -- node[below=0.07cm] {{$\sqrt{\eps}$}}
(\xx,0.007);

\node[] at (0.41,-0.02) {$\delta$};
\node[] at (-0.02,0.27) {$\sigma$};
 
\node[rotate=28] at (0.18,0.088) {$\sigma=\delta\eps^{1/4}$};
\node[rotate=27] at (0.11,0.12) {$\sigma=\sqrt{\delta\eps}$};
\node[rotate=42] at (0.35,0.225) {$\sigma=\delta^{3/2}$};

\draw [black,fill=white] ({\xx},0.05) circle [radius=0.018];
\node at ({\xx},0.05) {I};
\draw [black,fill=white] (0.06,0.08) circle [radius=0.018];
\node at (0.06,0.08) {II};
\draw [black,fill=white] (0.12,0.22) circle [radius=0.018];
\node at (0.12,0.22) {III};

\end{tikzpicture}
\vspace{-6mm}
\end{center}
\caption{$(\delta,\sigma)$-parameter plane for the stochastic FitzHugh--Nagumo  
SDE~\eqref{eq:FHN}, for a fixed $\eps>0$. The three regions correspond to (I) 
rare isolated spikes, (II) clusters of spikes, and (III) repeated spikes.}
\label{fig:NB_FHN_parameters}
\end{figure}
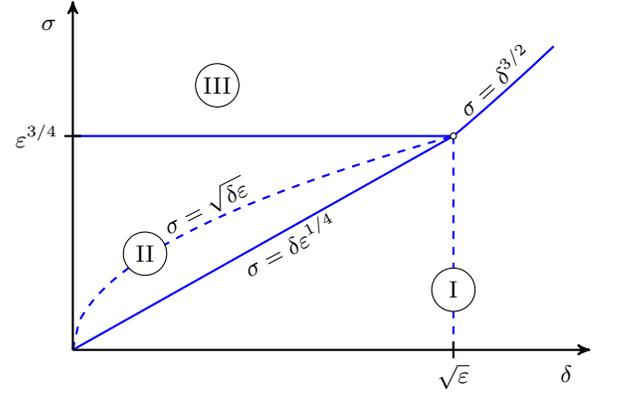

\begin{figure}
\begin{center}
\hspace{4mm}
\begin{overpic}[width=3.5cm,tics=10]{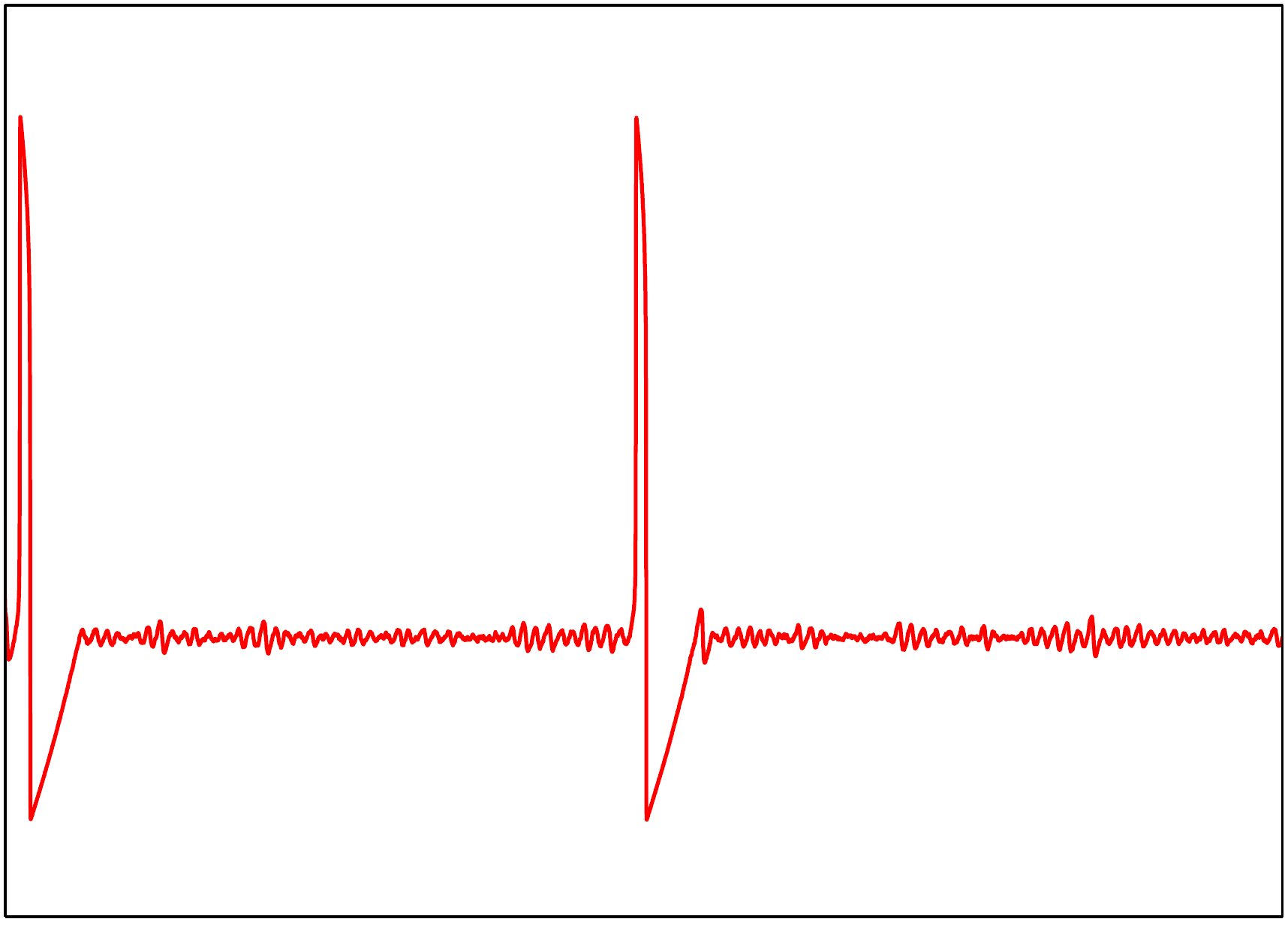}
\put(-12,62){(I)}
\end{overpic} 
\hspace{4mm}
\begin{overpic}[width=3.5cm,tics=10]{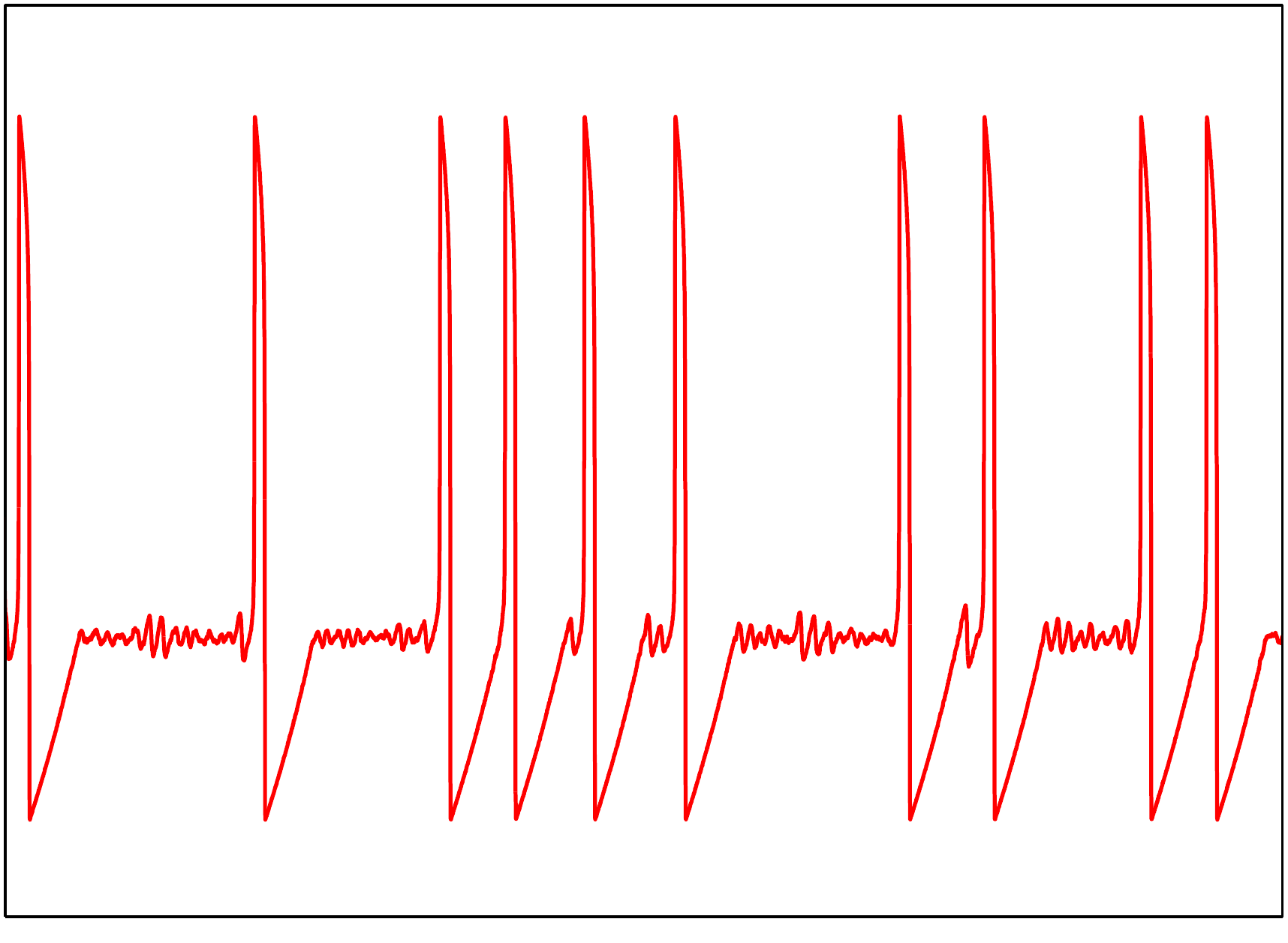}
\put(-15,62){(II)}
\end{overpic} \\
\vspace{1mm}\hspace{4mm}
\begin{overpic}[width=3.5cm,tics=10]{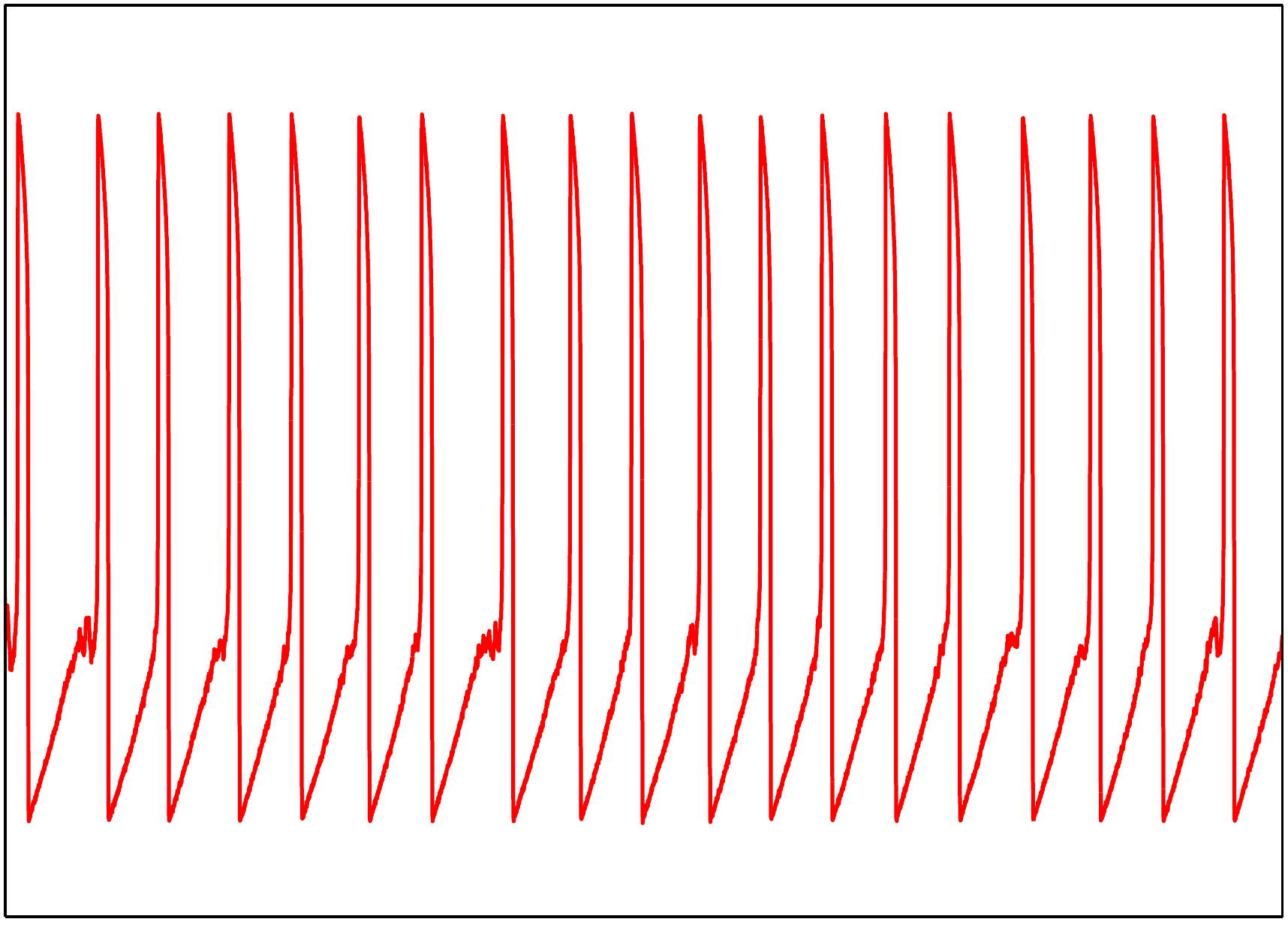}
\put(-19,62){(III)}
\end{overpic} 
\end{center}
\vspace{-6mm}
\caption{Time series $-x_t$ of solutions to the stochastic FitzHugh--Nagumo 
equation~\eqref{eq:FHN} in three different parameter regimes. Parameter values 
are $\eps=0.01$, $\delta=0.03$, and (I) $\sigma=0.001$, (II) $\sigma=0.0025$ 
and (III) $\sigma=0.01$.}
\label{fig:NB_FHN_timeseries}
\end{figure}

In~\cite{MuratovVandeneijnden2007}, the authors investigated the stochastic 
system~\eqref{eq:FHN} via formal computations, and found a large number of 
different parameter regimes. Some of these formal results have been proved 
rigorously in~\cite{BerglundLandon}. One can identify three main parameter 
regimes, as shown in Figure~\ref{fig:NB_FHN_parameters}:

\begin{itemize}
\item[(I)] 	If $0 < \delta < \sqrt{\eps}$ and $\sigma \ll \delta\eps^{1/4}$ 
or if $\sqrt{\eps} \leq \delta \ll 1$ and $\sigma \ll \delta^{3/2}$, the system 
displays rare isolated spikes (Figure~\ref{fig:NB_FHN_timeseries}-(I)). In 
particular, \cite[Theorem 3.2]{BerglundLandon} shows that if 
$\delta/\sqrt{\eps}$ is sufficiently small, then the expected number of small 
oscillations around $P^*$ between two consecutive spikes has order 
$\exp\{\delta^2\sqrt{\eps}/\sigma^2\}$. 

\item[(II)] 	If $0 < \delta < \sqrt{\eps}$ and $\delta\eps^{1/4} \le \sigma 
\le \eps^{3/4}$, one can observe clusters of spikes  
(Figure~\ref{fig:NB_FHN_timeseries}-(II)). In fact, what happens is that as 
$\sigma$ increases, the probability that a spike is immediately followed by 
another spike gradually increases like 
\begin{equation}
 \Phi\biggl(-\frac{\eps^{1/4}(\delta-\sigma^2/\eps)}{\sigma}\biggr)\;,
\end{equation} 
where $\Phi$ denotes the distribution function of a standard normal random 
variable. The dashed curve $\sigma=\sqrt{\delta\eps}$ in 
Figure~\ref{fig:NB_FHN_parameters} corresponds to this 
probability being close to $1/2$ (see~\cite[Section 5]{BerglundLandon}). 

\item[(III)] 	If $0 < \delta < \sqrt{\eps}$ and $\sigma \gg \eps^{3/4}$ 
or if $\sqrt{\eps} \leq \delta \ll 1$ and $\sigma \gg \delta^{3/2}$, the system 
displays repeated spikes (Figure~\ref{fig:NB_FHN_timeseries}-(III)), meaning 
that after having spiked, it is very likely to spike again immediately. 
\end{itemize}

Note that these three regimes actually use a probabilistic asymptotic
spiking pattern to define a property $\cP_{\textnormal{FHs}}$ to dissect the
(triple) singular limit parameter space. So the example nicely illustrates 
that also on a stochastic level, one can use macroscopic patterns, and that 
quite frequently even more than two small parameters are relevant.  

The behavior in regimes just described in (I)--(III) can be considered as a stochastic 
instance of mixed-mode oscillations (MMOs)~\cite{Desrochesetal}, in which 
small-amplitude and large-amplitude oscillations alternate; cf.~problem \eqref{eq:Olsen}. 
While deterministic MMOs often show a regular pattern, and sometimes a chaotic pattern, 
in the stochastic case considered here the number of small and large-amplitude 
oscillations are random variables. More intricate patterns can arise near 
folded-node bifurcations in three-dimensional SDEs, as for instance in the 
Koper model~\cite{BGK12,Berglund_Gentz_Kuehn_2015}.\medskip

The examples in this subsection have shown that the interplay between a
deterministic multiple time scale system with small noise provides a very
natural class of systems, and small noise induces a doubly singularly 
perturbed problem. Yet, stochastic differential equations provide many
other avenues to double limits, even without explicit time scale
separation for the drift. This is illustrated by the next subsection.

\subsection{Shear-Induced Chaos}
\label{ssec:shear}

In this section, we consider the interaction of shear forces and stochastic 
noise that can generate a switching from synchronization to chaotic behavior 
in stochastic oscillators. The onset of chaos by an interplay of shear and, typically small, noise has been broadly discussed within the context of stochastic Hopf bifurcation \cite{ArnoldSchenk96,Baxendale94,Baxendale03,Baxendale04,DoanEngelLambRasmussen, LinYoung08, Schenk96}, with important connections to coupled (neural) oscillators \cite{Blackbeardetal2011, LinSheaYoung09, LinYoung10, WedgwoodLin13} and questions around the role of noise and chaos in (turbulent) fluid flows \cite{Arnaudonetal18, HughesProctor, HughesProctor2, Farandaetal17}. Note that the idea of adding small noise to prove chaotic properties in the deterministic zero-noise limit has become an important tool in dynamical systems theory in recent years \cite{Blumenthaletal17, Blumenthaletal18, EngelGkogkasKuehn21, Young08}. 

As a basic toy model (cf.~\cite{EngelLambRasmussen1}), 
we consider the SDE, written in Stratonovich form,
\begin{align} 
\label{eq:ME_cylindermodel} \tag{{$\cX_{\textnormal{sic}}$}}
  \begin{array}{r@{\;\,=\;\,}l}
    \rmd y & - \alpha y \,\rmd t + \sigma \sum_{i=1}^m  f_i(\vartheta) 
		\circ \rmd W_t^i\,,\\
    \rmd \vartheta & (1 + b y) \,\rmd t
 \,,
  \end{array}
\end{align}
where $(y, \vartheta)\in\mathbb{R}\times \mathbb S^1$ are cylindrical amplitude-phase 
coordinates, $m \geq 1$ is a natural number, and $W_t^i$ for $i\in\{ 1, \dots, m\}$ 
denote independent one-dimensional Brownian motions. We will assume that 
$\alpha, \sigma, b \geq 0$, i.e.~all parameters are non-negative.

When there is no noise ($\sigma=0$), the \textcolor{black}{SDE~\eqref{eq:ME_cylindermodel} yields in its singular limit an} ODE, \textcolor{black}{which} has a globally 
attracting limit cycle at $y=0$ with contraction rate $\alpha > 0$; for $\alpha =0$, 
every trajectory is a periodic orbit at some $y \in \mathbb R$. In the presence of 
noise ($\sigma > 0$), the amplitude direction is driven by phase-dependent random 
perturbations. The real parameter $b$ induces an effect which is often 
called shear: if $b > 0$, the phase velocity depends on the amplitude $y$. \textcolor{black}{Note that for $\alpha=0$, the drift term of the $y$-component vanishes, while the second component has no noise component. This yields a very non-generic/singular coupling between a pure drift SDE and an ODE.}

In the tradition of random dynamical systems theory \cite{Arnold98}, and in contrast to
the sample paths approach in the last subsection, we now compare trajectories with different 
initial conditions but driven by the same noise. As trajectories depend on the noise 
realization, one cannot expect any convergent behavior of individual trajectories to 
a fixed attractor. An alternative point of view avoiding this problem is to consider, 
for a fixed noise realization in the past, the flow of a set of initial conditions 
from time $t=-T$ to a fixed endpoint in time, say $t=0$, and then take the (pullback) 
limit $T\to\infty$. If trajectories of initial conditions converge under this procedure 
to some set, then this set is called a random pullback attractor, or 
simply random attractor.

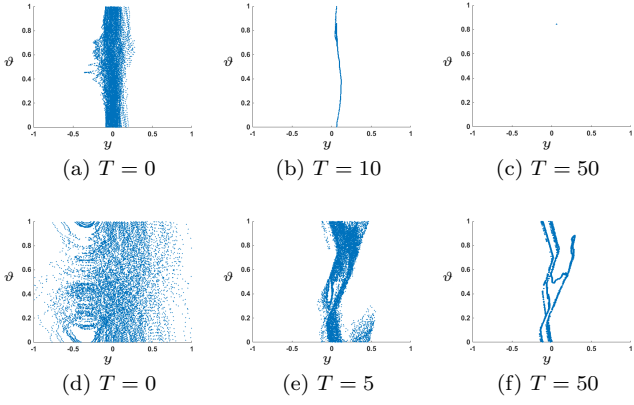
\begin{figure}[!ht]
\begin{center}
\subfloat[\label{fig:ME_synch1}$T=0$]{
\begin{overpic}[width=2.7cm,tics=10]{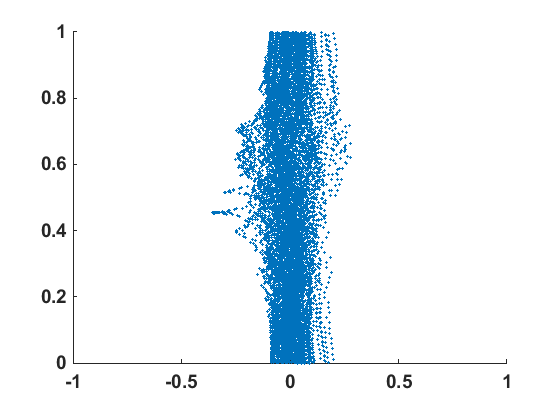}
\put(47,0){\tiny $y$}
	\put(-2,40){\tiny $\vartheta$}
\end{overpic} 
}
\subfloat[\label{fig:ME_synch2}$T=10$]{
\begin{overpic}[width=2.7cm,tics=10]{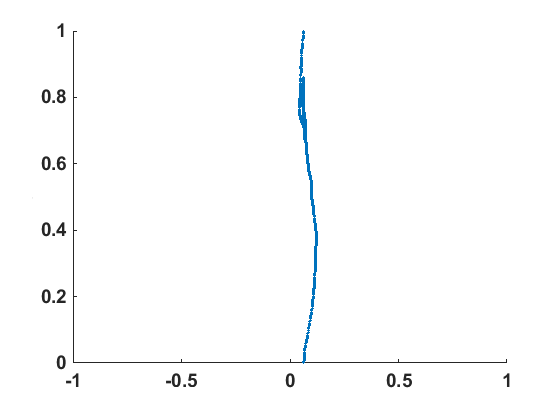}
\put(47,0){\tiny $y$}
	\put(-2,40){\tiny $\vartheta$}
\end{overpic} 
}
\subfloat[\label{fig:ME_synch3}$T=50$]{
\begin{overpic}[width=2.7cm,tics=10]{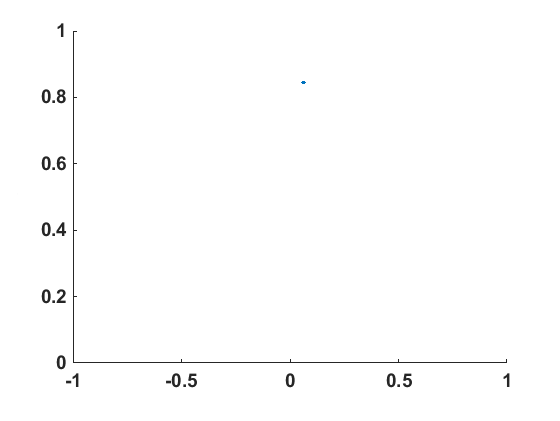}
\put(47,0){\tiny $y$}
	\put(-2,40){\tiny $\vartheta$}
\end{overpic} 
}\\
\subfloat[\label{fig:ME_strange1}$T=0$]{
\begin{overpic}[width=2.7cm,tics=10]{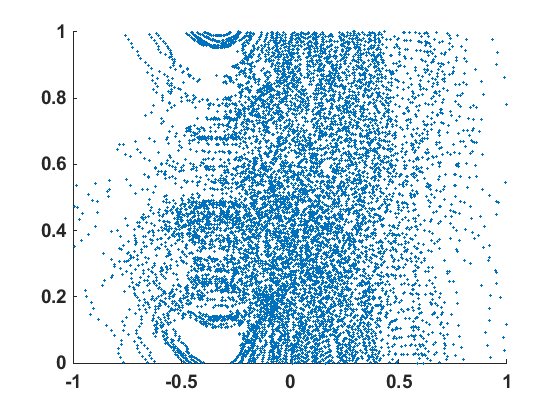}
\put(47,0){\tiny $y$}
	\put(-2,40){\tiny $\vartheta$}
\end{overpic} 
}
\subfloat[\label{fig:ME_strange2}$T=5$]{
\begin{overpic}[width=2.7cm,tics=10]{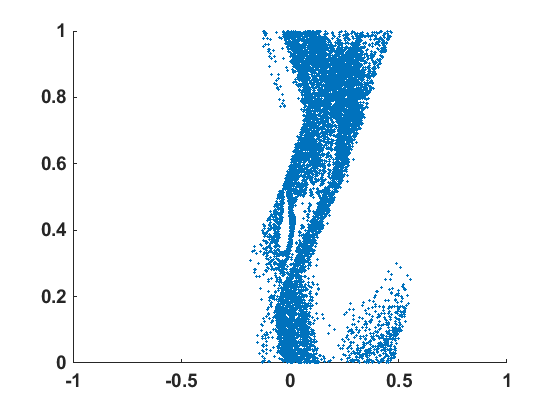}
\put(47,0){\tiny $y$}
	\put(-2,40){\tiny $\vartheta$}
\end{overpic} 
}
\subfloat[\label{fig:ME_strange3}$T=50$]{
\begin{overpic}[width=2.7cm,tics=10]{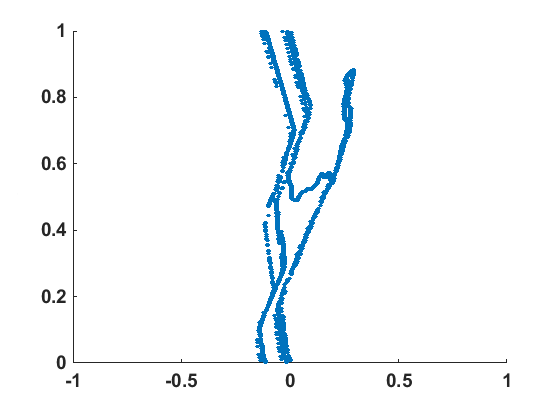}
\put(47,0){\tiny $y$}
	\put(-2,40){\tiny $\vartheta$}
\end{overpic} 
}
\end{center}
\caption{Pullback attraction to random equilibrium (a)-(c) for 
model~\eqref{eq:ME_cylindermodel} with $\sigma = 0.5, \alpha = 1.5, 
b = 3$ such that $\lambda_1 <0$, and to random strange attractor 
(d)-(f) for $\sigma = 2, \alpha = 1.5, b = 3$ such that $\lambda_1 >0$.}
\label{fig:ME_random_attractor}
\end{figure}

Typically, one can observe two different scenarios generated by the impact of 
noise on a stable limit cycle, as in model~\eqref{eq:ME_cylindermodel} with 
$\alpha > 0$: either synchronization of trajectories towards a random equilibrium 
(see Figure~\ref{fig:ME_random_attractor} (a)-(c)), or separation of trajectories 
within an attracting object, a random strange attractor with fractal properties 
(see Figure~\ref{fig:ME_random_attractor}(d)-(f)). The crucial quantity for 
determining the character of the dynamics is the sign of the first Lyapunov 
exponent $\lambda_1=\lambda_1(\alpha, b,\sigma)$ with respect to the ergodic 
invariant measure of the random system. The quantity $\lambda_1$ can be summarized 
as the dominant infinitesimal asymptotic growth rate of almost all trajectories.

The mechanism, whereby a combination of shear and noise leads to a positive 
Lyapunov exponent, was described as shear-induced chaos~\cite{LinYoung08}. The 
noise perturbations drive some points of the deterministic limit cycle up and 
some down on the cylinder. Due to the phase-amplitude coupling $b$, the points 
with larger $y$-coordinates move faster in the $\vartheta$-direction. At the 
same time, the dissipation force with strength $\alpha$ attracts the curve 
back to the limit cycles. This provides a mechanism for stretching and 
folding characteristic of chaos. The transition to chaos in the continuous 
time stochastic forcing is much faster than in the case of, e.g., periodic 
kicks~\cite{LinYoung08}. This is due to the fact that points end up in areas 
with arbitrarily large values of $y$ with positive probability such that already 
small shear can generate the described stretching and folding.
 
\begin{figure}[htpb]
\centering
\begin{tikzpicture}
[>=stealth',point/.style={circle,inner sep=0.035cm,fill=white,draw},
encircle/.style={circle,inner sep=0.07cm,draw},
x=4cm,y=4cm,declare function={f(\x) = {2*\x^(3/2)};}]
\draw[blue,thick,-,smooth,domain=0.001:0.52,samples=75,/pgf/fpu,/pgf/fpu/output
format=fixed] plot (\x, {f(\x)});

\draw[->,semithick,blue] (0,0) -> (1,0);
\draw[->,semithick] (0,0) -> (0,0.85);

\node[] at (0.9,-0.05) {$\alpha$};
\node[] at (-0.05,0.8) {$\sigma$};

\node[] at (0.7,0.57) {$\sigma = c_0^{-1/2} \mathbf \alpha^{3/2}$};

\draw [black] (0.5,0.3) circle [radius=0.06];
\node at (0.5,0.3) {I};
\draw [black] (0.2,0.72) circle [radius=0.06];
\node at (0.2,0.72) {II};
\draw [black] (0.9,0.1) circle [radius=0.06];
\node at (0.9,0.1) {III};
\draw[->] (0.85,0.07) -> (0.78,0.01);
\draw [black] (0.63,0.75) circle [radius=0.06];
\node at (0.63,0.75) {III};
\draw[->] (0.57,0.74) -> (0.5,0.68);

\end{tikzpicture}
\caption{Fixing $b=1$ in model~\eqref{eq:ME_cylindermodel}, the 
figure shows the areas of negative (I) and positive (II) 
$\lambda_1$ in the $(\alpha, \sigma)$-parameter space being 
separated by the curve $\{(\alpha, \sigma_0(\alpha,1))\}$ (III) 
of $\lambda_1$ being zero, using formula~\eqref{eq:ME_sigma0} 
for $\sigma_0$.}  \label{fig:ME_shear_cylinder}
\end{figure}
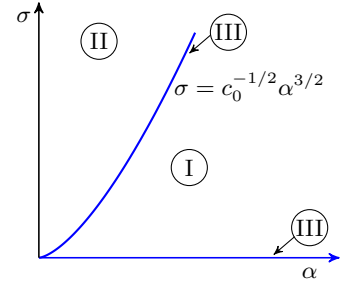

The validity of this mechanism has first been demonstrated 
analytically~\cite{WangYoung02,WangYoung03,OttStenlund10} in the case of periodically 
kicked limit cycles, including probabilistic characterizations of the dynamics. 
An analytical proof of shear-induced chaos in the stochastic setting was developed 
in~\cite{EngelLambRasmussen1}. Based on a specific machinery to explicitly express 
Lyapunov exponents for noisy oscillators~\cite{ImkellerLederer99, ImkellerLederer2001} 
one can provide the formula
\begin{align}
\lambda_1(\alpha,b, \sigma) &= - \frac{\alpha}{2} 
+ \frac{ b \sigma}{2} \int_{0}^{\infty}  v \ m_{\sigma, b, \alpha}(v)  
\,\rmd v\,, \label{eq:ME_topLyap}\\
\lambda_2(\alpha,b, \sigma) &= - \frac{\alpha}{2} 
- \frac{b \sigma}{2} \int_{0}^{\infty}  v \ m_{\sigma, b, \alpha}(v)  
\,\rmd v\,, \label{eq:ME_secLyap}
\end{align}
where 
\begin{equation}\label{eq:ME_msba}
m_{\sigma, b, \alpha}(v)= \frac{ \frac{1}{\sqrt{v}} 
\exp \left( - \frac{b \sigma}{6} v^3 + \frac{\alpha^2}{2b \sigma} v \right)  } 
{\int_{0}^{\infty} \frac{1}{\sqrt{u}}\exp \left( - \frac{b \sigma}{6} u^3 
+ \frac{\alpha^2}{2 b \sigma} u \right)  \rmd u},
\end{equation}
and $\lambda_2$ is the second Lyapunov exponent, which is always negative 
unless $\alpha = \sigma =0$. Furthermore, one can  
prove the following result~\cite{EngelLambRasmussen1}: assume 
the functions $f_i: \mathbb S^1\simeq [0,1) \to \mathbb{R}$  to be $C^{2,\kappa}$ 
for some $0 < \kappa \leq 1$ to guarantee differentiability of the random dynamical 
system (see \cite[Theorem 2.3.32] {Arnold98}), and, to make explicit calculations 
possible, assume $m\geq2$ with
\begin{equation} 
\label{eq:ME_sumcondition}
 \sum_{i=1}^m f_i'(\vartheta)^2 = 1 \ \fa \vartheta \in \mathbb S^1\,.
\end{equation}
Then there is $c_0 \approx 0.2823$ such that for all $\alpha, b > 0$, 
the number  
\begin{equation} 
\label{eq:ME_sigma0}
    \sigma_{0}(\alpha, b) = \frac{\alpha^{3/2}}{c_0^{1/2} b} > 0,
\end{equation} 
is the unique value of $\sigma$ where the top Lyapunov exponent 
$\lambda_1(\alpha, b,\sigma)$ of \eqref{eq:ME_cylindermodel} changes sign:
\begin{equation*}
\lambda_1(\alpha, b,\sigma) \begin{cases}
& < 0 \quad \text{if} \ 0 < \sigma < \sigma_0(\alpha, b)\,, \\
& = 0 \quad \text{if} \ \sigma = \sigma_0(\alpha, b)\,, \\
& > 0 \quad \text{if} \ \sigma > \sigma_0(\alpha, b)\,.
\end{cases}
\end{equation*}
In particular, we can just use the sign of the top Lyapunov exponent as a
definition of a property $\cP_{\textnormal{sic}}$ for the 
shear-induced chaos problem~\eqref{eq:ME_cylindermodel}. 
Figure~\ref{fig:ME_shear_cylinder} shows the graph of $\sigma_0$ for 
$0\leq \alpha \leq 1$ and fixed $b=1$. Note that for $b,\sigma\neq0$, we can 
always conduct a change of variables in the amplitude variable $y$ to rescale 
the shear parameter $b$ to $1$ and the effective noise amplitude to $\sigma b$. 
Hence, the above result and the corresponding illustration in 
Figure~\ref{fig:ME_shear_cylinder} hold in precisely the same way, when the 
roles of $\sigma$ and $b$ are exchanged. 

For all fixed $\alpha > 0$, if $\sigma =0$, i.e., in the zero-noise limit, 
we clearly have $\lambda_1 = 0$, now seen as the leading Lyapunov exponent 
associated with the attracting deterministic limit cycle. The convergence 
can also be seen by a different form of formula~\eqref{eq:ME_topLyap}, which 
is obtained by a change of variables as
\begin{equation}
\label{eq:ME_lambda1mainthm}
\lambda_1(\alpha, b, \sigma) = \frac{\alpha}{2} \left( \int_{0}^{\infty}  
u \ \tilde{m}_{\sigma,b, \alpha}(u)  \,\rmd u -1 \right),
\end{equation}
where
\begin{equation*}
\tilde{m}_{\sigma,b , \alpha}(u)= \frac{\frac{1}{\sqrt{u}} 
\exp \left( -  \frac{\alpha^3}{\sigma^2 b^2} \left[ \frac{1}{6} u^3 
- \frac{1}{2} u \right] \right)}{\int_0^{\infty} \frac{1}{\sqrt{w}} 
\exp \left( -  \frac{\alpha^3}{\sigma^2 b^2} \left[ \frac{1}{6} w^3 
- \frac{1}{2} w \right] \right)\,\rmd w}\,.
\end{equation*}
Hence, there is a continuous transition back to situation (III) at 
the $\alpha$-axis. When $\alpha = 0$ but $\sigma >0$, dissipativity and 
the existence of a random attractor are lost and the system becomes 
volume-preserving. Still, the associated first Lyapunov exponent $\lambda_1$ 
is positive and the $\sigma$-axis belongs to situation (II), as can be 
easily seen from formula~\eqref{eq:ME_topLyap}. The origin $(\sigma, \alpha) = (0,0)$ 
itself belongs to (III). This gives a full categorization of 
model~\eqref{eq:ME_cylindermodel} in terms of the first Lyapunov exponent 
under the double limit of the parameters $\sigma,b$ on the one side and 
$\alpha$ on the other.

Generally, shear-induced chaos can take more complicated forms with more 
nonlinearities. A paradigm problem is the normal form of a Hopf bifurcation 
with additive noise
\begin{equation}
\label{eq:ME_NormalForm}\tag{{$\cX_{\textnormal{sH}}$}}
\begin{array}{rl}
&\rmd x = (\alpha x - \beta y - (ax-by)(x^2 + y^2))\,\rmd t + \sigma \,\rmd W_t^1\,,\\
&\rmd y =  (\alpha y + \beta x - (bx+ay)(x^2 + y^2))\, \rmd t +  \sigma \,\rmd W_t^2\,,
\end{array}
\end{equation}
where $\sigma \geq 0$ is the strength of the noise, $\alpha\in\mathbb{R}$ 
equals the real part of eigenvalues of the linearization of the vector 
field at $(0,0)$, $b\in\mathbb{R}$ represents shear strength, $\beta\in \mathbb{R}$ 
is the linear component of rotational speed and $W_t^1, W_t^2$ denote independent 
one-dimensional Brownian motions. For $\alpha > 0$, the deterministic system 
($\sigma=0$) possesses a limit cycle at radius $\sqrt{\alpha}/a$, for any 
fixed $a >0$, with linear contraction rate $-2 \alpha$.

The model has been studied in \cite{DeVilleSriRapti11, DoanEngelLambRasmussen, Wiezcorek09} 
with various, predominantly numerical, approaches to describing shear-induced chaos. Hence,
it again makes sense to define $\cP_{\textnormal{sH}}$ via the sign of the first Lyapunov
exponent. For~\eqref{eq:ME_NormalForm}, only the case of synchronization, i.e.~$\lambda_1 < 0$, 
has been proven analytically \cite{DoanEngelLambRasmussen}. The change of sign of $\lambda_1$ to 
positive values is only proven in the particular context of the conditioned Lyapunov 
exponent~\cite{EngelLambRasmussen2}, considering the random dynamics on a bounded 
domain with killing at the boundary, by conducting a computer-assisted proof \cite{BredenEngel}. 
An explicit formula as before seems out of scope for system~\eqref{eq:ME_NormalForm} 
on the whole domain. 

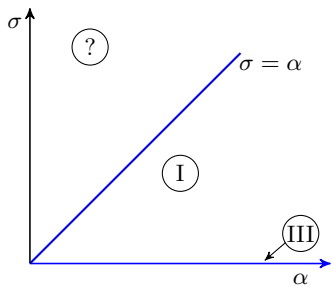
\begin{figure}[htpb]
\centering
\begin{tikzpicture}
[>=stealth',point/.style={circle,inner sep=0.035cm,fill=white,draw},
encircle/.style={circle,inner sep=0.07cm,draw},
x=4cm,y=4cm,declare function={f(\x) = \x;}]
\draw[blue,thick,-,smooth,domain=0.001:0.7,samples=75,/pgf/fpu,/pgf/fpu/output
format=fixed] plot (\x, {f(\x)});

\draw[->,semithick,blue] (0,0) -> (1,0);
\draw[->,semithick] (0,0) -> (0,0.85);

\node[] at (0.9,-0.05) {$\alpha$};
\node[] at (-0.05,0.8) {$\sigma$};

\node[] at (0.8,0.66) {$\sigma = \alpha$};

\draw [black] (0.5,0.3) circle [radius=0.06];
\node at (0.5,0.3) {I};
\draw [black] (0.2,0.72) circle [radius=0.06];
\node at (0.2,0.72) {?};
\draw [black] (0.9,0.1) circle [radius=0.06];
\node at (0.9,0.1) {III};
\draw[->] (0.85,0.07) -> (0.78,0.01);

\end{tikzpicture}
\caption{Fixing all other parameters in model~\eqref{eq:ME_NormalForm}, in particular $b \gg \sqrt{2}a$, 
we consider the $(\alpha, \sigma)$-parameter space for $\alpha, \sigma$ 
sufficiently small, and can associate the area beneath the diagonal with 
negative $\lambda_1$ (I) and the $\alpha$-axis, including the origin, 
with $\lambda_1=0$ (III).}  \label{fig:ME_noise_normalform}
\end{figure}

However, there are two small parameter results that give some indication concerning 
the question of double limits in this case and demonstrate the differences to the 
cylinder model. Firstly, it was shown in~\cite{DeVilleSriRapti11} and then further 
elaborated in~\cite{DoanEngelLambRasmussen} that for any fixed $a>0$, $b < \sqrt{2}a$ 
and $\alpha$ smaller than a given bound depending on all other parameters, the first 
Lyapunov exponent is negative, i.e.~$\lambda_1 < 0$. 
This means that for the case $b, \alpha \to 0$ we will always be in scenario (I), in 
contrast to model~\eqref{eq:ME_cylindermodel} where scenario (II) can happen in the 
double-limiting case, as illustrated in Figure~\ref{fig:ME_shear_cylinder} --- 
recall that $\sigma$ and $b$ are interchangeable in this case and the same formula 
and corresponding figure are also true for replacing $\sigma$ by $b$. This does not 
transfer to the more complicated, highly nonlinear situation of model~\eqref{eq:ME_NormalForm}. 
Secondly, Deville et al. \cite{DeVilleSriRapti11} demonstrate that $\lambda_1 < 0$ 
for $\sigma \frac{a}{\alpha} \to 0$. This allows us to give at least a partial 
picture of the small parameter situation for $\alpha, \sigma$ when the shear strength $b \gg \sqrt{2}a$ is large; Figure~\ref{fig:ME_noise_normalform} depicts such a sufficiently small domain in parameter space. Analytical approximation 
of other areas than the one beneath the diagonal seems out of reach with current 
methods.\medskip

The examples involving SDEs have shown clearly that small noise is a \textcolor{black}{common} source of
double limits. Yet, SDEs still carry some regularity due to the (almost $1/2$-H\"older) 
continuous input. The next subsection illustrates that even for stochastic switching 
problems one can frequently identify double limits.


\subsection{Piecewise Deterministic Processes}
\label{ssec:PDMP}

Piecewise deterministic processes are stochastic processes that evolve deterministically 
on most time intervals of short length; random events occur instantaneously and come 
for example in the shape of random switches between several driving vector fields, 
or jumps to randomly chosen sites of the phase space. In this subsection, we will consider 
two instances of piecewise deterministic processes, which are induced by a 
parameter-dependent differential equation with an intermittently-acting noise that 
depends itself on a small parameter. We work within the following basic framework: 
Let $M$ be an open subset of $\R^m$, $m \in \N$, and let $u_0$ and $u_1$ be smooth 
vector fields on $M$ that depend on a small positive parameter $\delta$. In addition, 
assume that for $i \in \{0,1\}$ and for every $x_0 \in M$, the initial-value problem 
\begin{align*}
\dot{x}(t) =& \ u_i(x(t)), \quad t > 0, \\
x(0) =& \ x_0
\end{align*}
has a unique solution $x(t) = \Phi_i^t(x_0)$ that is defined for all $t \geq 0$. 
Consider the differential equation 
\begin{equation}   
\label{eq:de} \tag{{$\cX_{\textnormal{pd}}$}}
\frac{\txtd x}{\txtd t} = U(\omega, x(t), t), 
\end{equation} 
where $\omega$ is a realization of a continuous-time Markov chain on $\{0,1\}$ with 
transition rates 
$$
\lambda_0 = \lim_{t \downarrow 0} \frac{\P(\omega_t = 1 \vert \omega_0 = 0)}{t}, 
\ \lambda_1 = \lim_{t \downarrow 0} \frac{\P(\omega_t = 0 \vert \omega_0 = 1)}{t}, 
$$
and where 
$$
U(\omega, x, t) := \begin{cases}
                                   u_0(x), & \quad \omega_t = 0, \\
                                   u_1(x), & \quad \omega_t = 1. 
                                   \end{cases}
$$
The differential equation in~\eqref{eq:de} is thus alternately driven by 
the vector fields $u_0$ and $u_1$, and switches between these vector fields 
correspond to the jumps of a continuous-time Markov chain. The latter being 
the only source of randomness, \textcolor{black}{we shall} assume that the transition 
rates $\lambda_0$ and $\lambda_1$ depend on a second small parameter $\varepsilon > 0$. 
For a typical choice of $\omega$, the equation in~\eqref{eq:de} has a unique 
solution $X(\omega)$ that is defined for all $t \geq 0$. The resulting stochastic 
process $X = (X_t)_{t \geq 0}$ on $M$ can be turned into a Markov process by 
adjoining the process $E = (E_t)_{t \geq 0}$ on $\{0,1\}$ defined by $E_t(\omega) 
:= \omega_t$. The resulting two-component process $(X,E)$ on the state space 
$M \times \{0,1\}$ belongs to the class of piecewise deterministic Markov 
processes~\cite{Davis_article}.  

In line with standard terminology, a stationary distribution for $(X,E)$ is 
a probability measure $\mu$ on $M \times \{0,1\}$ such that for every Borel 
set $A \subset M$, $i \in \{0,1\}$, and $t \geq 0$, 
$$
\mu(A \times \{i\}) = \sum_{j \in \{0,1\}} \int_M \mathfrak{P}^t(x,j; 
A \times \{i\}) \ \mu(\txtd x \times \{j\}), 
$$
where $(\mathfrak{P}^t)_{t \geq 0}$ denotes the Markov semigroup of $(X,E)$.\medskip

Consider the dynamical system induced by randomly switching between the 
two-dimensional linear vector fields $u_i(x) = U_i x$, $i \in \{0,1\}$, 
where 
\be
\label{eq:delin} \tag{{$\cX_{\textnormal{pdl}}$}}
U_0 := \begin{pmatrix}
        -\delta & 1 \\
        0 & -\delta 
        \end{pmatrix}, 
\quad U_1 := 
       \begin{pmatrix}
         -\delta & 0 \\
         -1 & -\delta 
       \end{pmatrix}.   
\ee
The switching rates are assumed to be $\lambda_0 = \lambda_1 = 
\varepsilon^{-1}$, i.e., for small $\varepsilon$ we are in the 
regime of fast switching. This system belongs to the class of 
switching systems studied in \cite{Lawley}. Here, we present some of 
the main findings from~\cite{Lawley} using the viewpoint of double 
limits in $\varepsilon$ and $\delta$. \textcolor{black}{Note that the problem is singularly perturbed since for $\varepsilon\ra 0$, one effectively obtains a single limit ODE governed by the average of $U_0$ and $U_1$, while for $\delta=0$, the individual linear vector fields give rise to ODEs whose solutions are constant in one component.}

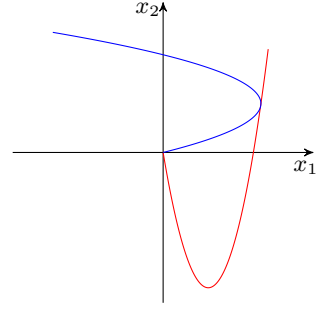
\begin{figure}[htbp]
	\centering
	\begin{tikzpicture}
[>=stealth',point/.style={circle,inner sep=0.01cm,fill=white,draw},
encircle/.style={circle,inner sep=0.07cm,draw},x=2cm,y=2cm]

\draw[->] (-1,0) -> (1,0);
\draw[->] (0,-1) -> (0,1);

\draw[red,-,smooth,domain=0:0.7,samples=200,/pgf/fpu,/pgf/fpu/output
format=fixed] plot (\x, {10*\x*(\x-0.6)});
\draw[blue,smooth,domain=0:0.8,samples=100] plot ({-0.65*\x*\x/(0.325)^2+1.3*\x/(0.325)},\x);

\node[] at (0.95,-0.1) {$x_1$};
\node[] at (-0.1,0.95) {$x_2$};

\end{tikzpicture}
	\caption{\label{fig:fig19TH} Sample trajectories for the vector fields $u_0$ and $u_1$ associated with $(\cX_{\textnormal{pdl}})$. The blue and red curves represent trajectories for $u_0$ and $u_1$, respectively. If one first flows along the blue curve towards the origin and then switches to the red one at the point where the curves touch, one can increase the distance to the origin.}
\end{figure}

Both $U_0$ and $U_1$ are defective matrices, meaning that the 
eigenspaces corresponding to their only eigenvalue $-\delta$ have 
dimension $1$. Since $-\delta < 0$, the equilibrium point $(0,0)$ 
shared by $u_0$ and $u_1$ is globally asymptotically stable for each 
individual ODE $\dot{x}(t) = u_i(x(t))$. However, as pointed out 
in \cite{Zitt}, \cite{Lawley} for the random case, and in \cite{Mason} 
for the deterministic case, switching between stable ODEs may cause 
instability. This phenomenon can be easily apprehended if switching 
takes place between two stable vector fields that admit an unstable 
average. As the switching rates tend to infinity, the random dynamics 
start to resemble the deterministic dynamics governed by the unstable 
average \cite{Gabrielli}. For the present system, however, the mechanism 
causing instability is more subtle (Figure \ref{fig:fig19TH}).  

Let us be more precise: We call the random dynamical system under 
consideration stable if the stochastic process $X$ on $\R^2$, induced by 
alternately flowing along $u_0$ and $u_1$, satisfies 
$$
\P_{x, i}\left(\lim_{t \to \infty} \| X_t \| = 0 \right) = 1 
$$
for every $x \in \R^2$ and $i \in \{0,1\}$. Here, $\P_{x,i}$ denotes the 
law of the Markov process $(X,E)$ starting at $(x,i)$, and $\| \cdot \|$ 
is the Euclidean norm on $\R^2$. The random dynamical system is said to 
be unstable if for every $x \in \R^2 \setminus \{(0,0)\}$ and $i \in \{0,1\}$, 
$$
\P_{x,i} \left(\lim_{t \to \infty} \| X_t \| = \infty \right) = 1. 
$$
A priori, there may be choices of $\varepsilon$ and $\delta$ for which the 
system is neither stable nor unstable. As we are about to see, this is, at 
least generically, not the case. We want to study the property 
$$
\cP_{\textnormal{pdl}}=\left\{\begin{array}{ll}
1\qquad & \textnormal{if the system is stable},\\
0\qquad & \textnormal{if the system is unstable.}
\end{array}\right.
$$
It is convenient to represent the stochastic process $X$ in polar 
coordinates (see \cite{Khasminskii} on the utility of polar 
decomposition for the study of Lyapunov exponents). 
Following \cite{Zitt}, one defines the radial process $R_t := \| X_t \|$ 
and the angular process $A_t := X_t/\|X_t\|$ whenever $X_t \neq (0,0)$. 
The two-component process $(A, E)$ on $S^1 \times \{0,1\}$ is then again 
a piecewise deterministic Markov process characterized by random switching 
between the vector fields $\theta \mapsto \sin^2(\theta)$ and 
$\theta \mapsto \cos^2(\theta)$, where $S^1$ is identified with the 
interval $[0, 2 \pi)$. According to \cite[Lemma~3.2]{Lawley}, $(A,E)$ 
admits a unique stationary distribution $\mu$ that is absolutely continuous 
with respect to the product of arc-length measure on $S^1$ and counting 
measure on $\{0,1\}$. In our example, $\mu$ only depends on the switching 
rate, i.e.\textcolor{black}{,} it is a function of $\varepsilon$ while being independent of 
$\delta$. Let $\rho$ be the probability density function of $\mu$ and 
let $\rho_i(\cdot) := \rho(\cdot, i)$ for $i \in \{0,1\}$. Since $\mu$ 
is $\varepsilon$-dependent, so are $\rho_0$ and $\rho_1$. Define 
\begin{equation}   
\label{eq:G_def} 
G(\varepsilon) := \int_0^{2 \pi} (\rho_0(\theta) - 
\rho_1(\theta)) \cos(\theta) \sin(\theta) \ \txtd \theta, 
\end{equation} 
which is set up in such a way that the integrand is positive for 
all $\theta \in [0, 2 \pi)$, and thus $G > 0$. From \cite[Lemma 3.3]{Lawley} 
one obtains the following cases: 
\begin{itemize}
 \item[(I)] If $\delta < G(\varepsilon)$, then $\cP_{\textnormal{pdl}} = 0$. 
 \item[(II)] If $\delta > G(\varepsilon)$, then $\cP_{\textnormal{pdl}} = 1$. 
 \end{itemize} 
There are explicit formulae for $\rho_0$ and $\rho_1$ \cite{Lawley}. 
Together with~\eqref{eq:G_def}, this yields a reasonably explicit 
representation for the threshold function $G$ that is in principle 
amenable to asymptotic analysis. 

\begin{figure}[htbp]
	\centering
\begin{tikzpicture}
[>=stealth',point/.style={circle,inner sep=0.01cm,fill=white,draw},
encircle/.style={circle,inner sep=0.07cm,draw},
x=4cm,y=4cm,declare function={f(\x) = {\x^(1/2)};}]
\draw[blue,thick,-,smooth,domain=0:0.9,samples=200,/pgf/fpu,/pgf/fpu/output
format=fixed] plot (\x, {f(\x)});
\draw[blue,thick] (0,0) -> (1,0);
\draw[dashed,thick] (0,0) -> (0,1);

\draw [fill=white] (0,0) circle [radius=0.02];

\node[] at (0.9,-0.05) {$\varepsilon$};
\node[] at (-0.04,0.95) {$\delta$};

\draw [black] (0.6,0.3) circle [radius=0.06];
\node at (0.6,0.3) {I};

\draw [black,fill=white] (0.2,0.7) circle [radius=0.06];
\node at (0.2,0.7) {II};

\end{tikzpicture}
	\caption{\label{fig:switching_01} Classification diagram with respect 
	to the property $\cP_{\textnormal{pdl}}$. The two regions (I) and (II) 
	correspond to the cases (I) (stable) and (II) (unstable). The blue curve 
	separating the regions (I) and (II) represents the graph of $G$. \textcolor{black}{We have not attempted to accurately render the asymptotic behavior for the graph of $G$ here.} On the 
	$\delta$-axis, the problem is not well-defined, which makes a classification 
	impossible.}
\end{figure}
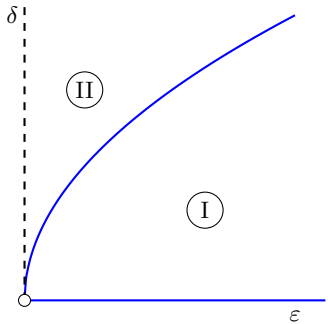 

If $\delta = 0$, $\varepsilon > 0$, the process $X$ alternately moves 
along lines parallel to the $x$-axis and lines parallel to the $y$-axis. 
It is not hard to see that \cite[Lemma 3.3]{Lawley} remains valid in this 
limiting case. Since $G(\varepsilon) > 0$, one has $\cP_{\textnormal{pdl}} = 0$. 

If $\varepsilon = 0$, the process $X$ is not well-defined because 
the switching rates are infinite. It does, however, make sense to study 
the limiting behavior of the random dynamical system as $\varepsilon \to 0$. 
According to~\cite[Thm.~2.5]{Lawley}, for $\varepsilon$ sufficiently small 
(with the required smallness depending on $\delta$), one has 
$\cP_{\textnormal{pdl}} = 1$. This implies that $\lim_{\varepsilon \to 0} 
G(\varepsilon) = 0$.  

Finally, we examine the situation when $\delta = 0$ and $\varepsilon \to 0$. 
In the classification diagram in Figure \ref{fig:switching_01}, this 
corresponds to approaching the origin along the $\varepsilon$-axis. By 
the averaging principle alluded to earlier~\cite[Thm.~2.1]{Gabrielli}, the 
process $X$ converges in probability, uniformly on compact time intervals, 
to the deterministic solution of the averaged problem 
\begin{align*}
\dot{x}(t) =& \frac{1}{2} (U_0 + U_1) x(t) = 
\begin{pmatrix} 
 0 & 1/2 \\
 -1/2 & 0 
\end{pmatrix}   x(t), \\
x(0) =& x_0.  
\end{align*}
The matrices $U_0$ and $U_1$ contribute equally to the averaged matrix 
$\tfrac{1}{2} (U_0 + U_1)$ because $\lambda_0 = \lambda_1$. The 
eigenvalues of the averaged matrix are $\pm \tfrac{i}{2}$, with zero 
real part. In this doubly singular situation, the previously observed 
dichotomy is broken: For every $x_0 \neq 0$, the trajectory of the solution 
to the averaged problem is a periodic orbit, more precisely a circle 
of radius $\| x_0 \|$ centered at the origin.\\ 

As second example for a piecewise deterministic Markov process, we are going to use a logistic growth model with random switching. Just as our first example, this Markov process is characterized by random switching between two vector fields with a critical point in common. Unlike the first example, though, the vector fields share a compact trapping region of positive Lebesgue measure that gives rise to a nontrivial stationary distribution.  

The logistic model 
is a classical model for the growth of a population that is limited by the 
capacity of the environment to sustain the population. The model is described 
by the logistic differential equation  
$\dot{x}(t) = \mathcal{U}(x(t), r, p)$, where 
$$
\mathcal{U}(x, r, p) := r x (1 - x/p).
$$
The time-dependent variable $x$ represents the population size. The parameters 
$r$ (the growth rate) and $p$ (the carrying capacity) are assumed to be positive. 

We consider the dynamical system induced by randomly switching between the 
logistic vector fields 
\be
\tag{{$\cX_{\textnormal{pdp}}$}}
\begin{array}{lcl}
u_0(x) &:=& \mathcal{U}(x,\delta, 1),\\
u_1(x) &:=& \mathcal{U}(x,1,2), 
\end{array}
\ee
at switching rates $\lambda_0 = \varepsilon$ and 
$\lambda_1 = 1$. Notice the asymmetry in the switching rates that will 
lead to the system spending more and more time in the regime governed 
by $u_0$ as $\varepsilon$ approaches $0$. In \cite{Hurth_Kuehn}, random 
switching between the vector fields $\mathcal{U}(\cdot, p_-, p_-)$ and 
$\mathcal{U}(\cdot, p_+, p_+)$ was studied in detail, for parameters 
$p_- < 0$ and $p_+ > 0$ to the left and to the right of the transcritical 
bifurcation at $p=0$. Even though the present setting is somewhat different, 
we will follow \cite{Hurth_Kuehn} quite closely.

For $r, p > 0$, the logistic vector field $\mathcal{U}(\cdot, r, p)$ has the equilibrium points $0$ and $p$, which are unstable and asymptotically stable, respectively. Stability of $1$ and $2$ for $u_0$ and $u_1$ implies that the compact interval $[1, 2]$ is positively invariant under the switching dynamics, i.e.\textcolor{black}{,} every switching trajectory starting in $[1, 2]$ stays in this interval for all positive times. Since, in addition, the Markov semigroup of $(X,E)$ is Feller (see Proposition 2.1 in \cite{Benaim}), the Krylov--Bogoliubov method (Theorem 3.1.1 in \cite{DaPrato}) yields the existence of a stationary distribution $\mu$ such that $\mu([1, 2] \times \{0,1\}) = 1$. Moreover, by \cite[Theorem 2]{Bakhtin} or by \cite[Theorem 4.4]{Benaim}, $\mu$ is the unique stationary distribution for $(X,E)$ that assigns full measure to $(0, \infty) \times \{0,1\}$. Finally, again by \cite[Theorem 2]{Bakhtin}, $\mu$ is absolutely continuous with respect to the product of Lebesgue measure on $(0, \infty)$ and counting 
measure on $\{0,1\}$. Hence, $\mu$ admits a density $\rho$ with respect to the latter measure.

For the 
invariant density $\rho_0(\cdot) := \rho(\cdot, 0)$, we consider 
the property 
$$
\cP_{\textnormal{bdd}} = \begin{cases}
  1, & \textnormal{if $\rho_0$ is bounded on $(1,2)$}, \\
  0, & \textnormal{if $\rho_0$ is unbounded on $(1,2)$.}
\end{cases}
$$ 
By \cite[Thm.~1]{Mattingly}, $\rho_0$ and $\rho_1 := \rho(\cdot, 1)$ 
are $C^{\infty}$ smooth in the open interval $(1, 2)$ because $u_0$ 
and $u_1$ are smooth vector fields with no equilibrium points in $(1,2)$. 
As a result, the corresponding probability fluxes 
$\varphi_i := \rho_i u_i$, $i \in \{0,1\}$, satisfy the Fokker--Planck 
equations \cite{Faggionato} 
\begin{equation}     
\label{eq:Fokker_Planck} 
\varphi_i'(x) = - \left(\frac{\varepsilon}{u_0(x)} 
+ \frac{1}{u_1(x)} \right) \varphi_i(x),  
\end{equation} 
for all $x \in (1, 2)$. 
The ODE in~\eqref{eq:Fokker_Planck} has the general solution
$$
\varphi_i(x) = C x^{- \frac{\varepsilon}{\delta} - 1} 
(x - 1)^{\frac{\varepsilon}{\delta}} (2 - x), \quad x \in (1,2), 
$$
hence 
\begin{align*}
\rho_0(x) =& c_1 x^{-\frac{\varepsilon}{\delta} -2} 
(x - 1)^{\frac{\varepsilon}{\delta} -1} (2 -x),  \\
\rho_1(x) =& c_2 x^{-\frac{\varepsilon}{\delta} -2} 
(x-1)^{\frac{\varepsilon}{\delta}}, 
\end{align*}
for positive normalizing constant $c_1$ and $c_2$. These formulae 
for $\rho_0$ and $\rho_1$ show that $\rho_1$ is always bounded 
on $(1,2)$. Furthermore, the invariant density $\rho_0$ has a 
singularity at the equilibrium point $1$ of $u_0$ if and only 
if $\varepsilon < \delta$. We obtain the following cases: 
\begin{itemize}
 \item[(I)] If $\delta \leq \varepsilon$, then $\cP_{\textnormal{bdd}} = 1$. 
 \item[(II)] If $\delta > \varepsilon$, then $\cP_{\textnormal{bdd}} = 0$. 
 \end{itemize} 
 
 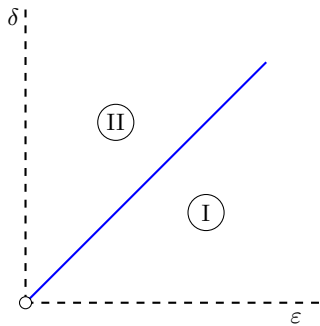
\begin{figure}[htbp]
	\centering
\begin{tikzpicture}
[>=stealth',point/.style={circle,inner sep=0.01cm,fill=white,draw},
encircle/.style={circle,inner sep=0.07cm,draw},
x=4cm,y=4cm]
\draw[dashed,thick] (0,0) -> (1,0);
\draw[dashed,thick] (0,0) -> (0,1);
\draw[blue,thick] (0,0) -> (0.8,0.8);
\draw [fill=white] (0,0) circle [radius=0.02];

\node[] at (0.9,-0.05) {$\varepsilon$};
\node[] at (-0.04,0.95) {$\delta$};

\draw [black] (0.6,0.3) circle [radius=0.06];
\node at (0.6,0.3) {I};

\draw [black,fill=white] (0.3,0.6) circle [radius=0.06];
\node at (0.3,0.6) {II};

\end{tikzpicture}
	\caption{\label{fig:switching_03} Classification diagram with 
	respect to the property $\cP_{\textnormal{bdd}}$. The two 
	regions I and II correspond to the cases (I) (bounded) and (II) (unbounded). 
	The blue ray separating the regions (I) and (II) belongs to region (I). 
	The property is not well-defined on the axes.}
\end{figure} 

This dichotomy admits the following heuristic explanation: 
If $\varepsilon$ (the rate of switching away from the vector 
field $u_0$) is small compared to $\delta$ (the contraction rate 
of $u_0$ at its equilibrium point $x=1$), then a large amount of 
probabilistic mass accumulates in the vicinity of the equilibrium 
point; a singularity at $x=1$ is formed. Conversely, if $\varepsilon$ 
is large in comparison with $\delta$, the system switches 
sufficiently often away from $u_0$ to prevent a strong accumulation 
of probabilistic mass near $x=1$; the invariant density $\rho_0$ 
stays bounded. 

In the singular case $\varepsilon = 0$, no switching away 
from $u_0$ takes place. The process $(X,E)$ still has a unique 
stationary distribution on $(0,\infty) \times \{0,1\}$, namely 
the product of the Dirac measure at $x=1$ and the measure $(1,0)$ 
on $\{0,1\}$. Of course, this distribution no longer has a 
probability density function with respect to the product of 
Lebesgue measure on $(0, \infty)$ and counting measure on $\{0,1\}$. 
It follows that the property $\cP_{\textnormal{bdd}}$ cannot 
be studied on the $\delta$-axis. 

If $\delta = 0$, the vector field $u_0$ is constantly equal to 
zero. As long as $\varepsilon > 0$, the system alternates between 
flowing along $u_1$ and staying put. The unique stationary 
distribution on $(0, \infty) \times \{0,1\}$ is then the product 
of the Dirac measure at the equilibrium point $x=2$ of the measure 
$(\tfrac{1}{1 + \varepsilon}, \tfrac{\varepsilon}{1 + \varepsilon})$ 
on $\{0,1\}$. Again, $\cP_{\textnormal{bdd}}$ cannot be 
meaningfully studied. Finally, in the doubly singular case 
$\varepsilon = \delta = 0$, one obtains an infinite family of 
stationary distributions $(\mu^x)_{x > 0}$, where $\mu^x$ is the 
product of the Dirac measure at $x$ and the measure $(1,0)$ on $\{0,1\}$. 

For switching systems 
in dimension greater than one, the set of singularities of 
invariant densities can have a much richer structure than the one 
exhibited here (see \cite{BHLM} for a simple yet nontrivial example in 2D). 
This can result in more complex classification diagrams with respect 
to a suitably defined version of $\cP_{\textnormal{bdd}}$.\medskip 

We conclude this subsection with some remarks on the two examples presented above. We also hint at additional topics in the field of piecewise deterministic processes where double limits may be fruitfully studied. 

In the first example, we saw that switching between vector fields of a certain kind (stable, in our example) can result in a dynamical system of a very different kind (unstable). In the same vein, for a Lotka--Volterra system of two competing species, it is shown in \cite{Lobry} that switching between two environments that both favor the same species can even lead to the extinction of this species. The articles \cite{Lobry}, \cite{Malrieu_Zitt}, and \cite{Phu} together provide a clear picture of which parameter choice results in which long-term behavior for the Lotka--Volterra system. It is thus possible to represent the interplay of the parameters by a double-limit diagram.    

The boundedness property for invariant densities is straightforward to study for piecewise deterministic processes of spatial dimension one \cite{Mattingly}. In higher dimensions, a regularity theory for invariant densities is still missing. However, the double-limits framework can also be meaningfully applied to other ergodic properties, e.g.\textcolor{black}{,} the number of stationary distributions, absolute continuity of stationary distributions with respect to a suitable reference measure, or exponential ergodicity. When it comes to the number of stationary distributions, an essential tool is the theory of stochastic persistence \cite{B17}, which gives criteria for the existence of a stationary distribution on the complement of an invariant closed subset of the phase space (the \textcolor{black}{so-called extinction set}). In \cite{Strickler_Benaim} and \cite{Strickler}, this theory -- originally devised for Markov processes in general -- has been further developed in the context of piecewise deterministic processes.  

In general, there is a lack of precise necessary conditions for absolute continuity and exponential ergodicity of the stationary distribution. Besides, neither of these properties is affected by the rates of switching, which makes it imperative to link both of the small parameters $\varepsilon$ and $\delta$ to the vector fields in order to obtain a nontrivial double-limit diagram. Apart from \cite{Benaim} and \cite{Bakhtin}, absolute continuity for piecewise deterministic processes was studied for instance in \cite{czapla} and \cite{Loecherbach}, where the process $X$ is allowed to have jumps. Sufficient conditions for exponential ergodicity in total-variation distance were given in \cite{Benaim}, \cite{Cui}, and \cite{BeHuSt2018}; and for exponential ergodicity in Wasserstein distance in \cite{Le_Borgne} and \cite{Cloez}. 

Other types of switching processes have been studied in the literature, some of them abundantly: switching between PDEs \cite{LawleyMattinglyReed2015}, non-Markovian switching \cite{Li}, switching between diffusions \cite{Yin}, etc. All of these classes of stochastic processes are amenable to the double-limit approach proposed in this article.   
\medskip 
 
We have already seen in the current context, that one expects double limit
problems for stochastic systems to be directly linked to double limits for
Fokker--Planck (or Kolmogorov) equations. We shall now continue with this
theme and focus in the next two subsections on problems arising from
various classes of partial differential equations (PDEs).

\subsection{Matched Asymptotic Expansions \& BVPs}
\label{ssec:matched}

Two-parameter singularly perturbed systems of differential equations have been widely studied from the analytical as well as from the numerical viewpoint (see \cite{Chen_1974,Gracia_2006,Herceg_2011,Kadalbajoo_2008,OMalley_1967,OMalley_1974,OMalley25,Roos_2003,Valarmathi_2003,Vulanovic_2001} and references therein). In most cases, the singularity is attributed to the presence of small parameters in front of the derivative terms; however, as shown in \cite{Popovic_2004}, this is not a necessary condition. This also applies to the problem presented in this section.

We start with a PDE problem, which still links to ODEs and classical
double limit fast-slow systems. We consider the following boundary value problem:
\begin{equation} 
\label{eq:AI_op} \tag{{$\cX_{\textnormal{mes}}$}}
\begin{aligned}
u_{XX}&=\frac\lambda{(1+u)^2}\bigg[1-\frac{\varepsilon^2}{(1+u)^2}\bigg],
 &&X\in[-1,1],\\
u&=0,  &&X=\mp1.
\end{aligned}
\end{equation}
Equation \eqref{eq:AI_op} describes the steady states associated to a 
second-order parabolic PDE problem arising in the context of Micro-Electro Mechanical 
Systems (MEMS) \cite{Li14}. In particular, the function $u(X)$ represents the 
deflection of an elastic membrane towards a ground plate under the action 
of an electric potential described by the parameter $\lambda > 0$, while 
$0 < \varepsilon \ll 1$ appears as a regularization parameter. \textcolor{black}{The problem is evidently singularly perturbed in $\lambda$, as for $\lambda=0$ it becomes just a trivial linear boundary value problem, while we shall see below that there is a hidden fast-slow singular perturbation structure with respect to $\varepsilon$.} 

\begin{figure}[!ht]
\begin{center}
\subfloat[\label{fig:AI_Lin:a}Bifurcation diagram.]{
\begin{overpic}[width=4.2cm,tics=10]{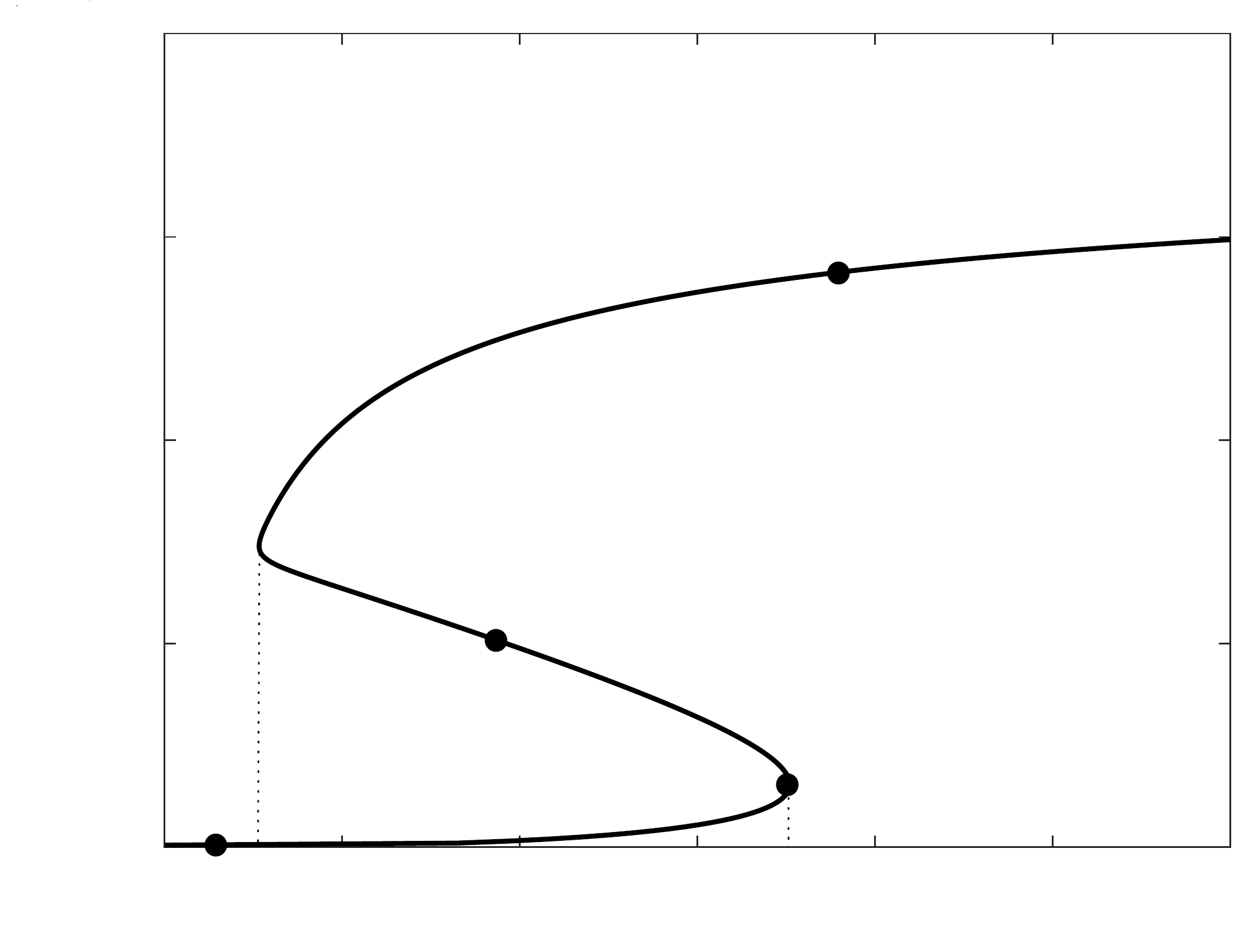}
\put(18,2){\scriptsize $\lambda_\ast$}
\put(60,2){\scriptsize $\lambda^\ast$}
\put(8,5){\tiny $0$}
\put(15.5,11.5){\scriptsize $a$}
\put(65,12){\scriptsize $b$}
\put(40,28){\scriptsize $c$}
\put(65,58){\scriptsize $d$}
\put(100,0){\scriptsize $\lambda$}
\put(-5,75){\scriptsize $\|u\|_2^2$}
\end{overpic} 
}
\subfloat[\label{fig:AI_Lin:b}Solutions in $(X,u)$-space.]{
\begin{overpic}[width=4.2cm,tics=10]{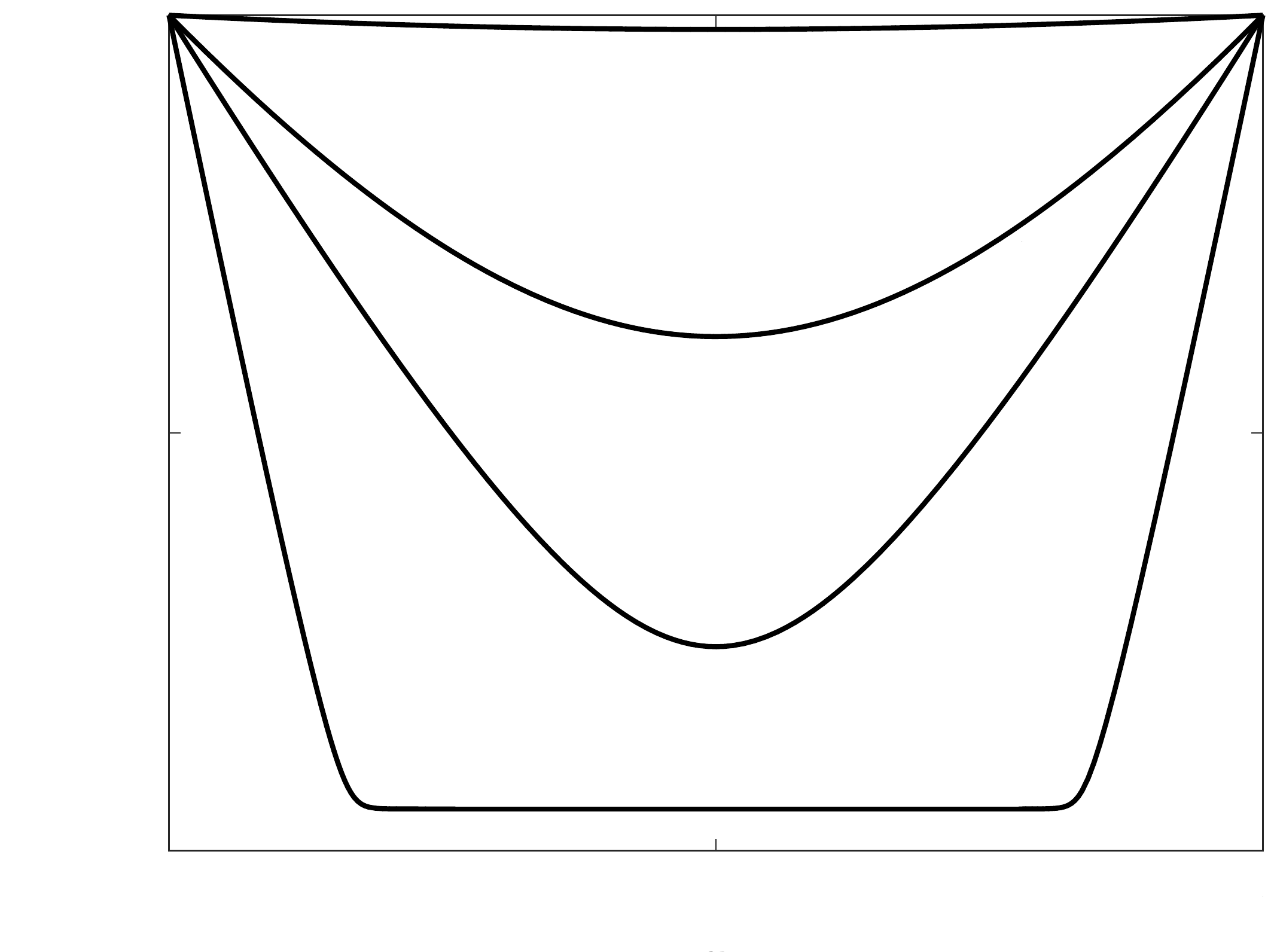}
\put(5,5){\tiny $-1$}
\put(80,68){\scriptsize $a$}
\put(80,52){\scriptsize $b$}
\put(80,41){\scriptsize $c$}
\put(80,14){\scriptsize $d$}
\put(100,0){\scriptsize $X$}
\put(5,75){\scriptsize $u$}
\end{overpic} 
}
\end{center}
\caption{(a) Numerically computed bifurcation diagram of the 
one-dimensional membrane model, \eqref{eq:AI_op}, for $\varepsilon=0.05$. 
(b) Corresponding solutions in $(X,u)$-space.}
\label{fig:AI_Lin}
\end{figure}
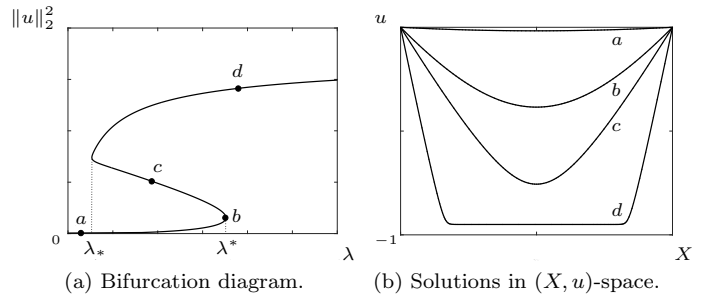

The bifurcation diagram associated to \eqref{eq:AI_op} consists of two branches 
of stable equilibria separated by a third, intermediate branch of unstable 
equilibria (see Figure \ref{fig:AI_Lin:a}). The middle and upper branch meet 
at a saddle-node bifurcation point $\lambda_\ast(\varepsilon)$. A steady-state 
solution exists for every $\lambda>0$, and the presence of the 
regularizing $\varepsilon$-dependent term in \eqref{eq:AI_op} guarantees 
that for any $\varepsilon$ the solution is bounded below by $u=-1+\varepsilon$; 
see Figure \ref{fig:AI_Lin:b}.

In \cite{Li14}, the authors have studied \eqref{eq:AI_op} both 
analytically, using matched asymptotic expansions to construct 
solutions, and numerically, investigating the structure of the 
$\varepsilon$-dependent bifurcation diagram. However, the analytical 
motivation behind logarithmic switchback terms in the expansions, as 
well as a detailed resolution of the bifurcation diagram for very small 
values of $\varepsilon$, were left as challenging open questions. 
In \cite{Iuorio_2019}, a detailed asymptotic resolution of 
\ref{fig:AI_Lin:a}, both in the singular limit of $\varepsilon=0$ 
and for $\varepsilon$ positive and sufficiently small, is accomplished 
through separate investigation of three distinct, yet overlapping, 
regions in the diagram, allowing us to tackle these questions.

To that end, we first reformulate the boundary value problem \eqref{eq:AI_op} 
in a dynamical systems framework; then, identification of two main parameters 
in the resulting equations yields a two-parameter singular perturbation 
problem. Careful asymptotic analysis of that problem allows us to identify 
the corresponding limiting solutions, and to show how the third branch in 
the diagram found for non-zero $\varepsilon$ emerges from the singular limit 
of $\varepsilon=0$, where only the lower and the middle branch are present. 
On that basis, we prove the existence and uniqueness of solutions close to 
these limiting solutions. 

We reformulate \eqref{eq:AI_op} as a first-order system by relabeling $u$ with $x$, introducing the variable $y=x_X$, and appending the trivial dynamics for the spatial variable $X$, which we relabel as $\xi$, and $\varepsilon$. Moreover, 
we desingularize the flow near $x=-1$ and define a shift in $x$ via 
$\tilde{x}=1+x$, which translates the singularity to $\tilde{x}=0$. 
Omitting the tilde and denoting differentiation with respect to the 
new independent variable by a prime, we obtain
\begin{subequations}
\label{eq:AI_sysl}
\begin{align}
x' &=x^4y, \label{eq:AI_sysl_a} \\
y' &=\lambda(x^2-\varepsilon^2), \label{eq:AI_sysl_b} \\
\xi' &=x^4, \label{eq:AI_sysl_c} \\
\varepsilon' &=0,
\end{align}
\end{subequations}
subject to the boundary conditions $x=1$ at $\xi=\mp1$. For $\varepsilon=0$, 
this systems admits the line of degenerate equilibria
\begin{equation}
\label{eq:AI_S0}
\mathcal{S}^0=\left\{(0,y,\xi,0) \, \big| \, y\in\mathbb{R},
\ \xi\in\mathbb{R} \right\}.
\end{equation}
When $\lambda=0$, there is an additional manifold of equilibria 
for~\eqref{eq:AI_sysl_a}-\eqref{eq:AI_sysl_b} given by
\begin{align}
\label{M0}
\mathcal{M}^0:=\left\{(x,0,\xi,0) \, \big| \, x\in\mathbb{R}^+,
\ \xi\in\mathbb{R} \right\}.
\end{align}
As it turns out, in two of the three regions we investigate it 
is useful to introduce a rescaled variable $\tilde{y}=\delta y$, where
\begin{align}
\label{eq:AI_delta}
\delta=\sqrt{\frac{\varepsilon}{\lambda}}.
\end{align}
Omitting the tilde for sake of simplicity, System \eqref{eq:AI_sysl} hence 
becomes
\begin{subequations}
\label{eq:AI_sysd}
\begin{align}
x' &=x^4y, \label{eq:AI_sysd_a} \\
y' &=\varepsilon(x^2-\varepsilon^2), \label{eq:AI_sysd_b} \\
\xi' &=\delta x^4, \label{eq:AI_sysd_c} \\
\varepsilon' &=0.
\end{align}
\end{subequations}
We observe that \eqref{eq:AI_sysd} is a fast-slow system, where 
$x$ is fast and $y$ is slow. The nature of $\xi$, however, depends 
on $\delta$: in particular, $\xi$ is fast when $\delta=\mathcal{O}(1)$, 
and it is slow when $\delta = \mathcal{O}(\varepsilon)$. For $\delta=0$, 
the manifolds $\mathcal{S}^0$ and $\mathcal{M}^0$ represent 
two branches of the critical manifold for \eqref{eq:AI_sysd}. Since 
$\mathcal{S}^0$ is not normally hyperbolic, and the reduced flow on 
it is highly degenerate, one can apply the blow-up method 
to describe the dynamics of \eqref{eq:AI_sysd} in its 
vicinity~\cite{Du93,DR96,KS01}. Such method has proved to be particularly useful when tackling two-parameter perturbed systems \cite{Kosiuk_2011,MiaoPopovicSzmolyan}. To this aim, we introduce the following 
blow-up transformation:
\begin{equation}\label{eq:AI_blowup}
x=\bar r\bar x,\quad y=\bar y,\quad\xi=\bar \xi,\quad
\text{and}\quad\varepsilon=\bar r\bar\varepsilon,
\end{equation}
where $(\bar y, \bar \xi) \in \mathbb{R}^2$ and $(\bar x, \bar \varepsilon) 
\in S^1$, { i.e.}, $\bar x^2+\bar\varepsilon^2=1$. The vector field induced 
by \eqref{eq:AI_blowup} on the cylindrical manifold in $(\bar x,\bar y,
\bar\xi,\bar\varepsilon,\bar r)$-space is best described in coordinate 
charts; in particular, to carry out our analysis we require the two following 
charts:
\begin{subequations}
\label{eq:AI_charts}
\begin{align}
K_1:\ & (x,y,\xi,\varepsilon)=(r_1,y_1,\xi_1,r_1\varepsilon_1),
\label{eq:AI_K1} \\
K_2:\ & (x,y,\xi,\varepsilon)=(r_2x_2,y_2,\xi_2,r_2). 
\label{eq:AI_K2}
\end{align}
\end{subequations}
We note that the phase-directional chart $K_1$ describes the 
``outer'' regime, which corresponds to the transient from $x=1$ to $x=0$ 
approaching $\mathcal{S}^0$, while the rescaling chart $K_2$ covers the 
``inner'' regime where $x \approx 0$, in the context of \eqref{eq:AI_sysd}. 
The corresponding dynamics are given by
\begin{equation}
  K_1: \, \left\{ \begin{aligned}
                    r_1' &=r_1y_1, \\
                    y_1' &=\varepsilon_1(1-\varepsilon_1^2), \\
                    \xi_1' &=\delta r_1, \\
                    \varepsilon_1' &=-\varepsilon_1y_1.
             \end{aligned} \right.
  \qquad K_2: \, \left\{ \begin{aligned}
                    x_2' &=x_2^4y_2, \\
                    y_2' &=x_2^2-1, \\
                    \xi_2' &=\delta r_2x_2^4, \\
                    r_2' &=0.
             \end{aligned} \right.
 \end{equation}
In order to construct singular solutions, we define the entry and exit 
sections in $K_1$
\begin{equation}
\label{eq:AI_s1in}
\begin{split}
\Sigma_1^{\rm in} :=&\left\{(\rho,y_1,\xi_1,\varepsilon_1) 
\, \big| \, y_1\in[y_-,y_+], \right. \\
&\left. \xi_1\in[\xi_-,\xi_+], \ \varepsilon_1\in[0,\sigma] \right\},
\end{split}
\end{equation}
\begin{equation}
\label{eq:AI_s1out}
\begin{split}
\Sigma_1^{\rm out} :=& \left\{(r_1,y_1,\xi_1,\sigma) \, \big| 
\, r_1\in[0,\rho],\ \right. \\
&\left. y_1\in[y_-,y_+], \ \xi_1\in[\xi_-,\xi_+] \right\},
\end{split}
\end{equation}
where $0<\rho<1$ and $0<\sigma<1$ are appropriately defined constants, 
while $y_\mp$ and $\xi_\mp$ are real constants, with $y_-<-\frac2{\sqrt3}$ 
and $y_+>\frac2{\sqrt3}$. Translating $\Sigma_1^{\rm out}$ in terms of 
$K_2$-coordinates, we obtain the section
\begin{equation}
\label{eq:AI_s2in}
\begin{split}
\Sigma_2^{\rm in}:=&\left\{(\sigma^{-1},y_2,\xi_2,r_2) \, \big| 
\, y_2\in[y_-,y_+],\  \right. \\
&\left. \xi_2\in[\xi_-,\xi_+], \ r_2\in[0,\rho\sigma] \right\}.
\end{split}
\end{equation}
In terms of matched asymptotics, such sections describe the transition 
between outer and inner regions. In particular, the outer regime 
corresponds to the area limited by $\Sigma_1^{\rm in}$ and 
$\Sigma_1^{\rm out}$ in $K_1$, while the inner regime is limited 
by $\Sigma_2^{\rm in}$ and the hyperplane $\left\{ y=0 \right\}$ in $K_2$.

\begin{figure}[!ht]
\begin{center}
\subfloat[\label{fig:AI_ss1}Type M1, $\lambda=\mathcal{O}(\varepsilon)$.]{
\begin{overpic}[width=4cm,tics=10]{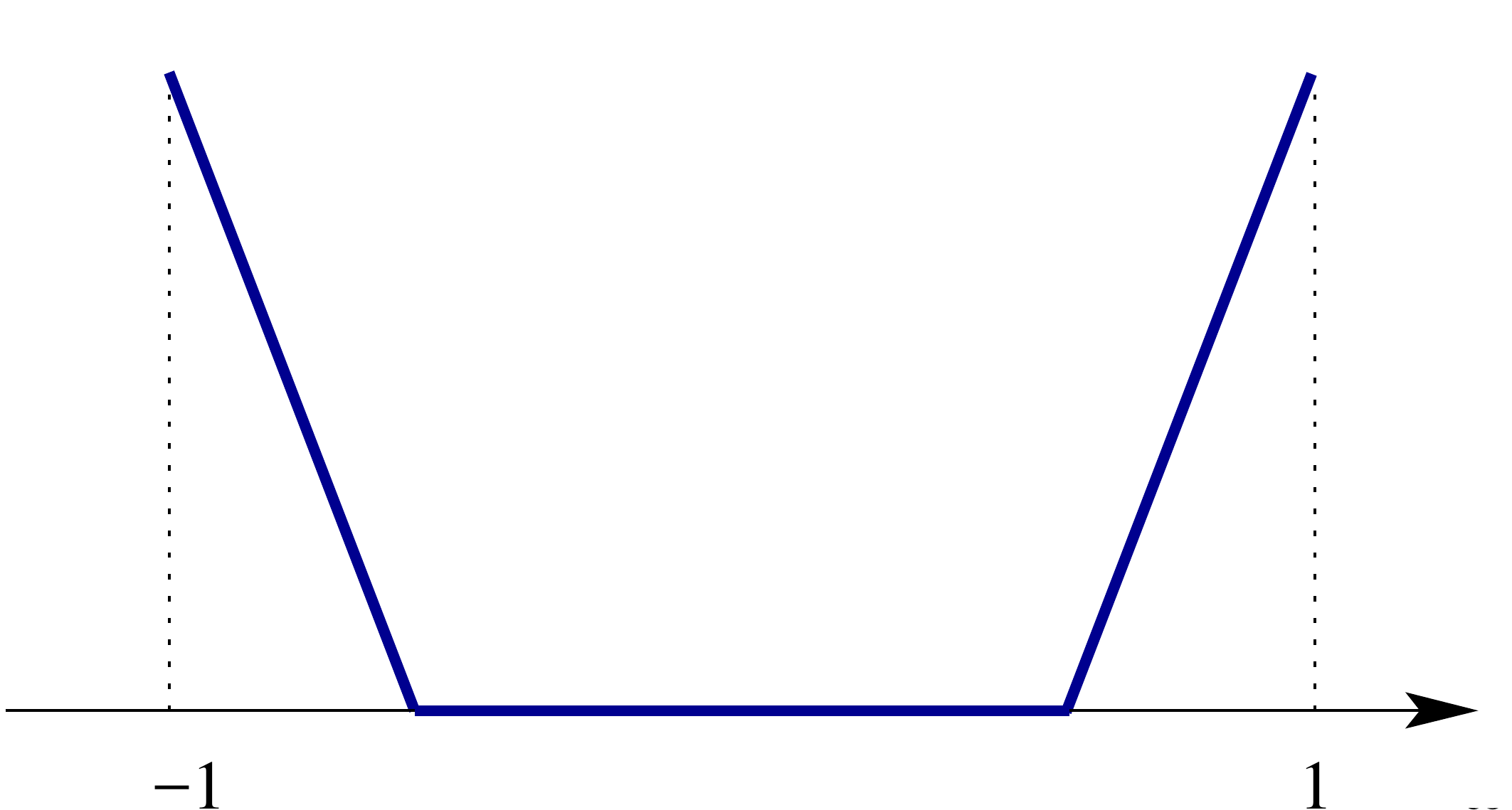}
\put(92,-1){\scriptsize $X$}
\put(2,50.5){\scriptsize $x$}
\end{overpic} 
}
\subfloat[\label{fig:AI_ss1bn}Type M1, $\lambda=\mathcal{O}(1)$.]{
\begin{overpic}[width=4cm,tics=10]{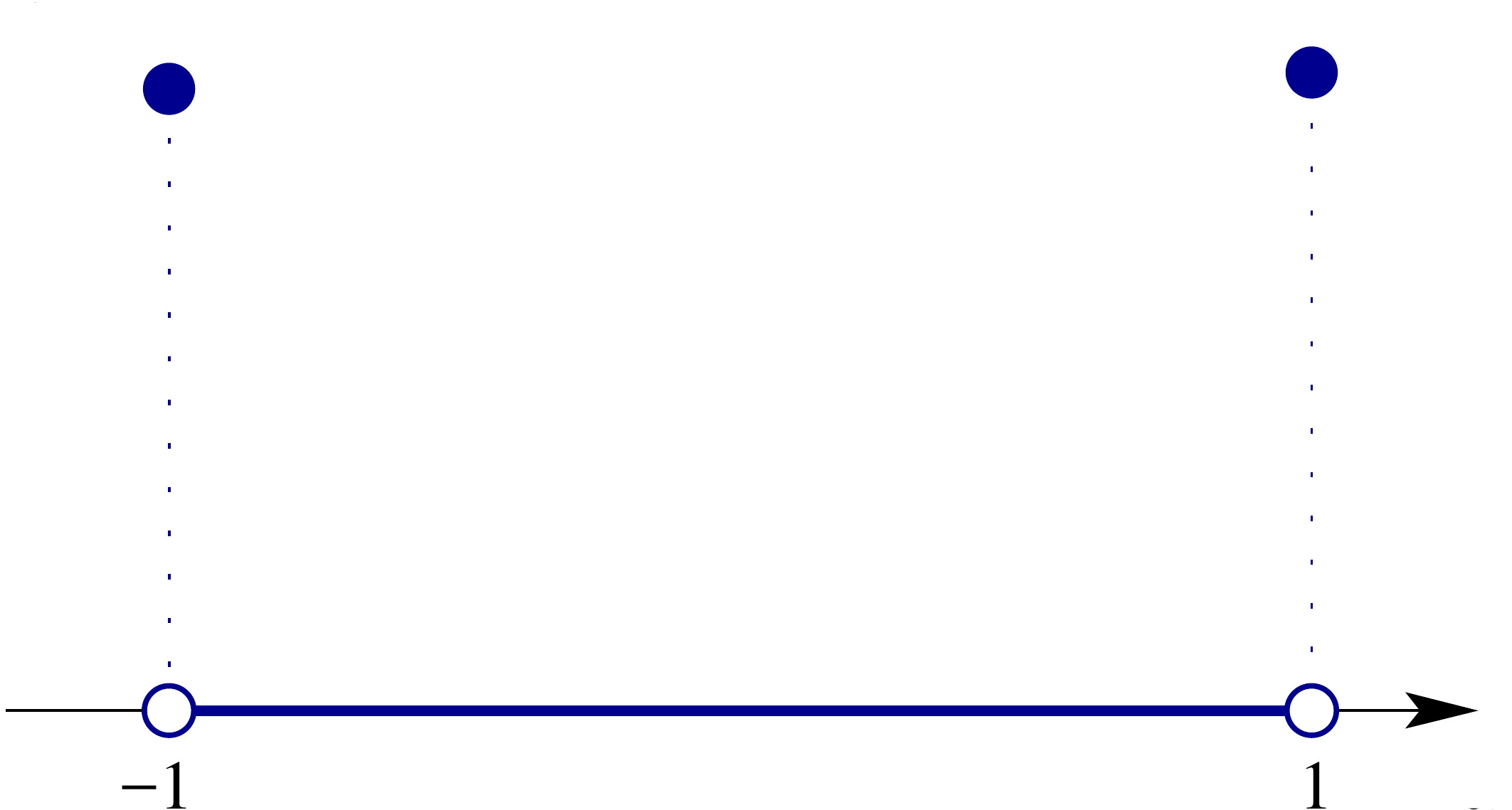}
\put(92,-1){\scriptsize $X$}
\put(2,50.5){\scriptsize $x$}
\end{overpic} 
}\\
\subfloat[\label{fig:AI_ss2}Type M2.]{
\begin{overpic}[width=4cm,tics=10]{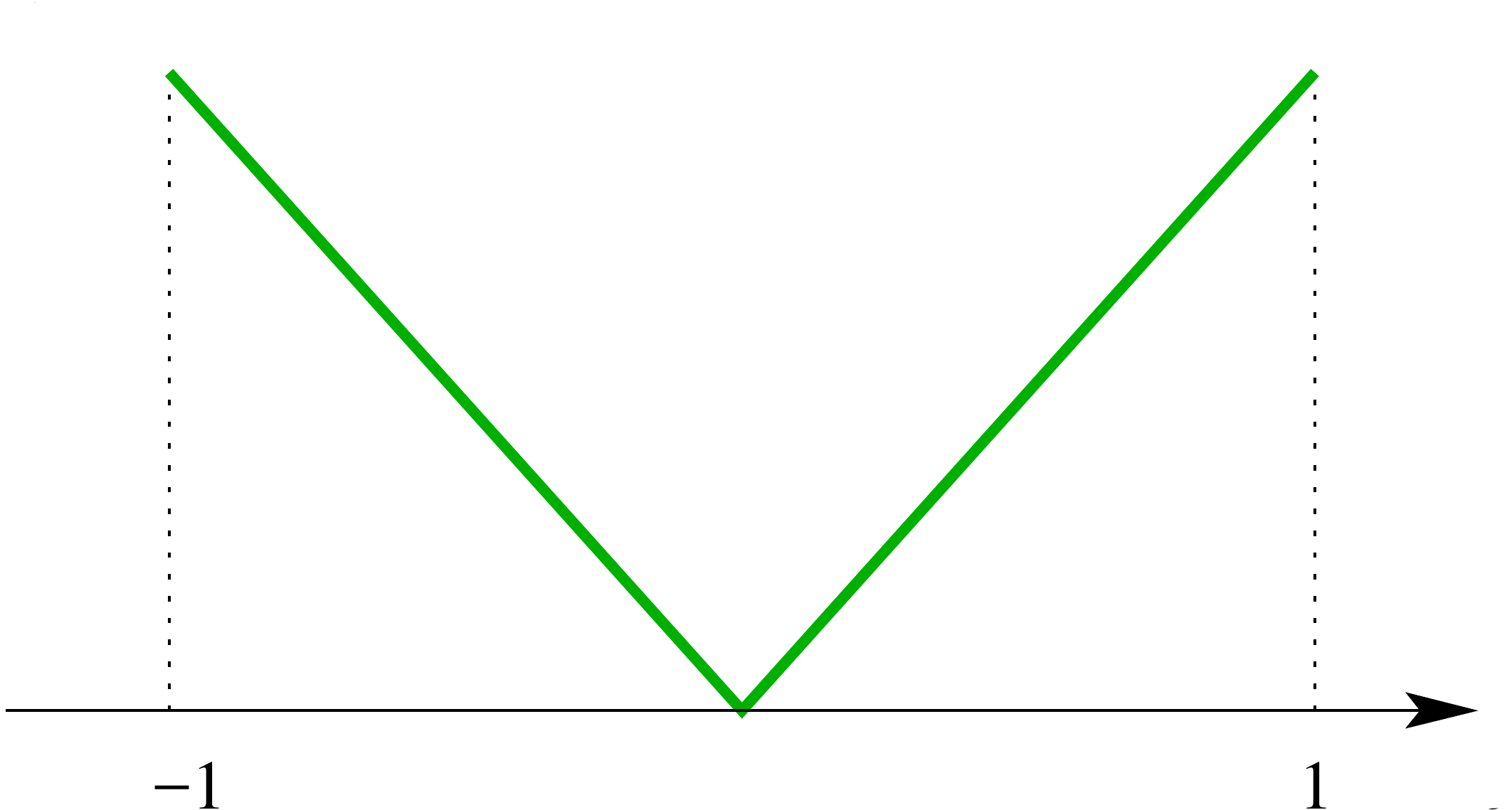}
\put(92,-1){\scriptsize $X$}
\put(2,50.5){\scriptsize $x$}
\end{overpic} 
}
\subfloat[\label{fig:AI_ss3}Type M3.]{
\begin{overpic}[width=4cm,tics=10]{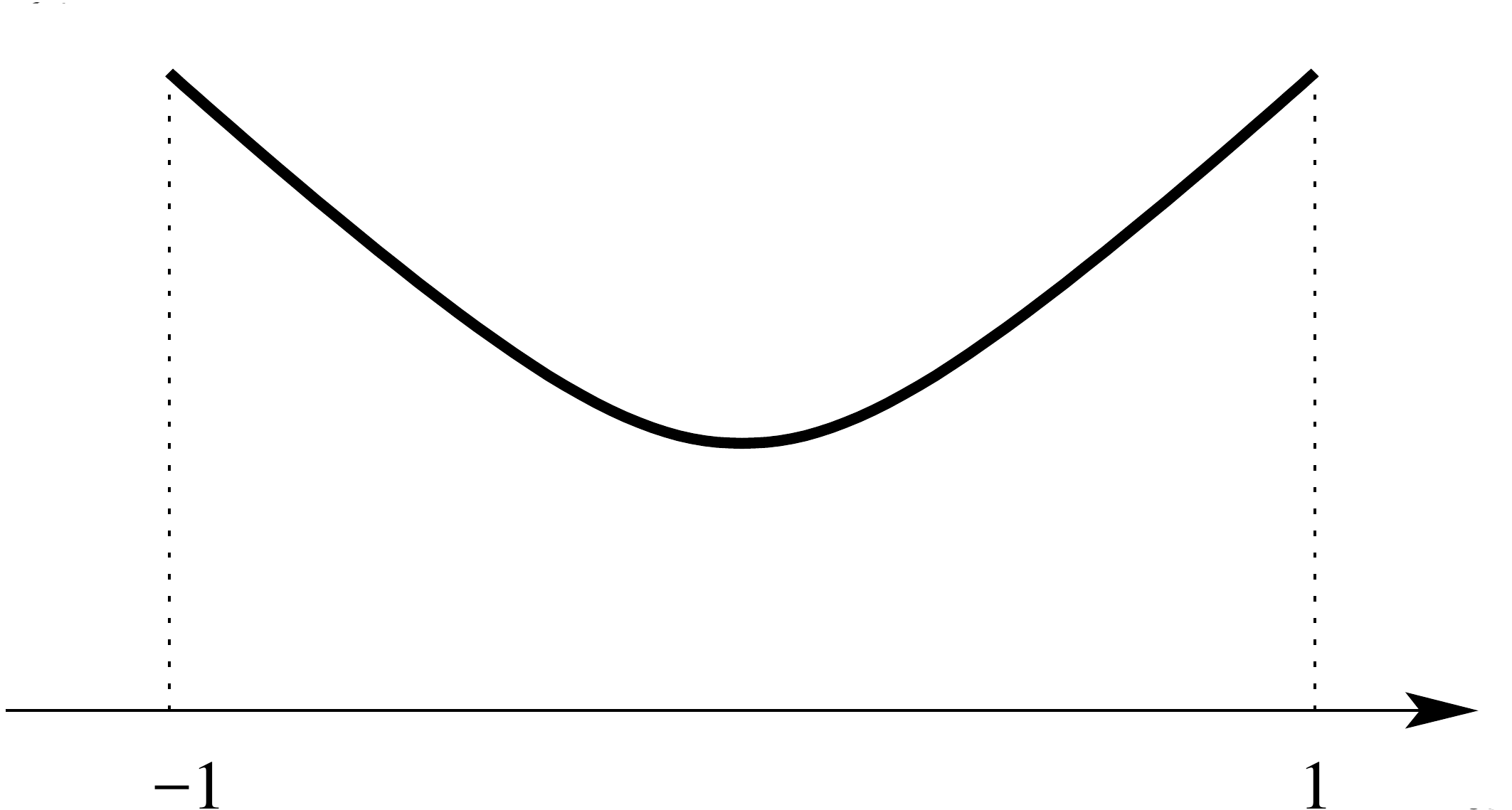}
\put(92,-1){\scriptsize $X$}
\put(2,50.5){\scriptsize $x$}
\end{overpic} 
}
\end{center}
\caption{Singular solutions to \eqref{eq:AI_op}.}
\label{fig:AI_soltyp}
\end{figure}

Solutions to \eqref{eq:AI_op} arise as perturbations of singular 
solutions obtained in the limit of $\varepsilon=0$. Such solutions
are constructed by analyzing the dynamics in charts $K_1$ and $K_2$ 
separately in the limit as $\varepsilon\to0$. In particular, solutions 
are constructed via two strategies:
\begin{description}
 \item[Strategy 1.] We consider two sets of boundary conditions, 
corresponding to suitable intervals of $y$-values that are defined 
at $\xi=-1$ and $\xi=1$, respectively. Flowing these two sets of 
boundary conditions forward and backward, respectively, we verify 
the transversality of the intersection of the two resulting manifolds 
at $\xi=0$. Each initial $y$-value $y_0$ for which these two manifolds 
intersect gives a solution to the boundary value problem \eqref{eq:AI_op}.
 \item[Strategy 2.] Since all such solutions are even, we can 
focus our attention on the $\xi$-interval $[-1,0]$, with boundary 
conditions $x(-1)=1$ and $y(0)=0$. The set of initial conditions 
at $\xi=-1$ and $x=1$, but with arbitrary initial $y$-value $y_0$, 
is then tracked forward up to the hyperplane $\{y=0\}$. The resulting 
manifold is parametrized by $x(y,\varepsilon,\delta,y_0)$ 
and $\xi(y,\varepsilon,\delta,y_0)$; the unique ``correct'' value 
$y_0(\varepsilon,\delta)$ corresponding to a solution to \eqref{eq:AI_op} 
is then obtained by solving $\xi(y_0,\varepsilon,\delta)=0$ under the 
constraint that $y(y_0,\varepsilon,\delta)=0$.
\end{description}
We distinguish three types of singular solutions to \eqref{eq:AI_op} 
(see Figure \ref{fig:AI_soltyp}):
\begin{description}
\item[\emph{Type M1.}] Solutions of type M1 (indicated in blue in 
the following figures) satisfy $x=0$ for $X \in I$, where $I$ is 
an interval centered at $X=0$. They occur in two subtypes: 
the ones corresponding to $\lambda=\mathcal{O}(\varepsilon)$ have 
constant finite slope $w$ outside of $I$, while the ones
corresponding to $\lambda=\mathcal{O}(1)$ vanish on $I=(-1,1)$.
\item[\emph{Type M2.}] Solutions of type M2 (indicated in green) 
are those of slope $y\equiv\mp1$. 
These solutions reach $\{x=0\}$ at one point only, namely at $X=0$. 
\item[\emph{Type M3.}] Solutions of type M3 (indicated in black) 
never reach $\{x=0\}$.
\end{description}

\begin{figure}[ht]
\centering
\begin{overpic}[scale=.3]{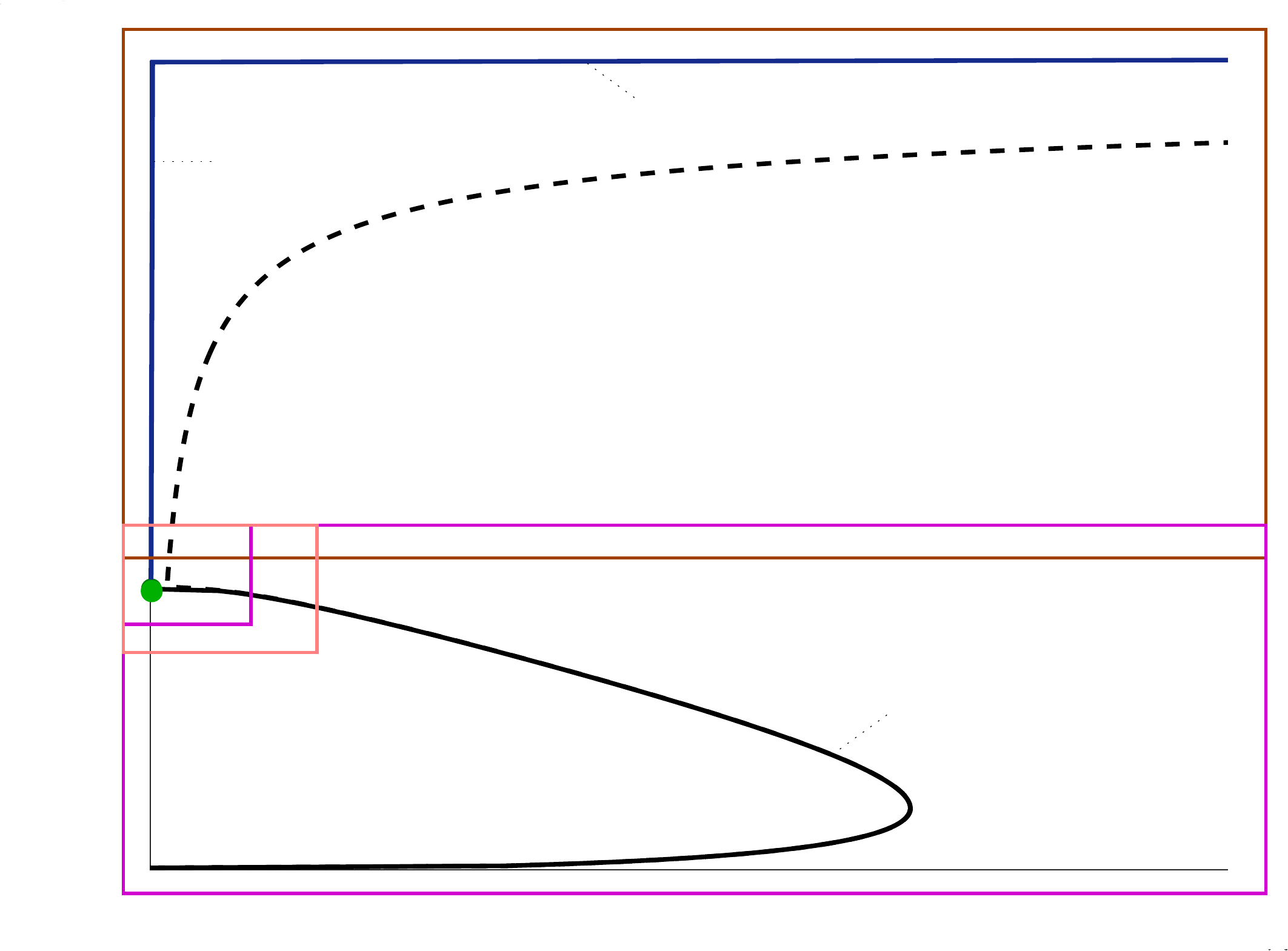}
\put(90,8.5){\scriptsize $\mathcal{R}_3$}
\put(70,18.5){\scriptsize $\mathcal{B}_3$}
\put(20,20){\scriptsize $\mathcal{R}_2$}
\put(90,35){\scriptsize $\mathcal{R}_1$}
\put(18,60){\scriptsize $\mathcal{B}_1$}
\put(50,63){\scriptsize $\mathcal{B}_2$}
\put(100,0){\footnotesize $\lambda$}
\put(-2,73){\footnotesize $\|u\|_2^2$}
\end{overpic} 
\caption{
Covering of the bifurcation diagram for \eqref{eq:AI_op} by 
regions $\mathcal{R}_1$ (brown), $\mathcal{R}_2$ (pink), and 
$\mathcal{R}_3$ (magenta). The branches of solutions to 
\eqref{eq:AI_op} for $\varepsilon=0.01$ (dotted curve) and 
$\varepsilon=0$ (solid curve) are displayed. For $\varepsilon=0$, 
the upper branch reduces to the union of a vertical part $\mathcal{B}_1$, 
corresponding to $\lambda = \mathcal{O}(\varepsilon)$, and a 
horizontal part $\mathcal{B}_2$ which corresponds to 
$\lambda = \mathcal{O}(1)$. The green dot at $B$ represents the 
singular solution of type M2 for $\lambda=0$. The black curve 
for type M3-solutions is labeled $\mathcal{B}_3$.
}
\label{fig:AI_bdreg}      
\end{figure}
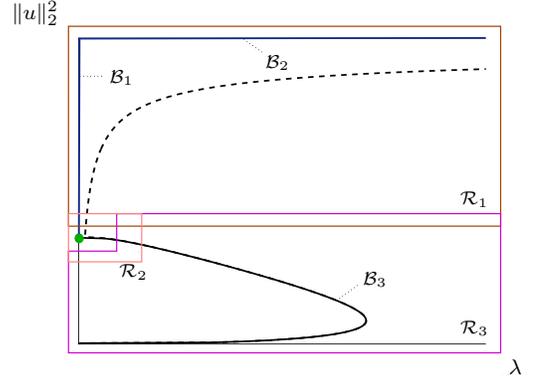

For $\varepsilon>0$, we divide the bifurcation diagram 
in $(\lambda, \Vert x \Vert_2^2)$, in terms of the original 
variable, into three overlapping regions, as shown 
in Figure \ref{fig:AI_bdreg}:
\begin{subequations}
\begin{align}
\mathcal{R}_1 &:=[0,1]\times\bigg[\frac23+\nu_1,2\bigg], 
\label{eq:R1} \\
\mathcal{R}_2 &:=[0,\varepsilon\lambda_2]\times
\bigg[\frac23-\nu_2,\frac23+\nu_2\bigg], \label{eq:R2} \\
\mathcal{R}_3 &:=[0,1]\times\bigg[0,\frac23+\nu_2\bigg]
\setminus[0,\varepsilon\lambda_3]
\times\bigg[\frac23-\nu_3,\frac23+\nu_2\bigg], 
\label{eq:R3}
\end{align}
\end{subequations}
with $\nu_2 > \nu_1 >0$, $\nu_2 > \nu_3 > 0$, and 
$\lambda_2 > \lambda_3 > 0$ large.

\begin{figure}[htbp]
 	\centering
 	\begin{tikzpicture}
[>=stealth',point/.style={circle,inner sep=0.01cm,fill=white,draw},
encircle/.style={circle,inner sep=0.07cm,draw},
x=4cm,y=4cm,declare function={f(\x) = {0.5*\x^(1/2)};}]
\draw[blue,thick,-,smooth,domain=0:1,samples=200,/pgf/fpu,/pgf/fpu/output
format=fixed] plot (\x, {f(\x)});

\draw[->,dashed,thick] (0,0) -> (1.1,0);
\draw[->,dashed,thick] (0,0) -> (0,1);
\draw[dashed,gray,thick] (0,0.9) -> (1,0.9);
\draw[dashed,gray,thick] (1,0.9) -> (1,0.5);
\draw[thick,blue] (0,0) -> (0,0.9);
\draw [blue,fill=blue] (0,0) circle [radius=0.01];

\node[] at (1,-0.05) {$\varepsilon_0$};
\node[] at (1.15,-0.05) {$\varepsilon$};
\node[] at (-0.04,1) {$\delta$};

\draw [black] (0.3,0.6) circle [radius=0.06];
\node at (0.3,0.6) {I};

\draw[->] (-0.12,0.4) -> (-0.02,0.5);
\draw [black,fill=white] (-0.12,0.4) circle [radius=0.06];
\node at (-0.12,0.4) {II};

\draw[->] (-0.1,-0.1) -> (-0.02,-0.02);
\draw [black,fill=white] (-0.12,-0.12) circle [radius=0.06];
\node at (-0.12,-0.12) {III};

\draw [black] (0.6,0.2) circle [radius=0.06];
\node at (0.6,0.2) {IV};

\end{tikzpicture}
 	\caption{Classification diagram of \eqref{eq:AI_op} with respect 
	to the property $\cP_{ss}$ in $\varepsilon \delta$-space. 
 	Regime (I) is bounded below by $\left\{ \delta=\sqrt{\varepsilon} 
	\right\}$ (blue curve) and above by $\left\{ \delta=\frac{2}{\sqrt{3}} 
	\right\}$ (dashed gray line). Regime (I): two singular solutions of 
	type M1 and M2 exist. Regime (II): singular solutions of type M1 exist. 
	Regime (III): singular solutions of type M2 and M3 exist. Regime 
	(IV): no singular solutions exist.}
 \label{fig:AI_epsdelta} 
 \end{figure}
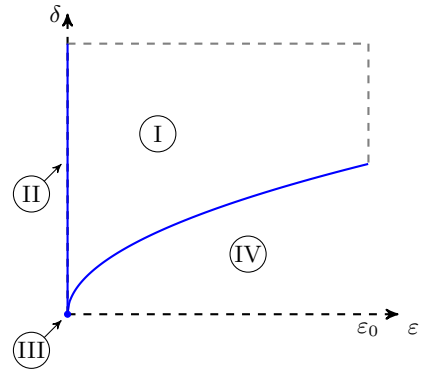

In our analysis, we consider $\lambda \in [0,1]$. In region 
$\mathcal{R}_3$, away from the point $B=\left(0,\frac23 \right)$, 
the perturbation with $\varepsilon$ is regular, and we consider 
$\lambda$ and $\delta$ as the two main parameters for our 
investigation. In regions $\mathcal{R}_1$ and $\mathcal{R}_2$, 
singular solutions exist only for $\lambda\geq\frac34 \varepsilon$ 
or, equivalently, for $\delta\leq\frac2{\sqrt3}$. Hence, in these 
regions, we need to take $\lambda\in\left[\frac34\varepsilon,1\right]$, 
i.e.~$\delta\in\left[\sqrt{\varepsilon},\frac2{\sqrt3}\right]$; 
see Figure \ref{fig:AI_epsdelta}. The two main parameters we consider 
in our proofs are here $\varepsilon$ and $\delta$. We define
\begin{equation}
    \cP_{ss} := \text{singular solutions of \eqref{eq:AI_sysd} exist.}
\end{equation}
In Regime $(I)$, such property is satisfied and singular 
solutions of type M1 and M2 exist, whereas in Regime $(IV)$ 
there are no singular solutions. Two special cases are represented 
by Regime $(II)$ (corresponding to $\mathcal{B}_1$), where 
singular solutions of type I exist, and Regime $(III)$ (corresponding to $\mathcal{B}_2 \cup \mathcal{B}_3$), where we 
recover singular solutions of type M1 and M3.

By definition, $\delta=0$ occurs only when $\varepsilon=0$. The corresponding, 
highly degenerate limit gives a singular orbit of type M1 with very singular 
structure, as shown in \ref{fig:AI_ss1bn}. Hence, the whole line 
$\mathcal{B}_2$ corresponds to that one singular solution. 

In $\mathcal{R}_1$, we construct singular solutions and show their 
persistence for \textcolor{black}{ $\varepsilon \in (0, \varepsilon_0)$ (with $\varepsilon_0$ small) using} Strategy 1 as follows. 
For a fixed $\lambda \in \left[ \frac{\varepsilon}{\delta_1^2},1 \right]$ 
with $0 < \delta_1 < \frac{2}{\sqrt{3}}$, i.e.~$\delta < \frac{2}{\sqrt{3}}$, 
the presence of a saddle-node equilibrium for the $(x_2,y_2)$-subsystem in 
chart $K_2$ at $(1,0)$, on which the reduced flow w.r.t. $\xi_2$ occurs, 
allows us to determine the unique, correct boundary value for $y$ at $\xi=-1$
 by following the stable manifold of such equilibrium, which does not depend 
on $\varepsilon$ and does therefore not change for $\varepsilon > 0$,
backwards until $\Sigma_2^{\rm in}=\Sigma_1^{\rm out}$, and then tracking 
the flow in chart $K_1$ backwards until $\xi_1=-1$. The intrinsic symmetry 
of the problem allows us to apply the same argument to the right part of 
the orbit, tracking the unstable manifold of the equilibrium $(x_2,y_2)=(1,0)$ 
and following the flow in chart $K_1$ until $\xi_1=1$. When $0 \leq \delta 
< \hat{\delta}$, the proof is analogous except for the fact that we must 
rescale $\delta=\sqrt{\varepsilon} \tilde{\delta}$ and obtain a slower 
reduced flow. The assumption that $\delta < \frac{2}{\sqrt3}$ ensures a 
non-trivial slow drift (i.e.~the portion of the orbit where $x=0$ does not 
reduce to a point), which allows us to apply the Exchange Lemma to infer 
persistence of solutions for $0 < \varepsilon \ll 1$. 

For $\frac{1}{\sqrt{\lambda_2}} \leq \delta \leq \delta_1$, i.e.~in $\mathcal{R}_2$, we show the existence of two unique type M1 and type M2 solutions which coincide when $\delta = \frac{2}{\sqrt{3}}$. 
The proof consists of two parts: we first consider a small 
neighborhood of $\delta_\ast=\frac{2}{\sqrt{3}}$, i.e.~of 
$\lambda=\frac34 \varepsilon$, where the saddle-node bifurcation 
occurs. We define a suitable bifurcation equation, which describes 
the transition from solutions which limit on type M1-solutions to 
those which limit on solutions of type M2. Such equation is constructed 
by imposing that $\xi_2(y_0,\varepsilon,\delta)=0$ when 
$y(y_0,\varepsilon,\delta)=0$, i.e.~using Strategy 2. Based on that 
equation, we infer the presence of the saddle-node bifurcation, and 
we calculate the expansion of the corresponding $\lambda$-value 
$\lambda_\ast$. This expansion presents logarithmic switchback terms 
due to both a resonance phenomenon in chart $K_1$ and the passage 
close to the saddle point $(x_2,y_2)=(1,0)$. In a second step, we 
consider the branch of solutions that limit on type M2-solutions for 
the remaining values of $\lambda$ in $\mathcal{R}_2$. That branch is 
then shown to connect to solutions that are covered by region 
$\mathcal{R}_3$, for which $\delta=0$. In that case, the type M2-solution 
constructed in $\mathcal{R}_2$ collapses onto the line $\{y_1=0\}$, 
which leads to singular dynamics in $K_1$. Since such singular nature 
is due to the $w$-rescaling introduced to obtain System \eqref{eq:AI_sysd}, 
this regime is better studied using System \eqref{eq:AI_sysl} and 
replacing $\varepsilon=\delta^2 \lambda$. Since this region contains 
a neighborhood of $(\delta,\lambda)=(0,0)$, we must perform an 
additional blow-up of $(u,\lambda)=(0,0)$ and split $\mathcal{R}_3$ into 
two sub-regions: for $\lambda \in [\tilde{\lambda},\lambda^\ast]$ with 
$\tilde{\lambda}>0$ and $\delta=0$, we can show the existence of a unique 
singular solution of type M3 which perturbs regularly when $0 < \delta \ll 1$ 
(in particular $\delta \leq \frac{1}{\sqrt{\lambda_3}}$ in $\mathcal{R}_3$). 
When $\lambda \in [0,\tilde{\lambda}]$, i.e.~when $\mathcal{R}_3$ and 
$\mathcal{R}_2$ overlap, we have a singular solution of type M2 as $\lambda 
\to 0$, and of type M3 as $\delta \to 0$.\medskip

In summary, even unfolding a rather innocent-looking PDE problem via 
spatial dynamics in one dimension leads to a highly interesting double limit
problem. In the next section, we continue this theme and consider a multi-component 
stationary PDE problem. 

\subsection{Fast Reaction Limits}
\label{ssec:fastreaction}

A variety of biological and ecological phenomena present different intrinsic time-scales, and typically some processes are faster than others. The singular limit, or fast reaction limit, expresses the fact that instantaneous dynamics is also included in the system. For instance, in a population, there can be a dichotomy of two groups, and switching between them may be possible. Compared to other interactions, the switch may seem \emph{instantaneous} and give rise to interesting effects such as an aggregation of individuals or a population density pressure \cite{hilhorst2009fast, brocchieri2020evolution}. Fast reaction limits have also been studied in other contexts, such as reversible and irreversible chemical reactions \cite{bothe2012instantaneous, bothe2003reaction}, bacteria proliferation \cite{hilhorst2007mathematical}, proteins localisation in stem cell division \cite{henneke2016fast}, but also to model the Neolithic spread of farmers in Europe \cite{eliavs2018well, eliavs2021singular}. 

In the context of predator--prey interactions, the expression of widely used functional responses can also come out of a systematic process in which one starts with a system of more than two equations with simple reaction terms and performs one \cite{metz2014dynamics, huisman1997formal, berardo2020interactions, lehtinen2019cyclic} or more limits \cite{geritz2012mechanistic, desvillettes2019non}.

We consider here the cross-diffusion system, known as Shigesada--Kawasaki--Teramoto (SKT) model~\cite{shigesada1979spatial}, proposed to account for stable inhomogeneous steady states exhibiting spatial segregation between two species competing for resources. We refer to~\cite{breden2019influence, kuehn2020numerical, soresina2021hopf} and references therein for more details. The system is given by
\be
\label{eq:CS_fullX}
\begin{array}{rcl}
\partial_t u-\Delta_x\left((d_1+d_{12} v) u\right) &=&  f(u,v)u,\\
\partial_t v-\Delta_x\left((d_2+d_{21} u) v\right) &=&  g(u,v)v,
\end{array}
\ee
endowed with initial conditions and homogeneous Neumann boundary conditions. 
The quantities~$u(t,x),\, v(t, x) \geq 0$ represent the population densities 
of two species at time $t$ and position $x$, confined on a bounded and 
connected domain~$\Omega \subset \mathbb{R}^N$. The movements of the 
individuals on the domain are described by non-linear cross-diffusion 
terms: the positive coefficients $d_1,\,d_2$ refer to the (standard) 
diffusion, while the non-negative cross-diffusion coefficients~$d_{12},\, 
d_{21}$ stand for competition pressure. The reaction terms describe 
the growth and the interaction of the two species, where
\be
\label{eq:CS_reactionpart}
\begin{array}{rcl}
f(u,v)&=&r_1-a_1 u-b_1 v,\\
g(u,v)&=&r_2-b_2 u-a_2 v,
\end{array}
\ee
with the non-negative coefficients~$r_i,\, a_i,\, b_i \,(i = 1,2)$ 
being the intrinsic growth, the intra-specific competition and the 
inter-specific competition rates. 

Model~\eqref{eq:CS_fullX} falls into the class of quasilinear 
parabolic systems for which even the existence problem of solutions 
is not trivial. When $d_{21}=0$ (triangular cross-diffusion system), 
it has been shown~\cite{iida2006diffusion, izuhara2008reaction} that 
the solutions of~\eqref{eq:CS_fullX} can be approximated in a finite 
time interval by those of a three-component reaction--diffusion system 
if the solutions are bounded and provided that a suitable parameter 
is small enough. The rigorous proof of the convergence of solutions 
of the three-component reaction--diffusion system towards the solutions 
of a triangular cross-diffusion system of two equations has been 
initially given in dimension~$N=1$~\cite{conforto2014rigorous}, 
and later generalized to a wider set of admissible reaction terms 
and in any dimension~\cite{desvillettes2015new}.

The convergence of the stationary steady states of the fast-reaction 
system towards the ones of the cross-diffusion system has been also 
investigated by looking at bifurcation diagrams with respect to 
different bifurcation 
parameters~\cite{izuhara2008reaction,kuehn2020numerical}. In particular,
it has been observed that the bifurcation structure of the fast-reaction
expands and converges as the time scale parameter becomes smaller, 
sometimes going through major qualitative changes. 

When~$d_{21}>0$, the full cross-diffusion system~\eqref{eq:CS_fullX} 
can be obtained, at least formally, as the singular-limit of a 
four-component fast-reaction system involving two small time scale 
parameters $\varepsilon,\, \delta$\textcolor{black}{, so the problem has a doubly singular perturbation structure}. In this case both species 
are split into quiet and active states, denoted by~$u_1,\, v_1$ 
and $u_2,\, v_2$ respectively. Hence, we have that~$u:=u_1+u_2,\,
v:=v_1+v_2$. The resulting reaction--diffusion system is
\be
\label{eq:CS_4e}\tag{{$\cX_{\textnormal{fr}}$}}
\begin{array}{rl}
\partial_t u_1-d_1\Delta_x u_1&\hspace{-0.2cm} = 
f(u,v)u_1+\frac{1}{\varepsilon}h(u_1,u_2,v),\\[0.2cm]
\partial_t u_1-\hat{d}_1\Delta_x u_2 & \hspace{-0.2cm} =
f(u,v)u_2-\frac{1}{\varepsilon}h(u_1,u_2,v),\\[0.2cm]
\partial_t v_1-d_2\Delta_x v_1&\hspace{-0.2cm}  =g(u,v)
v_1+\frac{1}{\delta}k(u,v_1,v_2),\\[0.2cm]
\partial_t v_2-\hat{d}_2\Delta_x v_2 &\hspace{-0.2cm}= 
g(u,v)v_2-\frac{1}{\delta}k(u,v_1,v_2),
\end{array}
\ee
together with initial conditions and homogeneous Neumann boundary 
conditions. Active states are supposed to have a larger diffusion 
coefficient than the corresponding quiet state. In particular, we 
assume that the diffusion coefficients of the active states are 
given by~$\hat{d}_1:=d_1+d_{12}M_2$ and~$\hat{d}_2:=d_2+d_{21}M_2$, 
where~$M_1,\,M_2$ are positive constants such that~$0\leq u(t,x)\leq 
M_1$,~$M_1\geq r_1/a_1$ and~$0\leq v(t,x)\leq M_2$,~$M_2\geq r_2/a_2$ 
in~$\mathbb{R}\times\Omega$. The functions~$h,\, k$ describing the 
switch between the states are
\be
\label{eq:CS_components}
\begin{array}{rcl}
h(u_1,u_2,v)&=&\left( 1-\dfrac{v}{M_2}\right)u_2-u_1\dfrac{v}{M_2},\\[0.25cm]
k(u,v_1,v_2)&=&\left( 1-\dfrac{u}{M_1}\right)v_2-v_1\dfrac{u}{M_1},
\end{array}
\ee
while the time scale parameters~$\varepsilon,\, \delta$ describe 
that the switch between the two different states happens much 
faster than the other processes. 

\begin{figure}
\begin{center}
\begin{tikzpicture}

\draw [black, ultra thick, <->] (0.5,4.5) -- (0.5,0.3) -- (5.3,0.3);
\node at (5.15,0) {$\varepsilon$};
\node at (0.2,4.15) {$\delta$};
\draw [black, fill=gray!30!white] (0.5,0.3) circle [radius=0.1];

\node at (-0.85,-0.2) {2 eqs.};
\node at (-0.85,-0.55) {(full cross-d)};

\node at (-0.85,3.25) {3 eqs.};
\node at (-0.85,2.75) {(1st cross-d, $\delta$)};

\node at (3.707,-0.2) {3 eqs.};
\node at (3.707,-0.55) {(3rd cross-d, $\varepsilon$)};

\node at (3.707,3.25) {4 eqs.};
\node at (3.707,2.75) {$\varepsilon, \, \delta$};

\draw[gray!50!white,thick,->] (1.8,-0.3) -- (1.3,-0.3);
\node at (1.55,-0.55) {{\color{gray!80!white}$\varepsilon\to 0$}};

\draw[gray!50!white,thick,->] (2,2.75) -- (1.5,2.75);
\node at (1.75,3.25) {{\color{gray!80!white}$\varepsilon\to 0$}};

\draw[gray!50!white,thick,->] (-0.55,1.8) -- (-0.55,1.3);
\node at (-1.3,1.55) {{\color{gray!80!white}$\delta\to 0$}};

\draw[gray!50!white,thick,->] (3.75,1.8) -- (3.75,1.3);
\node at (4.5,1.55) {{\color{gray!80!white}$\delta\to 0$}};
\end{tikzpicture}
\end{center}
\caption{Schematic representation of the systems of PDEs in 
the~$\varepsilon \delta$-plane. 
The first quadrant corresponds to the four-equation system~\eqref{eq:CS_4e}. 
The~$\varepsilon$-axis corresponds to the reduced system with three 
equations for~$u_1,\, u_2,\, v$, being the last equation with 
cross-diffusion. The~$\delta$-axis corresponds to the reduced system 
with 3 equations for~$u, \, v_1,\, v_2$, being the first equation with 
cross-diffusion. Finally, the origin corresponds to 
system~\eqref{eq:CS_fullX} with two cross-diffusion equations.
}
\label{fig:CS_parplane}
\end{figure}
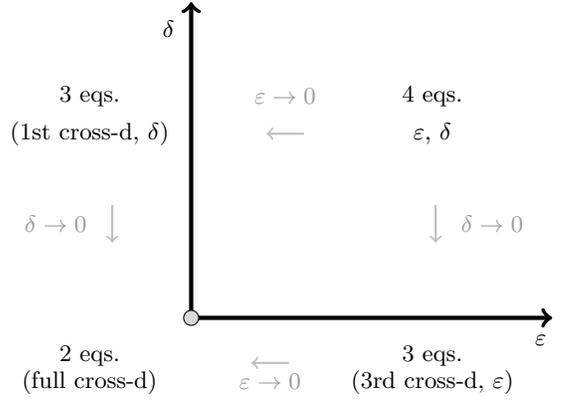

At a formal level, when~$\varepsilon\to0$, system~\eqref{eq:CS_4e} 
reduces to an intermediate three-component reaction--cross-diffusion 
system in the variables~$u,\, v_1,\, v_2$. The equation for~$u$ represents 
cross-diffusion, while the other time scale parameter~$\delta$ is 
still present in the equations for~$v_1,\, v_2$. Letting~$\delta\to0$, 
the intermediate three-component system reduces to the full cross-diffusion 
system. The same considerations hold if we let~$\delta\to 0$ first, 
and~$\varepsilon\to0$. In the time scale parameter plane 
(Figure~\ref{fig:CS_parplane}), the first quadrant corresponds to the 
four-equation system~\eqref{eq:CS_4e}. The~$\varepsilon$-axis corresponds 
to the reduced system with three equations for~$u_1,\, u_2,\, v$, being 
the last equation with cross-diffusion. The~$\delta$-axis corresponds 
to the reduced system with three equations for~$u, \, v_1,\, v_2$, being 
the first equation with cross-diffusion. Finally, the origin corresponds 
to system~\eqref{eq:CS_fullX} with two cross-diffusion equations.

To the best of the authors' knowledge, there are currently no rigorous results of 
convergence of solutions of the four-component reaction--diffusion systems 
to the solutions of the full cross-diffusion system. From a numerical 
point of view, despite a greater number of equations, the structure of 
system~\eqref{eq:CS_4e} is simpler than the cross-diffusion 
system~\eqref{eq:CS_fullX}, since it presents standard diffusion terms. 
For suitable small values of the time scale parameters~$\varepsilon, \, 
\delta$ that leads to a ``good'' approximation of the cross-diffusion 
system~\eqref{eq:CS_fullX}, the four-component fast-reaction system tends 
to be more tractable. In order to establish how accurate is the 
approximation, we look at the bifurcation structure of stationary 
solutions when~$\varepsilon, \, \delta$ become small.

On the one hand, system~\eqref{eq:CS_fullX} admits the homogeneous 
coexistence state~$(u_*,v_*)$ where
$$u_*=\dfrac{r_1a_2-r_2b_1}{a_1a_2-b_1b_2},\quad 
v_*=\dfrac{r_2a_1-r_1b_2}{a_1a_2-b_1b_2},$$
which is positive for suitable parameter values 
(see~\cite{breden2019influence, kuehn2020numerical}). It is known 
that the homogeneous solution undergoes some bifurcations under 
parameter variation, and branches of non-homogeneous solutions 
originate at these bifurcation points which correspond to different 
spatial distributions (patterns) of the species on the domain. 

\begin{figure}[!ht]
\begin{center}
\subfloat[\label{CS_bif_1}$\varepsilon=5\cdot 10^{-2},\, \delta=\cdot 10^{-1}$]{
\begin{overpic}[width=4cm,tics=10]{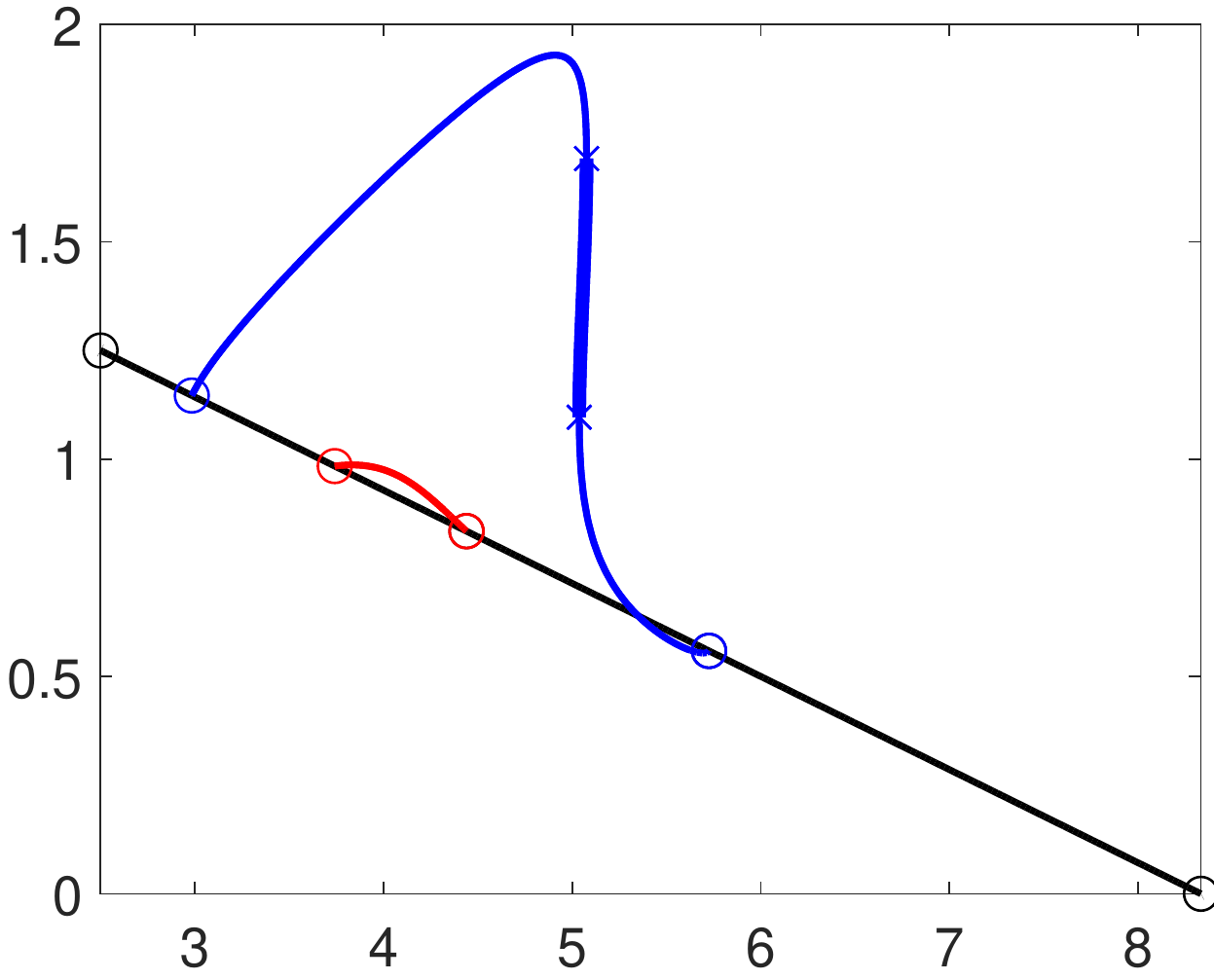}
\end{overpic} 
}
\subfloat[\label{CS_bif_2}$\varepsilon=5\cdot 10^{-2},\, \delta=10^{-2}$]{
\begin{overpic}[width=4cm,tics=10]{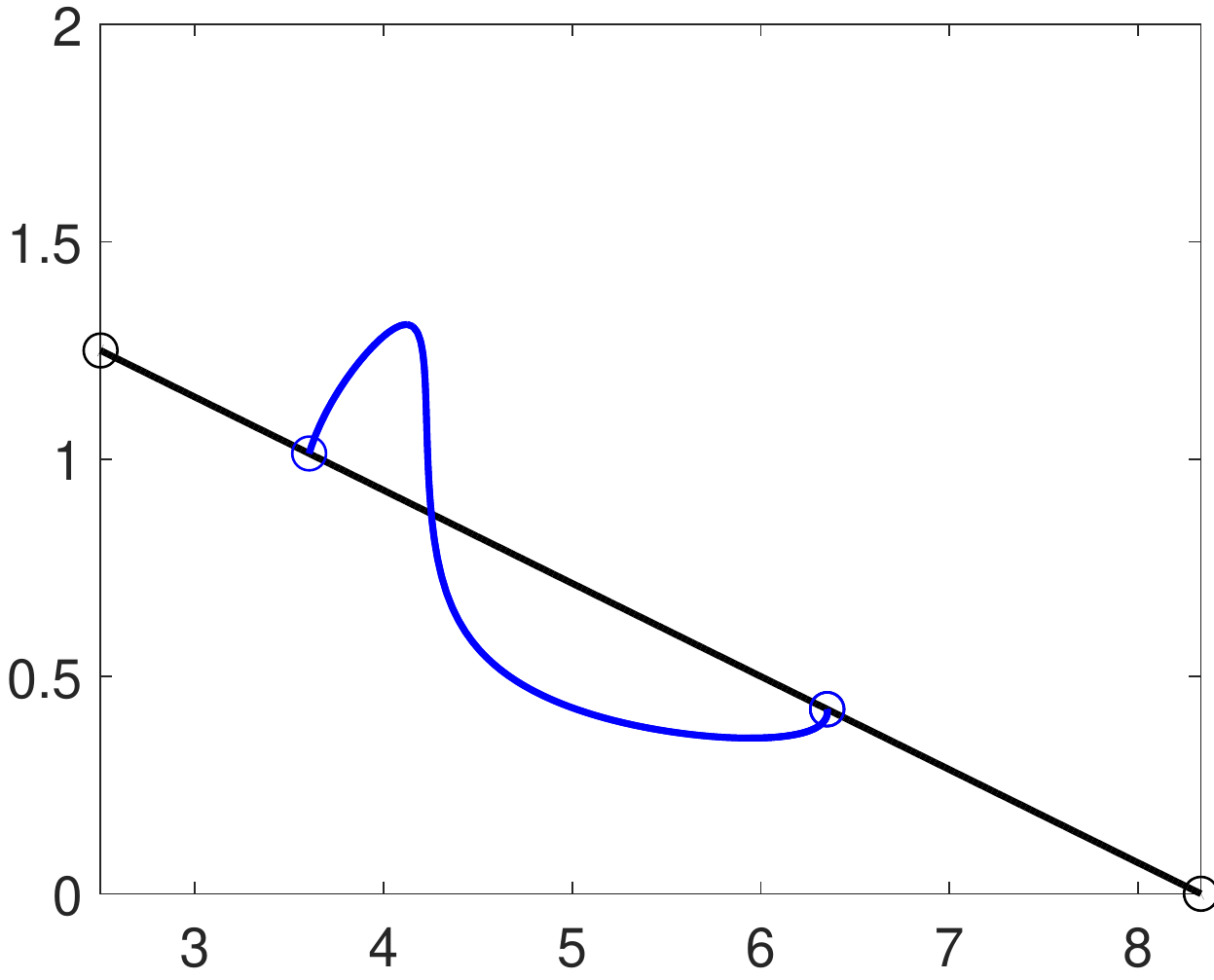}
\end{overpic} 
}\\
\subfloat[\label{CS_bif_3}$\varepsilon=5\cdot 10^{-2},\, \delta=5\cdot 10^{-3}$]{
\begin{overpic}[width=4cm,tics=10]{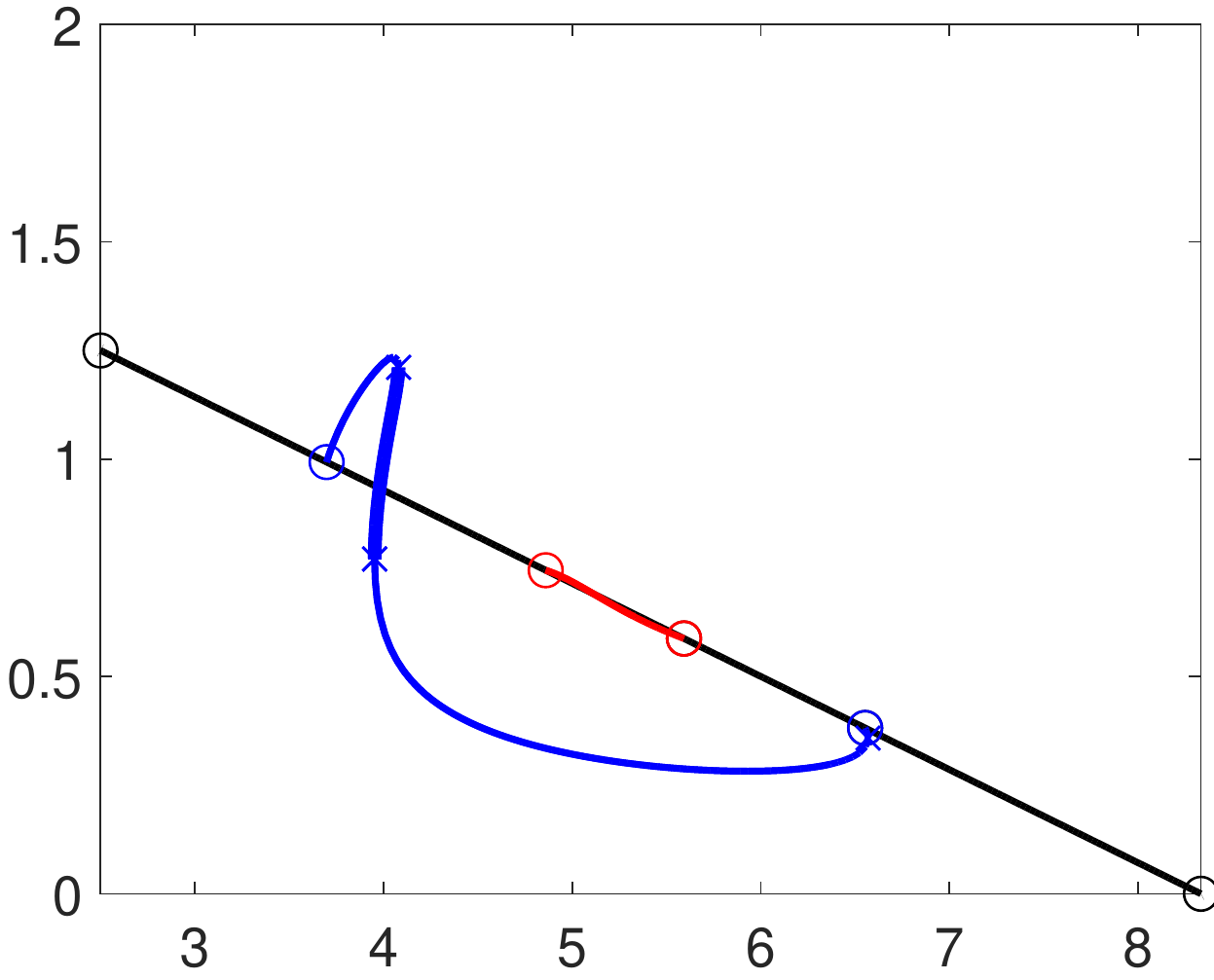}
\end{overpic} 
}
\subfloat[\label{CS_bif_4}$\varepsilon=10^{-3},\, \delta=5\cdot 10^{-3}$]{
\begin{overpic}[width=4cm,tics=10]{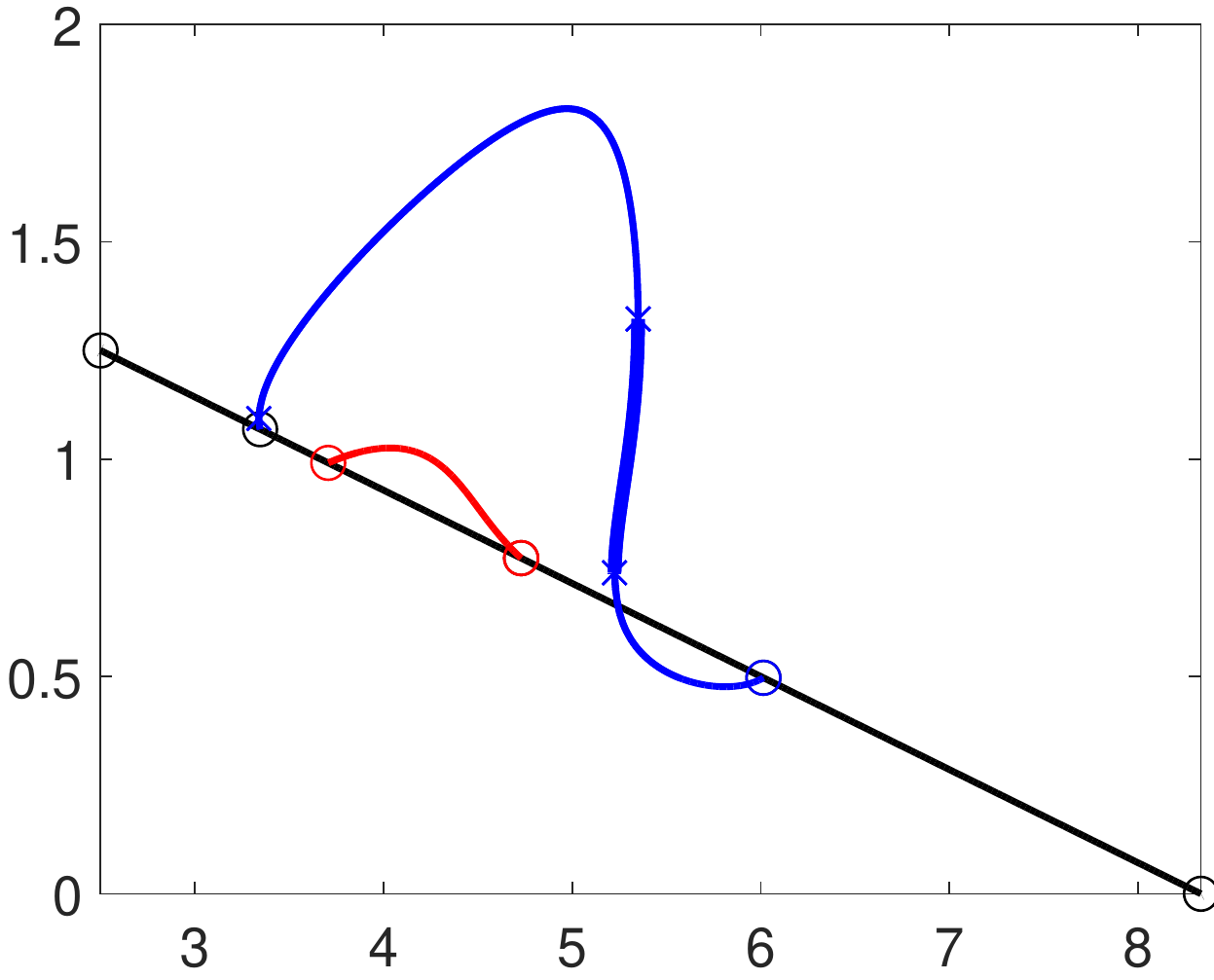}
\end{overpic} 
}\\
\subfloat[\label{CS_bif_5}$\varepsilon=10^{-5},\, \delta=10^{-5}$]{
\begin{overpic}[width=4cm,tics=10]{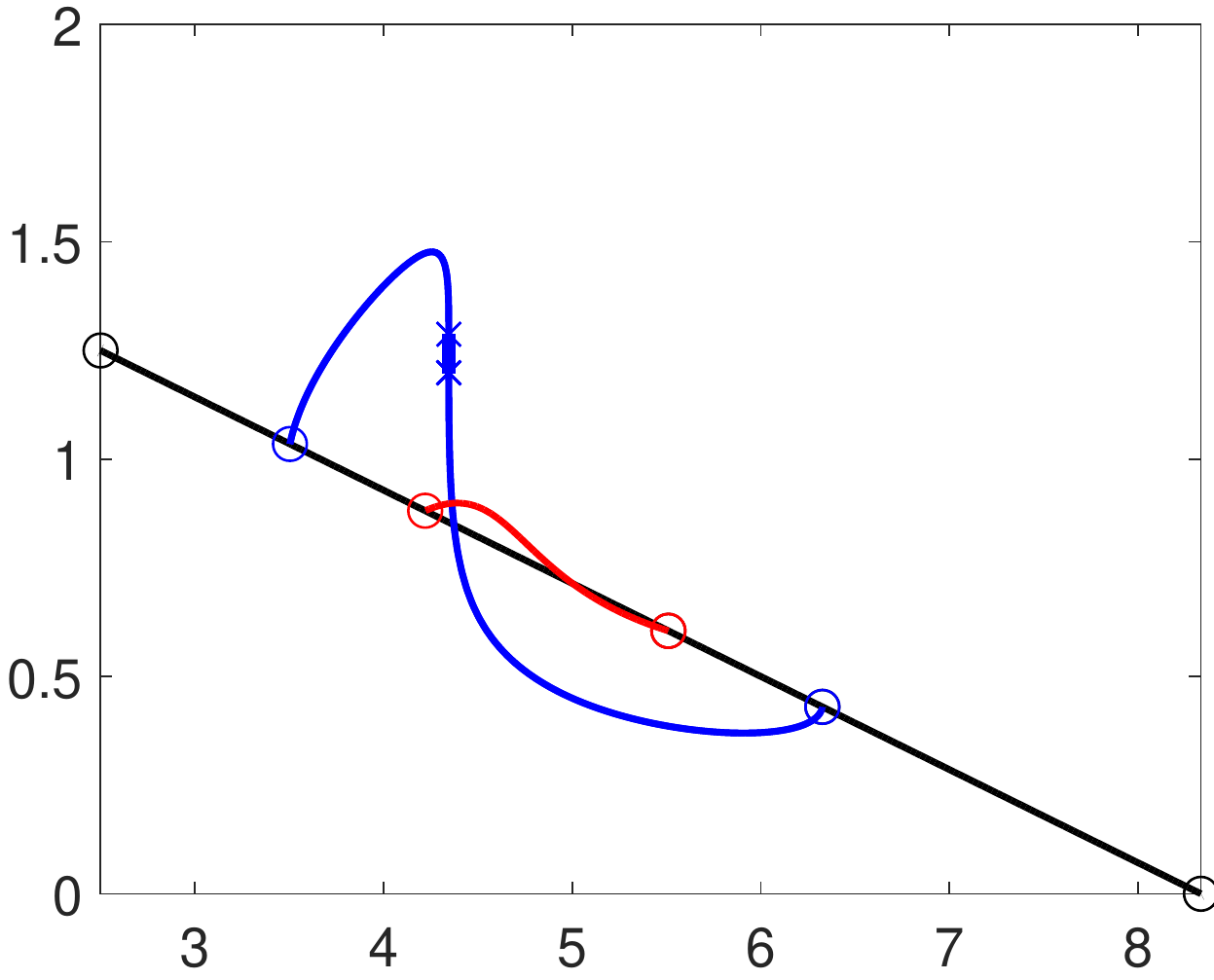}
\end{overpic} 
}
\subfloat[\label{CS_bif_6}cross-diffusion]{
\begin{overpic}[width=4cm,tics=10]{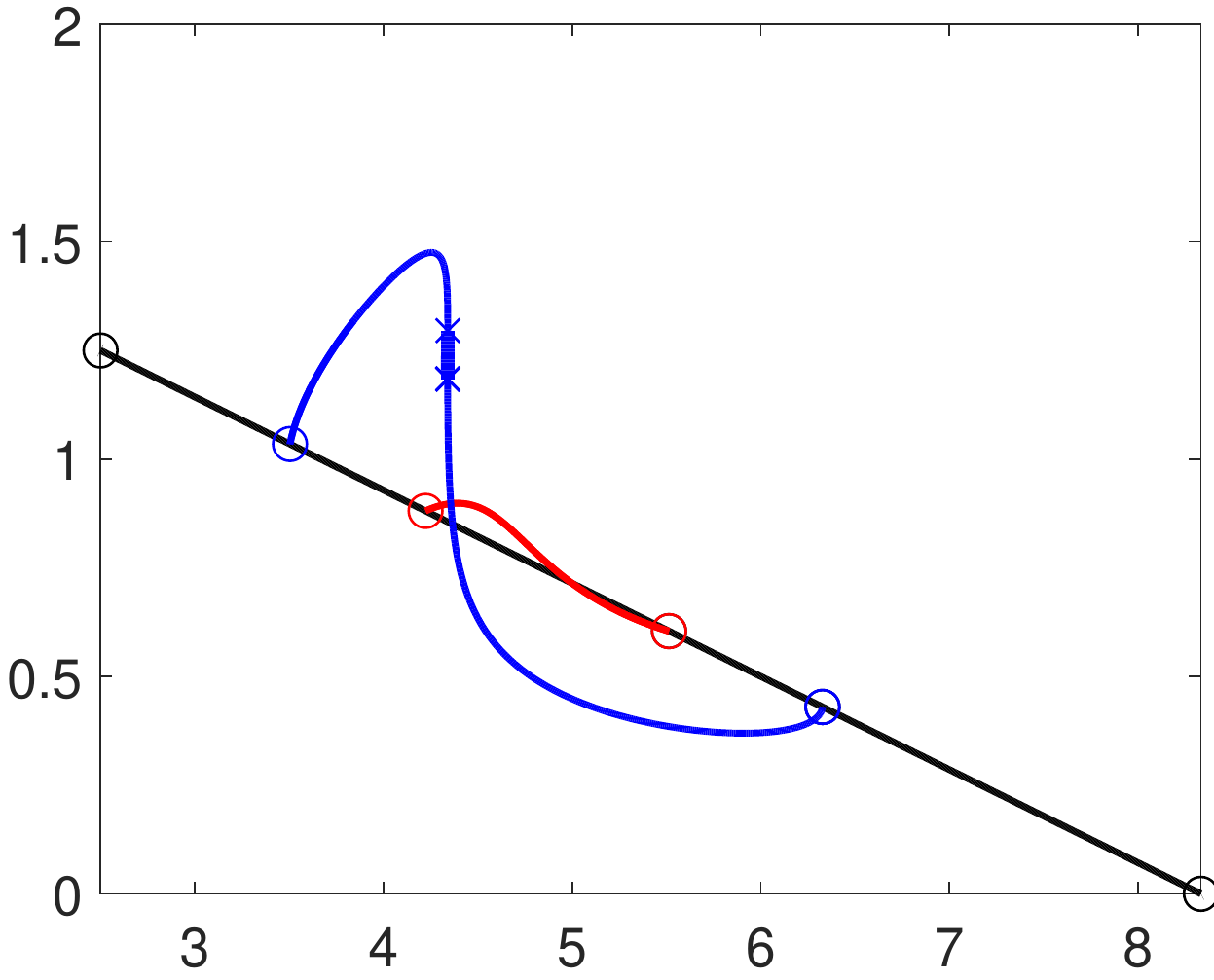}
\end{overpic} 
}
\end{center}
\caption{Bifurcation diagrams with respect to the bifurcation 
parameter $r_1$ and the parameter set in Table \ref{tab:CS_param} 
corresponding to different values of $\varepsilon$ and $\delta$. 
The black line corresponds to the homogeneous branch, while blue 
and red lines denote the bifurcating branches of non-homogeneous 
solutions. Bifurcation points are marked by circles.}
\label{fig:CS_bifdiag}
\end{figure}

On the other hand, also system~\eqref{eq:CS_4e} admits the homogeneous 
coexistence state~$(u_{1*},u_{2*},v_{1*},v_{2*})$, given by 
\benn
\begin{array}{rr}
u_{1*}=u_*\left(1-\dfrac{v_*}{M_2}\right),\qquad& u_{2*}=u_*\dfrac{v_*}{M_2},\\[0.2cm]
v_{1*}=v_*\left(1-\dfrac{u_*}{M_1}\right),\qquad& v_{2*}=v_*\dfrac{u_*}{M_1}.
\end{array}
\eenn
The homogeneous coexistence state turns out to be independent 
of the parameters~$\varepsilon,\,\delta$. However, the number and 
the position of the bifurcation points on the homogeneous branch, 
and hence the global bifurcation structure, changes depending on 
the time scale parameters. Then, we say that the cross-diffusion 
system~\eqref{eq:CS_fullX} and the four-component fast-reaction 
system~\eqref{eq:CS_4e} are equivalent if they have the same property
\benn
\begin{array}{rcl}
\cP_{fr}&:=&\text{number of bifurcation points}\\
&& \text{on the homogeneous branch w.r.t.} \\
&& \text{the bifurcation parameter}.
\end{array}
\eenn
In the following, we select a set of parameters already used 
in~\cite{breden2019influence} and reported in Table~\ref{tab:CS_param}. 
It corresponds to the strong competition case~$a_1a_2-b_1b_2<0$, in 
which the homogeneous coexistence state is unstable in absence of 
diffusion. However, stable non-homogeneous solutions arise on branches 
originating from bifurcation points on the homogeneous branch. In 
Figures~\ref{CS_bif_1}--\ref{CS_bif_5} we show different bifurcation 
diagrams obtained for smaller values of the parameters~$\varepsilon,\,\delta$, 
considering~$r_1$ as bifurcation parameter and fixing the other 
parameter values, while Figure \ref{CS_bif_6} corresponds to the 
non-triangular cross-diffusion system \eqref{eq:CS_fullX}.

\begin{figure}
\begin{center}
\begin{tikzpicture}

\draw [blue, ultra thick, <->] (0,5) -- (0,0) -- (5,0);
\node at (4.85,-0.3) {$\varepsilon$};
\node at (-0.3,4.85) {$\delta$};

\draw [blue, ultra thick] (0.2,0.2) arc (310:347:8cm);
\draw [blue, ultra thick] (0.2,0.2) arc (120:93:10cm);

\draw [black] (0.75,2) circle [radius=0.25];
\node at (0.75,2) {II};
\draw [black] (3,0.5) circle [radius=0.25];
\node at (3,0.5) {II};
\draw [black] (3,2) circle [radius=0.25];
\node at (3,2) {I};

\draw [gray, fill=gray] (2,3) circle [radius=0.05];
\node at (1.7,3) {{\color{gray}(a)}};
\draw [gray, fill=gray] (2,1.6) circle [radius=0.05];
\node at (2.3,1.6) {{\color{gray}(b)}};
\draw [gray, fill=gray] (2,0.6) circle [radius=0.05];
\node at (1.7,0.6) {{\color{gray}(c)}};
\draw [gray, fill=gray] (0.4,0.6) circle [radius=0.05];
\node at (0.4,0.8) {{\color{gray}(d)}};
\draw [gray, fill=gray] (0.1,0.1) circle [radius=0.05];
\node at (0.4,0.1) {{\color{gray}(e)}};
\draw [gray, fill=gray] (0,0) circle [radius=0.05];
\node at (-0.15,-0.2) {{\color{gray}(f)}};
\end{tikzpicture}
\end{center}
\caption{Qualitative classification diagram of system~\eqref{eq:CS_4e} 
with respect to the property~$\cP_{fr}$ in the~$\varepsilon \delta$-plane. 
Region I: two bifurcation points. Region II: four bifurcation points. 
Grey points correspond to the bifurcation diagrams in Figure \ref{fig:CS_bifdiag}.}
\label{fig:CS_parplane_zones_strong}
\end{figure}
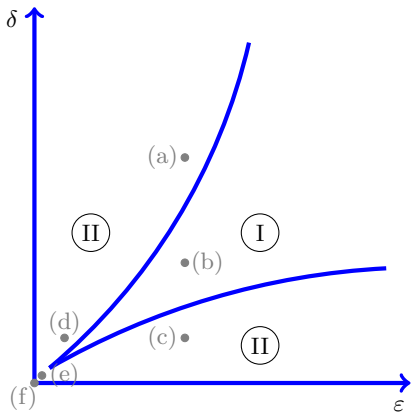

As shown in Figure~\ref{fig:CS_bifdiag}, considering the parameter 
set in Table~\ref{tab:CS_param}, we have that~$\cP_{fr}=2,4$. 
In Figure~\ref{fig:CS_parplane_zones_strong} the qualitative 
classification diagram of system~\eqref{eq:CS_4e} with respect 
to the property~$\cP_{fr}$ in the~$\varepsilon \delta$-plane is 
shown. The~$\varepsilon \delta$-plane can be split into two 
regions.  
Note also that in general the~$\varepsilon\delta$-diagram is not 
symmetric with respect to the diagonal~$\varepsilon=\delta$, 
but the intersections of the separation curves with the axis 
depend on the parameter set, in particular on the cross-diffusion 
coefficients. With different parameter sets, mainly with smaller 
standard diffusion coefficients $d$, one can obtain more bifurcation 
points on the homogeneous branch, and more zones in the~$\varepsilon 
\delta$-plane, but its structure remains qualitatively similar to 
Figure~\ref{fig:CS_parplane_zones_strong}. 

\begin{table}
\centering
\begin{tabular}{cccccccccc}
\hline\\[-0.3cm]
$r_2$&$a_1$&$a_2$&$b_1$&$b_2$&$d$&$d_{12}$&$d_{21}$&$M_1$&$M_2$\\
\hline\\[-0.3cm]
5&2&3&5&4&0.03&3&3&$5$&$2$\\
\hline
\end{tabular}
\caption{Set of parameter values relevant 
to~\ref{fig:CS_parplane_zones_strong}. The set~$r_i,\;a_i,\;b_i,\; (i=1,2)$ 
corresponds to the strong-competition case~($a_1a_2-b_1b_2<0$), 
namely the homogeneous coexistence state is unstable for the reaction part.}
\label{tab:CS_param}
\end{table}

The same study can be performed for other fast-reaction systems 
with multiple time scales and their cross-diffusion 
limits~\cite{conforto2018reaction,desvillettes2019non}.

\subsection{Coupled Oscillators}
\label{ssec:networks}

As the last example, we proceed to systems on networks. As discussed above, the presence of multiple time scales can lead to oscillations that are relevant in a variety of physical contexts; whether it is simple relaxation oscillations~\cite{Ginoux2012}, mixed mode oscillations~\cite{Desrochesetal}, or other examples of oscillatory deterministic dynamics discussed in Sections~\ref{ssec:mts} and~\ref{ssec:shear}.
However, it is not only the oscillations themselves but also the \emph{interaction} between different oscillatory processes that play an important role in many physical systems: These range from Huygens' synchronizing clocks~\cite{Huygens1888} to coupled oscillatory dynamics in neuroscience~\cite{Hoppensteadt1997,Ashwin2015}. From a mathematical perspective, such systems can be understood as networks of coupled oscillators. In isolation, each node oscillator has state~$z\in\R^d$ whose evolution is determined by a smooth ODE 
\begin{equation}
\label{eq:Osc}
z' := \frac{\txtd}{\txtd t} z = F(z)
\end{equation}
that gives rise to an asymptotically stable limit cycle~$\gamma\subset\R^d$. 
In a network, nodes interact non-trivially if there is an edge between the two nodes. Despite the dynamics of each node being fairly simple, the network dynamics, namely the dynamics of joint state of all nodes in the 
network, can be rich. 
While synchronization is probably one of the 
best understood dynamical phenomena in networks of coupled oscillators~\cite{Strogatz2004,Pikovsky2003,Bick2018c}, even networks consisting of just a few fully symmetric nodes can give rise to complicated dynamics~\cite{Bick2011}.
The network dynamics depend on both the intrinsic dynamical properties of each node and the network interactions. 
Here we will consider networks of weakly coupled relaxation oscillators, which have two small parameters: The time scale separation~$\delta$ as an intrinsic property of the oscillators themselves and the coupling constant~$\eps$ that is small by the assumption of weak coupling.

A network of~$N$ identical all-to-all coupled oscillators consists of~$N$ copies of~\eqref{eq:Osc} whose states~$z_k\in\R^d$, 
$k\in\{1, \dotsc N\}$, evolve according to
\begin{align}\label{eq:OscNet}
 z_k' &= F(z_k) + \frac{\eta}{N} \sum_{j=1}^N H(z_j, z_k),
\end{align}
where~$H$ is a smooth interaction function and~$\eta$ the 
coupling strength. If the coupling is weak, then the dynamics 
of this system on~$\R^{Nd}$ can be reduced to a lower-dimensional 
system~\cite{Ashwin1992}: If $\eta=0$ then~\eqref{eq:OscNet} has 
a normally hyperbolic invariant torus~$\gamma^N$ which persists 
for small coupling~\cite{Hoppensteadt1997}. Specifically, there 
exists an $\eta_0>0$ such that for any $\eta<\eta_0$ the 
system~\eqref{eq:OscNet} has an attracting normally hyperbolic 
invariant torus~$\mathbb{T}$ as a perturbation of~$\gamma^N\subset\R^{Nd}$.
In the following assume that~$\eta_0$ is maximal with this property; 
note that, depending on~$H$, this may allow for~$\eps_0=\infty$, 
for example, for trivial coupling $H=0$.
The dynamics of~\eqref{eq:OscNet} reduce to the interaction of~$N$ 
circular phase variables that evolve on~$\mathbb{T}$, a phase oscillator 
network. The dynamics on the invariant torus are typically referred 
to as a phase reduction of~\eqref{eq:OscNet}; cf.~\cite{Nakao2015,Pietras2019} 
for more details on how to compute these.

While a phase reduction is possible for any smooth oscillator, in 
many contexts the oscillators have particular properties. Relaxation 
oscillators are characterized by two time scales leading to a combination 
of slow quasi-static and fast transitions. The most famous examples 
include the van der Pol oscillator~\cite{vanderPol1} and FitzHugh--Nagumo 
oscillator~\cite{FitzHugh1961,Nagumo1962}. Consider a planar 
system~\eqref{eq:Osc} with state $z = (x,y)$ that evolves according to
\begin{subequations}
\label{eq:OscSF}
\begin{align}
\eps x' &= f(x,y)\\
y' &= g(x,y)
\end{align}
\end{subequations}
where $f,g:\R^2\to\R$ are smooth and~$\eps>0$ is the time scale separation 
of the fast variable~$x$ and the slow variable~$y$. Now assume 
that~\eqref{eq:OscSF} gives rise to a family of relaxation oscillators, 
that is there is a family of asymptotically stable limit 
cycles~$\gamma_\delta\subset\R^2$ that converge in the limit $\eps\to 0$ 
to a union of orbit segments consisting of part of the critical manifold 
$\{(x,y)\mid f(x,y)=0\}$ and line segments that correspond to the fast transitions.

In a series of papers~\cite{Somers1993,Kopell1995}, Somers and Kopell developed a theory to explain rapid synchronization in networks of coupled 
relaxation oscillators motivated by computational neuroscience. Write $z_k = (x_k, y_k)$ for the state of 
oscillator~$k$ which evolves according to~\eqref{eq:OscSF} when uncoupled. 
The networks analyzed in~\cite{Somers1993,Kopell1995} include systems of the 
form
\be
\label{eq:OscSFNet} \tag{{$\cX_{\textnormal{net}}$}}
\begin{array}{rcl}
\eps x_k' &=& f(x_k,y_k) + \frac{\delta}{N} \sum_{j=1}^N h(x_j,x_k),\\
 y_k' &=& g(x_k,y_k),
\end{array}
\ee
for $k\in\{1, \dotsc, N\}$ and coupling function~$h$ without 
specific assumptions on the coupling strength~$\eps$. \textcolor{black}{Note that \eqref{eq:OscSFNet} is clearly singularly perturbed in $\varepsilon$ due to its fast-slow structure, while $\delta=0$ yields a singular limit since we go from a networked dynamical system to an uncoupled case without network structure.}
The analysis~\cite{Somers1993,Kopell1995} considers the singular limit $\eps\to 0$ for network interactions such that the input from one node to the other is constant on each segment of the critical manifold and evaluates the ``compression'' of time it takes a singular trajectory to traverse segments of the critical manifold.
But even in the context of coupled neurons, other forms of network interactions~$h$---such as pulsatile coupling---are relevant.

\begin{figure}
	\centering
	\begin{overpic}[width=0.2\textwidth]{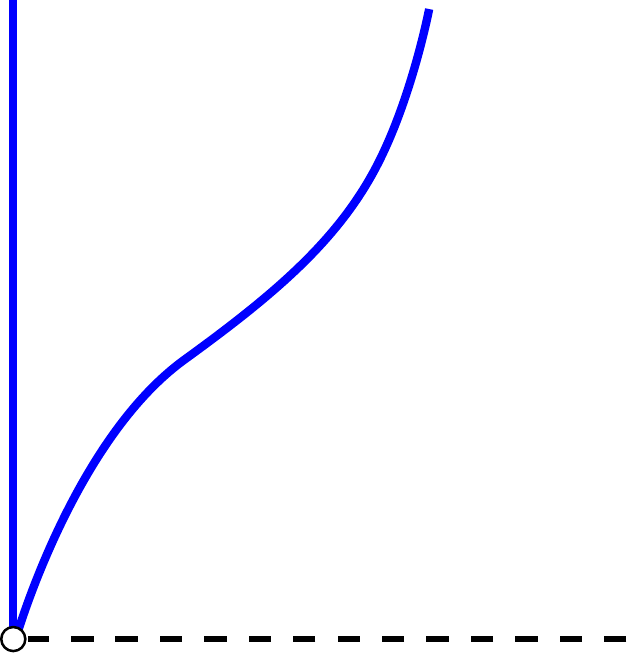}
	\put(90,-5){$\delta$}
	\put(-5,90){$\eps$}
	\put(20,72){\circle{13} \makebox(-32,0){I}}
	\put(74,32){\circle{13} \makebox(-32,0){II}}
	\end{overpic}
\caption{\label{fig:CoupledOsc}
\textcolor{black}{Sketch of a typical classification diagram expected for a phase reduction of a network of relaxation oscillators.}
The property~$\cP_{\textnormal{net}}$ 
divides the parameter space for~\eqref{eq:OscSFNet} into a region~(I), 
where a phase reduction is possible and a region~(II), where we expect a 
torus breakdown for a generic coupling function~$h$. The line 
dividing the region is given by~$\delta_0(\eps)$.}
\end{figure}
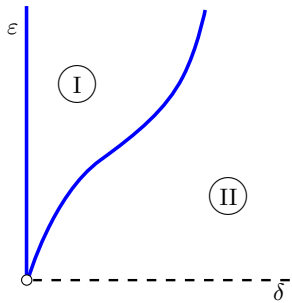

If both the time scale separation~$\eps$ for the relaxation 
oscillator and the coupling strength~$\delta$ are small, then the 
qualitative dynamics of~\eqref{eq:OscSFNet} can be understood in 
terms of the unified framework above.
\textcolor{black}{Consider the property
\benn
\begin{array}{rcl}
\cP_{\textnormal{net}}&:=&\text{a phase reduction is possible}.
\end{array}
\eenn}
We obtain a system of the form~\eqref{eq:OscNet} by 
dividing the fast equations~\eqref{eq:OscSFNet} by~$\eps$ and setting 
$\eta = {\delta}/{\eps}$. By fixing~$\eps$ we obtain 
an~$\eta_0(\eps)$ such that~$\cP_{\textnormal{net}}$ holds 
for all $\eta<\eta_0$. Hence, there is~$\delta_0(\eps)$ such 
that~$\cP_{\textnormal{net}}$ holds for $\delta<\delta_0(\eps)$ 
in~\eqref{eq:OscSFNet}. This leads to the classification of the 
parameter space~$\cK$ into a region~(I) where~$\cP_{\textnormal{net}}$ 
holds and its complement~(II).
Depending on the coupling function~$h$, we may have 
$\lim_{\eps\to 0}\delta_0(\eps)\neq 0$ (for example if 
$h=0$ as mentioned above). However, for a generic interaction 
function one would expect $\eta_0(\eps)<C$ for some 
constant~$C$. \textcolor{black}{In this case, we have $\lim_{\eps\to 0}
\delta_0(\eps) = 0$. The resulting classification diagram is sketched in Fig.~\ref{fig:CoupledOsc}.}

Izhikevich~\cite{Izhikevich2000} derived explicit expressions 
for the dynamics on the invariant torus in the relaxation limit. 
As noted there, these expressions only describe the doubly 
singular limit for paths in parameter space converging to the 
limit point $(\eps,\delta)=(0,0)$ that lie entirely within 
region~(I). A first-order truncation of the phase dynamics---as commonly 
considered---does not describe the dynamics of the full oscillator 
network~\eqref{eq:OscSFNet} for all points in~(I) since higher-order 
terms may play a nontrivial role in the dynamics~\cite{Bick2016b,Leon2019a,Bick2021}.

While we focused on the interplay of small parameter in the intrinsic oscillator dynamics and the network coupling, interacting small parameters also arise in different ways in networked systems.
In contrast to coupled relaxation oscillators, one can also consider the case of coupled oscillators close to a Hopf bifurcation where oscillations are almost sinusoidal. Considering both small bifurcation parameter and weak coupling, one obtains explicit phase reductions~\cite{Ashwin2015a} that can---depending on the order of the approximation---contain nonpairwise interaction terms as mentioned above.
Limits involving multiple small parameters also occur if the network connections are adaptive~\cite{Gross2008}. 
This includes for example networks of neurons~\cite{Dan2004,Cooke2006,Markram2011} or adaptation in epidemic networks~\cite{Gross2006}.
Indeed, oscillator networks with adaptive interactions on have received renewed attention recently, whether the adaptation is slow (see, e.g.,~\cite{Seliger2002,Berner2019,Kasatkin2019}) or fast~\cite{Ashwin2019} relative to the oscillatory dynamics.
However, there are only few approaches taking into account distinct time scales explicitly (cf.~\cite{Jardon-Kojakhmetov2020}) in particular when multiple small parameters interact.
Thus, for adaptive networks with multiple time scales, the framework presented here may help classify the dynamics of such coupled oscillator networks.

\section{Comparison}
\label{sec:comparison}

In Section~\ref{sec:doublysingular}, we have described a wide variety of doubly-singular limit problems arising in differential equations. Yet, from the different examples, several themes emerge for the future of multiple singular limit systems.\medskip

\textbf{Property Types:} We have seen various ways of defining properties $\cP$ to obtain 
double limits which, however, share quite surprising similarities:
\begin{itemize}
 \item \textit{Individual Pattern Classification:} It turned out to be extremely useful to 
define $\cP$ via important types of patterns, e.g., the number of solutions/roots of 
an algebraic equation in Section~\ref{ssec:algebra}, the slow manifold shapes near
the transcritical point as well as the oscillation patterns 
for the Olsen model in Section~\ref{ssec:mts}, the stochastic excitable patterns for 
FitzHugh--Nagumo \textcolor{black}{SDEs} in Section~\ref{ssec:stochfss}, the types of stationary 
patterns for MEMS in Section~\ref{ssec:matched}, and the number of bifurcation points
for fast reaction PDEs in Section~\ref{ssec:fastreaction}.
 \item \textit{Phase Space Structure:} A strongly related class of properties 
emerges once one investigates pattern-forming properties more on a global level, \textcolor{black}{by} studying
the entire phase space at once. Examples are probabilistic 
quantifiers such as escape probabilities in Section~\ref{ssec:stochfss}, the sign of the 
first Lyapunov exponent in Section~\ref{ssec:shear} for oscillators with shear, or the 
global stability for linear PDMPs in Section~\ref{ssec:PDMP}.   
 \item \textit{Mathematical Features:} A last important class of properties has 
emerged corresponding to elements of proofs or mathematical properties. This includes
convexity from Section~\ref{ssec:algebra}, the exchange of partial derivatives in 
Section~\ref{ssec:analysis}, the existence of a stationary distribution in 
Section~\ref{ssec:PDMP}, or the applicability of phase reduction for networks of oscillators 
in Section~\ref{ssec:networks}.  
\end{itemize}

In view of all the preceding examples, it seems difficult to imagine that, for 
practical problems in singularly perturbed differential equations, there are highly useful
properties that do not fit within the three classes mentioned above. In fact, we see that
each class asks a different type of question, namely: How to understand individual/observed
patterns? How to understand the global structure of phase space? What are the technical 
ingredients for proofs? Looking forward, it always seems useful directly at the start of
a work on double- (or multiple-) limits to state carefully the major type of property one
is interested in for dissecting the non-negative parameter cone $\cK$.\medskip

\textbf{Diagram Structures:} Even if one has obtained a suitable partitioning of $\cK$,
one can now ask, comparing to other typical double limit problems, whether this partitioning
via $\cP$ is ``typical'' or ``common''? Quite surprisingly, a cohesive and well-founded answer 
to this problem is possible as many common features seem to emerge in 
$(\eps,\delta)\ra (0,0)$ double-limit diagrams:

\begin{itemize}
  \item \textit{Origin Ill-Posedness:} \textcolor{black}{Sometimes} it turns out that classifying
	the origin $(\eps,\delta)=(0,0)$ is ill-posed as $\cP$ is not well-defined or virtually
	impossible to evaluate at the origin. This situation \textcolor{black}{may still be} completely satisfactory from an
	applied mathematical perspective. \textcolor{black}{Indeed, if the important regime for practical applications only occurs for small positive values, and we can analyze this regime,} we do not really lose major information if we exclude
	the origin in \textcolor{black}{certain} problems.
	\item \textit{Special Axes:} Another common theme is that the two axes $\{\varepsilon=0,\delta>0\}$
	$\{\varepsilon>0,\delta=0\}$ have special or degenerate properties with respect to $\cP$. \textcolor{black}{These axes are often crucial in proofs to construct perturbation results, i.e., to infer the scaling laws in the small parameters via singular limit constructions.} Hence, it is often a suitable
	strategy to \textcolor{black}{first} understand the axes, and then aim for a perturbation, homotopy, or extension of
	the results to the interior of the cone~$\cK$.
	\item \textit{Polynomial Dissection:} As expected from classical scaling law results in
	physics as well as from the mathematical viewpoint of singularity/regularity theory, we often 
	find curves $\delta=\delta(\eps)$ (resp.~$\eps=\eps(\delta)$) with $\delta(0)=0$ 
	(resp.~$\eps(0)=0$), which provide a partitioning of the cone $\cK$. Indeed, local Taylor 
	(or H\"older-type) expansions should appear, and one can then classify the 
	partitioning of $\cK$ via the critical powers/exponents of the leading-order terms of the
	curves.
	\item \textit{Special Features:} Certain problems, either due to their inherent 
	problem formulation or due to dynamical effects, may lead to non-polynomial or otherwise 
	special dissection. Examples are exponential terms arising in stochastic metastability as well 
	as for canard problems, or curves without $\eps(0)=0$ as for fast-reaction bifurcation points.
\end{itemize}  
 
In summary, it seems clear that a complete unifying classification is impossible but in many cases
a rather exhaustive description can be provided within a common framework. First, one can aim to 
classify the behaviour on the axis for a single limit problem. Second, one can aim to obtain 
a set of (polynomial) curves partitioning the interior of $\cK$ including the leading-order scaling
exponents. Third, one aims to check whether there are any special cases occurring for the polynomial
scaling or lack of connectivity of the curves to the origin; these special cases are then treated on a
case-by-case basis and/or using a suitable shift or re-scaling to obtain polynomial order and/or
connecting curves.\medskip

\textbf{Mathematical Techniques:} Another important lesson from the comparison of the different
examples of doubly-singularly perturbed problems is that the analytical and numerical techniques 
tend to look very different at first sight. Yet, this seems to be a superficial view if one 
delves deeper into each methodology. There are many common themes appearing. First, numerical methods
tend to become more ``stiff'' near singular limits, yet analytical methods become far more feasible
the closer we are to the origin within $\cK$. This implies that a natural approach is to combine 
both approaches within $\cK$ by locally using analytical techniques and then extend the results 
beyond a small neighborhood of $(\eps,\delta)=(0,0)$ via numerical computations. Second, analytical
methods are always based upon similar principles, regardless of the differential equation studied:

\begin{itemize}
 \item \textit{Limit equations:} In a simple limit with one parameter fixed, i.e., on the coordinate
axes in the two-parameter plane, we can often obtain a reduced problem from which to start.
 \item \textit{Relative scaling:} It frequently makes sense to assume the existence of a relative
scaling $\eps=\eps(\delta)$ (resp.~$\delta=\delta(\eps)$), which provides again one-parameter
families of sub-problems lying on curves in the interior of $\cK$.
 \item \textit{Desingularization:} It often makes sense via geometric desingularization such as 
blow-up, or just via purely algebraic scaling, to generate a more complicated differential
equation, which better splits the relative scalings.
 \item \textit{Regularization:} Some problems become significantly easier if another singular
parameter is added, e.g., noise is well-known to regularize the dynamics in many instances. In fact,
we have seen this effect for excitable systems as well as for stochastically perturbed limit cycles where a non-hyperbolic structure collapses.     
\end{itemize}

In summary, also the mathematical techniques to attack very distant-looking singular perturbation
problems are more deeply related than one might anticipate. 

\section{Outlook}
\label{sec:outlook}

In this review, we have only been able to illustrate a more general framework for differential equations with multiple small parameters for certain classes of problems. It is evident that many important questions still remain. To illustrate the diversity of remaining problems, we present a few crucial questions that seem tractable within the next couple of decades:

\begin{itemize}
    \item[(Q1)] For many double-singular perturbations, multiple methodological approaches exist and we definitely need a better understanding how these approaches can be compared more directly in concrete double-limit test problems. This approach is very common in other mathematical disciplines, e.g., in numerical analysis, which often provides sharp and precise comparisons of algorithms, or even in classical analysis, where many problems involve the derivation of best-possible upper a-priori bounds. As a concrete example for the case of double limits, consider the case of multiple time scale stochastic problems discussed in Section~\ref{ssec:stochfss}. We have shown a sample-paths approach to estimate probabilities, but alternatively one could also use a distributional approach via the Fokker-Planck equation, non-autonomous dynamics techniques such as skew-product flows, quasi-stationary distributions, matched asymptotic expansions, numerical methods, as well as many other methods to 
study the double limit. The same remark applies to all other examples we have discussed. A detailed discussion of the advantages and disadvantages of every method for double limits is clearly an open problem.
    \item[(Q2)] For many double-limit problems, there are concrete conjectures left to be proven for certain regions in the two-parameter diagrams. A good example is the Olsen model in Section~\ref{ssec:mts}, where the case of non-classical relaxation oscillations is solved. Yet, rigorous proofs for mixed-mode/bursting-type oscillations as well as chaotic dynamics are missing, although the geometry of the orbits has been \textcolor{black}{well} illuminated via singular limits as well as via numerics. This is actually a common theme for all the problems, i.e., even though certain scaling regimes are tractable, it is often extremely challenging to cover the entire parameter space via rigorous proofs. An excellent goal for future research could be to develop better first-principles mathematical indicators, which tell us much quicker about the difficulties of certain scaling regions. Currently, trial-and-error is still often our best approach in this regard. 
    \item[(Q3)] Another question to follow within future work is the role played by low regularity in singular-perturbation problems. An astonishing variety of small-parameter problems in differential equations are connected to trajectories, which may have low regularity. Beyond this, even the important dynamical invariant structures (such as attractors) have low regularity. One example has been presented in Section~\ref{ssec:shear}, as shear-induced chaos for stochastic differential equations is connected to relatively rough individual sample paths and simultaneously to a chaotic attractor. Since chaotic attractors often have fractal dimension, they contribute another aspect of low regularity. In more generality, the same theme also appears for chaotic deterministic switching problems or in a completely different setting in large-scale network limits, where the regularity of the finite-dimensional problem may not always transfer to the mean-field or continuum limit.
    \item[(Q4)] From a numerical perspective, many crucial challenges are posed by double-limit dynamics. In fact, even very classical stiff differential equations with a single small parameter constitute a vast area already. Having two different, yet possibly connected, singular parameters tends to make the situation much worse. It seems wise to combine analytical pre-processing, i.e., re-writing the differential equations first into the best possible numerical problem, and careful a-priori error estimates, to avoid spurious solutions. A good example of re-writing the numerical setting has been discussed in Section~\ref{ssec:fastreaction}, where numerical continuation in the small parameters leads to well-conditioned boundary-value problems instead of quite poorly conditioned initial-value problems. It is a very worthwhile general goal to develop as many numerical methods as possible that have robustness/well-conditioning against small-parameter limits.
    \item[(Q5)] Another aspect where many open questions remain is the interplay between double limits and areas usually quite far from classical singular-limit problems for differential equations. An illustrating example are limits in coupled oscillators as discussed in Section~\ref{ssec:networks}. More generally, one can assume that the oscillators are coupled on a graph, on a simplicial complex, or a general hypergraph~\cite{Battiston2020,Bick2021}. In these cases, methods from graph theory, combinatorics, and geometry/topology are going to enter the mathematical challenge, and double-limit problems are not as classical in these areas as they are for differential equations. Yet, exploring whether it is possible to translate open questions in double-limit problems into new areas seems to be promising.
    \item[(Q6)] We have often assumed throughout this work that the studied differential equations have quite a high degree of regularity in their defining equations as this is often the most natural starting point, e.g., by invoking a more microscopic modelling approach to retain smoothness. Even in the case of SODEs with classical white noise, we have H\"older regularity in Sections~\ref{ssec:stochfss} and~\ref{ssec:shear}. Only for the PDMP case in Section~\ref{ssec:PDMP}, we have less regularity as discontinuous jumps occur. Of course, if one allows for arbitrary degeneracy in terms of input regularity, then this leads already to very intriguing questions on the level of existence of a suitable dynamical system, even for ODEs~\cite{LongoNovoObaya,BossoliniBronsKristiansen}. Already for ODEs the number of possibilities for bifurcations in non-smooth systems is extremely large~\cite{PSDS,Jeffrey} and their unfoldings via multiple small parameters involving a regularization is still under active 
development~\cite{BuzzidaSilvaTeixeira,JelbartKristiansenWechselberger}. For non-smooth SDEs and PDEs, the situation will be even more complicated. In summary, identifying principles to derive universally valid and sufficiently low-dimensional double-limit problems is already challenging once regularity assumptions are relaxed.  
    \item[(Q7)] The biggest, and \textcolor{black}{practically} most pressing, remaining challenge is to broaden the applicability of double-limit results. In fact, the steps (S1)--(S3) in the introduction apply to many other problems. For example, double-limit differential equations occur in homogenization of PDEs~\cite{Menon}, in homogenization of fast chaos~\cite{EngelGkogkasKuehn}, in rate-independent systems modeling viscoelasticity~\cite{MielkeTruskinovsky}, in bursting oscillations in neuroscience~\cite{TekaTabakBertram}, in oscillators from systems biology~\cite{MiaoPopovicSzmolyan}, in plasma physics \cite{DonatelliMarcati}, in mean-field analysis of particle systems~\cite{BodnarVelazquez}, in stochastic optimization~\cite{BorkarMitter}, and in fluid dynamics~\cite{SteinrueckSchneiderGrillhofer}. Of course, this list could be continued with many additional fields. 
    \item[(Q8)] \textcolor{black}{From a theoretical perspective, one of the most challenging conceptual open problems is how to delineate the class of singular perturbation problems, where one has to carefully apply steps (S1)--(S3), from those differential equations where direct abstract techniques allow us to neglect the small parameters easily. For ODEs, several approaches have been proposed, and one might intuitively think that it should} \textcolor{black}{be easy} \textcolor{black}{to} \textcolor{black}{ sharpen or restrict our definition of singular perturbation}\textcolor{black}{, and use this improvement to} \textcolor{black}{transfer certain results to }\textcolor{black}{ other classes of differential equations. Unfortunately, this is not} \textcolor{black}{simple}\textcolor{black}{. As} \textcolor{black}{an} \textcolor{black}{ example consider the commonly used definition that a problem is ``singular'' if a small parameter multiplies the highest derivative. Now consider an SDE. If a small parameter makes the entire drift term vanish in the limit,} \textcolor{black}{then we view the problem as singularly perturbed}. \textcolor{black}{Yet, if one re-writes the SDE via the Fokker-Planck PDE, then the drift term is generically not the highest derivative. Similar struggles appear with} \textcolor{black}{other approaches to find more restrictive} \textcolor{black}{ definitions for ``singular perturbation'' if one wants to transport them across classes.} 
\end{itemize}

Finally, we would like to point out that our general view on double-limit problems in differential equations might also have a general impact in several respects, not only within the areas of the examples we have presented, for the questions (Q1)--(Q\textcolor{black}{8}), but also well beyond:

\begin{itemize}
   \item The diagram structure, which we have utilized to summarize the main results for each case, seems to be well-adapted to the basic case of two parameters but, using suitable projections, higher-dimensional generalizations are certainly conceivable. 
   \item Although a complete classification of all possible scaling laws in all double-limit problems seems out of reach, a classification into generic cases via an abstract universality theory, analogous to critical exponents in physics, may very well exist.  
    \item It seems very promising to consistently reconsider double-limit problems that might have looked too challenging in the past. With a more coherent data base and a more structured classification, one might be able to search for new methods in virtually any other doubly-singular limit problem. 
\end{itemize}

\section*{Acknowledgments}

CK has been supported by a Lichtenberg Professorship of the VolkswagenStiftung. CK also acknowledges inspiring discussions with Grigorios A. Pavliotis regarding limit problems in differential equations, which were made possible by a TUM John von Neumann Visiting Professorship.
NB has been supported by the ANR project PERISTOCH, ANR--19--CE40--0023. CB has been supported by the Institute for Advanced Study at the Technical University of Munich through a Hans Fischer fellowship and the Engineering and Physical Sciences Research Council (EPSRC) through the grant EP/T013613/1.
ME has been supported by Germany's Excellence Strategy -- The Berlin Mathematics Research Center MATH+ (EXC-2046/1, project ID: 390685689). AI acknowledges support by an FWF Hertha Firnberg Research Fellowship (T 1199-N).
CS has received funding from the European 
Union's Horizon 2020 research and innovation program under the Marie Sk\l odowska--Curie grant 
agreement No. 754462. Support by INdAM-GNFM is gratefully acknowledged by CS. TH gratefully 
acknowledges support through SNF grant $200021-175728/1$. \textcolor{black}{We also thank two anonymous referees, whose comments and suggestions have helped to improve the presentation of this work.}

\bibliographystyle{plain}
\bibliography{refs}

\begin{thebibliography}{100}

\bibitem{Arnaudonetal18}
A.~Arnaudon, A.~L. De~Castro, and D.~D. Holm.
\newblock Noise and dissipation on coadjoint orbits.
\newblock {\em J. Nonlinear Sci.}, 28(1):91--145, 2018.

\bibitem{Arnold98}
L.~Arnold.
\newblock {\em Random Dynamical Systems}.
\newblock Springer, Berlin, 1998.

\bibitem{ArnoldSchenk96}
L.~Arnold, N.~Sri~Namachchivaya, and K.R. Schenk-Hopp\'{e}.
\newblock Toward an understanding of stochastic {H}opf bifurcation: a case
  study.
\newblock {\em Internat. J. Bifur. Chaos Appl. Sci. Engrg.}, 6(11):1947--1975,
  1996.

\bibitem{Arrhenius}
S.~Arrhenius.
\newblock On the reaction velocity of the inversion of cane sugar by acids.
\newblock {\em J.~Phys.\ Chem.}, 4:226, 1889.
\newblock In German. Translated and published in: Selected Readings in Chemical
  Kinetics, M.H. Back and K.J. Laider (eds.), Pergamon, Oxford, 1967.

\bibitem{Ashwin1992}
P.~Ashwin and J.W. Swift.
\newblock {The dynamics of $n$ weakly coupled identical oscillators}.
\newblock {\em J. Nonlinear Sci.}, 2(1):69--108, 1992.

\bibitem{Ashwin2019}
Peter Ashwin, Christian Bick, and Camille Poignard.
\newblock {State-dependent effective interactions in oscillator networks
  through coupling functions with dead zones}.
\newblock {\em Philosophical Transactions of the Royal Society A: Mathematical,
  Physical and Engineering Sciences}, 377(2160):20190042, 2019.

\bibitem{Ashwin2015}
Peter Ashwin, Stephen Coombes, and Rachel Nicks.
\newblock {Mathematical Frameworks for Oscillatory Network Dynamics in
  Neuroscience}.
\newblock {\em The Journal of Mathematical Neuroscience}, 6(1):2, 2016.

\bibitem{Ashwin2015a}
Peter Ashwin and Ana Rodrigues.
\newblock {Hopf normal form with S{\_}N symmetry and reduction to systems of
  nonlinearly coupled phase oscillators}.
\newblock {\em Physica D}, 325:14--24, 2016.

\bibitem{Bakhtin}
Y.~Bakhtin and Tobias H.
\newblock Invariant densities for dynamical systems with random switching.
\newblock {\em Nonlinearity}, 25(10):2937--2952, 2012.

\bibitem{BHLM}
Y.~Bakhtin, T.~Hurth, S.D. Lawley, and J.C. Mattingly.
\newblock Singularities of invariant densities for random switching between two
  linear odes in 2{D}.
\newblock {\em arXiv:2009.01299}, 2020.

\bibitem{Mattingly}
Y.~Bakhtin, T.~Hurth, and J.C. Mattingly.
\newblock Regularity of invariant densities for 1d-systems with random
  switching.
\newblock {\em Nonlinearity}, 28:3755--3787, 2015.

\bibitem{Mason}
M.~Balde, U.~Boscain, and P.~Mason.
\newblock A note on stability conditions for planar switched systems.
\newblock {\em Internat. J. Control}, 82(10):1882--1888, 2009.

\bibitem{Battiston2020}
Federico Battiston, Giulia Cencetti, Iacopo Iacopini, Vito Latora, Maxime
  Lucas, Alice Patania, Jean-gabriel Young, and Giovanni Petri.
\newblock {Networks beyond pairwise interactions: Structure and dynamics}.
\newblock {\em Physics Reports}, 874:1--92, 2020.

\bibitem{Baxendale_Greenwood_11}
Peter~H. Baxendale and Priscilla~E. Greenwood.
\newblock Sustained oscillations for density dependent {M}arkov processes.
\newblock {\em J. Math. Biol.}, 63(3):433--457, 2011.

\bibitem{Baxendale94}
P.H. Baxendale.
\newblock A stochastic {H}opf bifurcation.
\newblock {\em Probab. Theory Related Fields}, 99(4):581--616, 1994.

\bibitem{Baxendale03}
P.H. Baxendale.
\newblock Lyapunov exponents and stability for the stochastic {D}uffing-van der
  {P}ol oscillator.
\newblock In {\em I{UTAM} {S}ymposium on {N}onlinear {S}tochastic {D}ynamics},
  volume 110 of {\em Solid Mech. Appl.}, pages 125--135. Kluwer Acad. Publ.,
  Dordrecht, 2003.

\bibitem{Baxendale04}
P.H. Baxendale.
\newblock Stochastic averaging and asymptotic behavior of the stochastic
  {D}uffing-van der {P}ol equation.
\newblock {\em Stochastic Process. Appl.}, 113(2):235--272, 2004.

\bibitem{B17}
M.~Bena{\"\i}m.
\newblock Stochastic persistence (part {I}).
\newblock {\em Available at https://arxiv.org/abs/1806.08450}, 2018.
\newblock preprint.

\bibitem{Zitt}
M.~Bena{\"{\i}}m, S.~Le~Borgne, F.~Malrieu, and P.-A. Zitt.
\newblock On the stability of planar randomly switched systems.
\newblock {\em Ann. Appl. Probab.}, 24(1):292--311, 2014.

\bibitem{Benaim}
M.~Bena\"{\i}m, S.~Le~Borgne, F.~Malrieu, and P.-A. Zitt.
\newblock Qualitative properties of certain piecewise deterministic {M}arkov
  processes.
\newblock {\em Ann. Inst. Henri Poincar\'{e} Probab. Stat.}, 51(3):1040--1075,
  2015.

\bibitem{BeHuSt2018}
Michel Bena\"im, Tobias Hurth, and Edouard Strickler.
\newblock A user-friendly condition for exponential ergodicity in randomly
  switched environments.
\newblock {\em Electron. Commun. Probab.}, 23:1--12, 2018.

\bibitem{Le_Borgne}
Michel Bena{\"{\i}}m, St{\'e}phane Le~Borgne, Florent Malrieu, and
  Pierre-Andr{\'e} Zitt.
\newblock Quantitative ergodicity for some switched dynamical systems.
\newblock {\em Electron. Commun. Probab.}, 17:no. 56, 14, 2012.

\bibitem{Strickler_Benaim}
Michel Bena\"{\i}m and Edouard Strickler.
\newblock Random switching between vector fields having a common zero.
\newblock {\em Ann. Appl. Probab.}, 29(1):326--375, 2019.

\bibitem{Lobry}
Michel Benaïm and Claude Lobry.
\newblock Lotka–volterra with randomly fluctuating environments or “how
  switching between beneficial environments can make survival harder”.
\newblock {\em Ann. Appl. Probab.}, 26(6):3754--3785, 12 2016.

\bibitem{BenderOrszag}
C.M. Bender and S.A. Orszag.
\newblock {\em Asymptotic Methods and Perturbation Theory}.
\newblock Springer, 1999.

\bibitem{BenoitCallotDienerDiener}
E.~Beno\^{i}t, J.L. Callot, F.~Diener, and M.~Diener.
\newblock Chasse au canards.
\newblock {\em Collect. Math.}, 31:37--119, 1981.

\bibitem{BensoussanLionsPapanicolaou}
A.~Bensoussan, J.-L. Lions, and G.~Papanicolaou.
\newblock {\em Asymptotic analysis for periodic structures}.
\newblock Chelsea, 2011.

\bibitem{berardo2020interactions}
C.~Berardo, S.~Geritz, M.~Gyllenberg, and G.~Raoul.
\newblock Interactions between different predator--prey states: a method for
  the derivation of the functional and numerical response.
\newblock {\em J.~Math.~Biol.}, 80:2431--2468, 2020.

\bibitem{BG1}
N.~Berglund and B.~Gentz.
\newblock Pathwise description of dynamic pitchfork bifurcations with additive
  noise.
\newblock {\em Probab. Theory Rel.}, 122(3):341--388, 2002.

\bibitem{BG2}
N.~Berglund and B.~Gentz.
\newblock A sample-paths approach to noise-induced synchronization:
  {S}tochastic resonance in a double-well potential.
\newblock {\em {Ann.}\ {Appl.}\ {Probab.}}, 12:1419--1470, 2002.

\bibitem{BG6}
N.~Berglund and B.~Gentz.
\newblock Geometric singular perturbation theory for stochastic differential
  equations.
\newblock {\em J.~ Differ. {Equations}}, 191:1--54, 2003.

\bibitem{BGbook}
N.~Berglund and B.~Gentz.
\newblock {\em Noise-induced phenomena in slow--fast dynamical systems. A
  sample-paths approach}.
\newblock Probability and its Applications. Springer-Verlag, London, 2006.

\bibitem{BGK12}
N.~Berglund, B.~Gentz, and C.~Kuehn.
\newblock Hunting {F}rench ducks in a noisy environment.
\newblock {\em J. Differ. Equations}, 252(9):4786--4841, 2012.

\bibitem{Berglund_Gentz_Kuehn_2015}
N.~Berglund, B.~Gentz, and C.~Kuehn.
\newblock From random {P}oincar\'e maps to stochastic mixed-mode-oscillation
  patterns.
\newblock {\em J. Dyn. Differ. Equ.}, 27(1):83--136, 2015.

\bibitem{BerglundLandon}
N.~Berglund and D.~Landon.
\newblock Mixed-mode oscillations and interspike interval statistics in the
  stochastic {F}itz{H}ugh-{N}agumo model.
\newblock {\em Nonlinearity}, 25(8):2303--2335, 2012.

\bibitem{Berner2019}
R.~Berner, E.~Sch{\"{o}}ll, and S.~Yanchuk.
\newblock Multiclusters in networks of adaptively coupled phase oscillators.
\newblock {\em SIAM J. Appl. Dyn. Sys.}, 18(4):2227--2266, 2019.

\bibitem{Beutler}
G.~Beutler.
\newblock {\em Methods of Celestial Mechanics ({volume I}): physical,
  mathematical, and numerical principles}.
\newblock Springer, 2004.

\bibitem{Bick2016b}
C.~Bick, P.~Ashwin, and A.~Rodrigues.
\newblock Chaos in generically coupled phase oscillator networks with
  nonpairwise interactions.
\newblock {\em Chaos}, 26(9):094814, 2016.

\bibitem{Bick2011}
C.~Bick, M.~Timme, D.~Paulikat, D.~Rathlev, and P.~Ashwin.
\newblock Chaos in symmetric phase oscillator networks.
\newblock {\em Phys. Rev. Lett.}, 107(24):244101, 2011.

\bibitem{Bick2018c}
Christian Bick, Marc Goodfellow, Carlo~R. Laing, and Erik~A. Martens.
\newblock {Understanding the dynamics of biological and neural oscillator
  networks through exact mean-field reductions: a review}.
\newblock {\em The Journal of Mathematical Neuroscience}, 10(1):9, 2020.

\bibitem{Bick2021}
Christian Bick, Elizabeth Gross, Heather~A. Harrington, and Michael~T. Schaub.
\newblock {What are higher-order networks?}
\newblock {\em arXiv:2104.11329}, apr 2021.

\bibitem{Blackbeardetal2011}
N.~Blackbeard, H.~Erzgr\"{a}ber, and S.~Wieczorek.
\newblock Shear-induced bifurcations and chaos in models of three coupled
  lasers.
\newblock {\em SIAM J. Appl. Dyn. Syst.}, 10(2):469--509, 2011.

\bibitem{Blumenthaletal17}
A.~Blumenthal, J.~Xue, and L.-S. Young.
\newblock Lyapunov exponents for random perturbations of some area-preserving
  maps including the standard map.
\newblock {\em Ann. of Math. (2)}, 185(1):285--310, 2017.

\bibitem{Blumenthaletal18}
A.~Blumenthal, J.~Xue, and L.-S. Young.
\newblock Lyapunov exponents and correlation decay for random perturbations of
  some prototypical 2{D} maps.
\newblock {\em Comm. Math. Phys.}, 359(1):347--373, 2018.

\bibitem{BodnarVelazquez}
M.~Bodnar and J.J.L. Vel{\'a}zquez.
\newblock An integro-differential equation arising as a limit of individual
  cell-based models.
\newblock {\em J. Differen. Equat.}, 222(2):341--380, 2006.

\bibitem{BorkarMitter}
V.S. Borkar and S.K. Mitter.
\newblock A strong approximation theorem for stochastic recursive algorithms.
\newblock {\em J. Optim. Theor. Appl.}, 100(3):499--513, 1999.

\bibitem{Borowski_Kuske_etal_2011}
Peter Borowski, Rachel Kuske, Yue-Xian Li, and Juan Luis~Cabrera.
\newblock Characterizing mixed mode oscillations shaped by noise and
  bifurcation structure.
\newblock {\em Chaos}, 20(4):043117, 2010.

\bibitem{BossoliniBronsKristiansen}
E.~Bossolini, M.~Br{\o}ns, and K.U. Kristiansen.
\newblock A stiction oscillator with canards: on piecewise smooth nonuniqueness
  and its resolution by regularization using geometric singular perturbation
  theory.
\newblock {\em SIAM Rev.}, 62(4):869--897, 2020.

\bibitem{bothe2003reaction}
D.~Bothe and D.~Hilhorst.
\newblock A reaction--diffusion system with fast reversible reaction.
\newblock {\em J. Math. Anal. Appl.}, 286(1):125--135, 2003.

\bibitem{bothe2012instantaneous}
D.~Bothe and M.~Pierre.
\newblock The instantaneous limit for reaction-diffusion systems with a fast
  irreversible reaction.
\newblock {\em Discrete Contin. Dyn. Syst. Ser. S}, 5(1):49, 2012.

\bibitem{BoudjellabaSari}
H.~Boudjellaba and T.~Sari.
\newblock Dynamic transcritical bifurcations in a class of slow-fast
  predator-prey models.
\newblock {\em J. Diff. Eq.}, 246:2205--2225, 2009.

\bibitem{BredenEngel}
M.~Breden and M.~Engel.
\newblock Computer-assisted proof of shear-induced chaos in stochastically
  perturbed {H}opf systems.
\newblock {\em arXiv:2101.01491}, pages 1--39, 2020.

\bibitem{breden2019influence}
M.~Breden, C.~Kuehn, and C.~Soresina.
\newblock On the influence of cross-diffusion in pattern formation.
\newblock {\em J. Comput. Dyn.}, 8(2):213--240, 2021.

\bibitem{brocchieri2020evolution}
E.~Brocchieri, L.~Corrias, H.~Dietert, and Y.-J. Kim.
\newblock Evolution of dietary diversity and a starvation driven
  cross-diffusion system as its singular limit.
\newblock {\em arXiv preprint arXiv:2011.10304}, 2020.

\bibitem{BuzzidaSilvaTeixeira}
C.A. Buzzi, P.R. da~Silva, and M.A. Teixeira.
\newblock A singular approach to discontinuous vector fields on the plane.
\newblock {\em J. Diff. Eq.}, 231:633--655, 2006.

\bibitem{CardinTeixeira}
P.T. Cardin and M.A. Teixeira.
\newblock Fenichel theory for multiple time scale singular perturbation
  problems.
\newblock {\em SIAM J. Appl. Dyn. Syst.}, 16(3):1425--1452, 2017.

\bibitem{Chen_1974}
J.~Chen and R.E. O’Malley, Jr.
\newblock On the asymptotic solution of a two-parameter boundary value problem
  of chemical reactor theory.
\newblock {\em {SIAM} J. Appl. Math.}, 26(4):717--729, 1974.

\bibitem{Cloez}
B.~Cloez and M.~Hairer.
\newblock Exponential ergodicity for {M}arkov processes with random switching.
\newblock {\em Bernoulli}, 21:505--536, 2015.

\bibitem{conforto2014rigorous}
F.~Conforto and L.~Desvillettes.
\newblock Rigorous passage to the limit in a system of reaction--diffusion
  equations towards a system including cross diffusions.
\newblock {\em Commun. Math. Sci.}, 12(3):457--472, 2014.

\bibitem{conforto2018reaction}
F.~Conforto, L.~Desvillettes, and C.~Soresina.
\newblock About reaction--diffusion systems involving the {H}olling-type {II}
  and the {B}eddington--{D}e{A}ngelis functional responses for predator--prey
  models.
\newblock {\em Nonlinear Differ. Equ. Appl.}, 25(3):24, 2018.

\bibitem{Cooke2006}
S.~F. Cooke and T.~V.~P. Bliss.
\newblock {Plasticity in the human central nervous system}.
\newblock {\em Brain}, 129(7):1659--1673, 2006.

\bibitem{Crauel_Flandoli_98}
Hans Crauel and Franco Flandoli.
\newblock Additive noise destroys a pitchfork bifurcation.
\newblock {\em J. Dynam. Differential Equations}, 10(2):259--274, 1998.

\bibitem{czapla}
Dawid Czapla, Katarzyna Horbacz, and Hanna Wojew\'odka-\'Sciazko.
\newblock On absolute continuity of invariant measures associated with a
  piecewise-deterministic {M}arkov processes with random switching between
  flows.
\newblock {\em Available at https://arxiv.org/abs/2004.06798}, 2021.

\bibitem{DaPrato}
G.~Da~Prato and J.~Zabczyk.
\newblock {\em Ergodicity for infinite-dimensional systems}, volume 229 of {\em
  London Mathematical Society Lecture Note Series}.
\newblock Cambridge University Press, Cambridge, 1996.

\bibitem{Dan2004}
Yang Dan and Mu-ming Poo.
\newblock {Spike Timing-Dependent Plasticity of Neural Circuits}.
\newblock {\em Neuron}, 44(1):23--30, 2004.

\bibitem{Davis_article}
M.H.A. Davis.
\newblock Piecewise-deterministic {M}arkov processes: a general class of
  nondiffusion stochastic models.
\newblock {\em J. Roy. Statist. Soc. Ser. B}, 46(3):353--388, 1984.
\newblock With discussion.

\bibitem{DegnOlsenPerram}
H.~Degn, L.F. Olsen, and J.W. Perram.
\newblock Bistability, oscillation, and chaos in an enzyme reaction.
\newblock {\em Ann. N. Y. Acad. Sci.}, 316(1):623--637, 1979.

\bibitem{Desrochesetal}
M.~Desroches, J.~Guckenheimer, C.~Kuehn, B.~Krauskopf, H.~Osinga, and
  M.~Wechselberger.
\newblock Mixed-mode oscillations with multiple time scales.
\newblock {\em SIAM Rev.}, 54(2):211--288, 2012.

\bibitem{DesrochesKrauskopfOsinga1}
M.~Desroches, B.~Krauskopf, and H.M. Osinga.
\newblock {The geometry of mixed-mode oscillations in the Olsen model for the
  perioxidase-oxidase reaction}.
\newblock {\em DCDS-S}, 2(4):807--827, 2009.

\bibitem{desvillettes2019non}
L.~Desvillettes and C.~Soresina.
\newblock Non-triangular cross-diffusion systems with predator--prey reaction
  terms.
\newblock {\em Ric. Mat.}, 68(1):295--314, 2019.

\bibitem{desvillettes2015new}
L.~Desvillettes and A.~Trescases.
\newblock New results for triangular reaction cross diffusion system.
\newblock {\em J. Math. Anal. Appl.}, 430(1):32--59, 2015.

\bibitem{DeVilleSriRapti11}
L.~DeVille, N.~Sri~Namachchivaya, and Z.~Rapti.
\newblock Stability of a stochastic two-dimensional non-{H}amiltonian system.
\newblock {\em SIAM J. Appl. Math.}, 71(4):1458--1475, 2011.

\bibitem{PSDS}
M.~di~Bernardo, C.J. Budd, A.R. Champneys, and P.~Kowalczyk.
\newblock {\em Piecewise-smooth Dynamical Systems}, volume 163 of {\em Applied
  Mathematical Sciences}.
\newblock Springer, 2008.

\bibitem{Ditlevsen_Greenwood_12}
Susanne Ditlevsen and Priscilla Greenwood.
\newblock The {M}orris-{L}ecar neuron model embeds a leaky integrate-and-fire
  model.
\newblock {\em Journal of Mathematical Biology}, 67(2):239--259, 2013.

\bibitem{DoanEngelLambRasmussen}
T.S. Doan, M.~Engel, J.S.W. Lamb, and M.~Rasmussen.
\newblock Hopf bifurcation with additive noise.
\newblock {\em Nonlinearity}, 31(10):4567--4601, 2018.

\bibitem{DonatelliMarcati}
D.~Donatelli and P.~Marcati.
\newblock A quasineutral type limit for the {Navier-Stokes-Poisson} system with
  large data.
\newblock {\em Nonlinearity}, 21(1):135--148, 2008.

\bibitem{Du93}
F.~Dumortier.
\newblock Techniques in the theory of local bifurcations: blow-up, normal
  forms, nilpotent bifurcations, singular perturbations.
\newblock In {\em Bifurcations and Periodic Orbits of Vector Fields
  ({M}ontreal, {PQ}, 1992)}, volume 408 of {\em NATO Adv. Sci. Inst. Ser. C
  Math. Phys. Sci.}, pages 19--73. Kluwer Acad. Publ., Dordrecht, 1993.

\bibitem{DR96}
F.~Dumortier and R.~Roussarie.
\newblock Canard cycles and center manifolds.
\newblock {\em Mem. Amer. Math. Soc.}, 121(577):x+100, 1996.

\bibitem{eliavs2021singular}
J.~Elia{\v{s}}, D.~Hilhorst, M.~Mimura, and Y.~Morita.
\newblock Singular limit for a reaction-diffusion-{ODE} system in a neolithic
  transition model.
\newblock {\em J~Differ.~Equ.}, 295:39--69, 2021.

\bibitem{eliavs2018well}
J.~Elia{\v{s}}, M.H. Kabir, and M.~Mimura.
\newblock On the well-posedness of a dispersal model for farmers and
  hunter--gatherers in the {N}eolithic transition.
\newblock {\em Mathematical Models and Methods in Applied Sciences},
  28(02):195--222, 2018.

\bibitem{EngelGkogkasKuehn21}
M.~Engel, M.A. Gkogkas, and C.~Kuehn.
\newblock Homogenization of coupled fast-slow systems via intermediate
  stochastic regularization.
\newblock {\em J. Stat. Phys.}, April 2021. [Online]. doi:
  https://doi.org/10.1007/s10955-021-02765-7.

\bibitem{EngelKuehn}
M.~Engel and C.~Kuehn.
\newblock Discretized fast-slow systems near transcritical singularities.
\newblock {\em Nonlinearity}, 32(7):2365--2391, 2019.

\bibitem{EngelLambRasmussen1}
M.~Engel, J.S.W. Lamb, and M.~Rasmussen.
\newblock Bifurcation analysis of a stochastically driven limit cycle.
\newblock {\em Comm. Math. Phys.}, 365(3):935--942, 2019.

\bibitem{EngelLambRasmussen2}
M.~Engel, J.S.W. Lamb, and M.~Rasmussen.
\newblock Conditioned {L}yapunov exponents for random dynamical systems.
\newblock {\em Trans. Amer. Math. Soc.}, 372(9):6343--6370, 2019.

\bibitem{Faggionato}
A.~Faggionato, D.~Gabrielli, and M.~Ribezzi~Crivellari.
\newblock Non-equilibrium thermodynamics of piecewise deterministic {M}arkov
  processes.
\newblock {\em J. Stat. Phys.}, 137(2):259--304, 2009.

\bibitem{Gabrielli}
A.~Faggionato, D.~Gabrielli, and M.~Ribezzi~Crivellari.
\newblock Averaging and large deviation principles for fully-coupled piecewise
  deterministic {M}arkov processes and applications to molecular motors.
\newblock {\em {M}arkov Process. Relat.}, 16:497--548, 2010.

\bibitem{Farandaetal17}
D.~Faranda, Y.~Sato, B.~Saint-Michel, C.~Wiertel, V.~Padilla, B.~Dubrulle, and
  F.~Daviaud.
\newblock Stochastic chaos in a turbulent swirling flow.
\newblock {\em Phys. Rev. Lett.}, 119:014502, Jul 2017.

\bibitem{Fenichel4}
N.~Fenichel.
\newblock Geometric singular perturbation theory for ordinary differential
  equations.
\newblock {\em J. Differential Equat.}, 31:53--98, 1979.

\bibitem{FitzHugh1961}
R.~FitzHugh.
\newblock Impulses and physiological states in theoretical models of nerve
  membrane.
\newblock {\em Biophys. J.}, 1(6):445--466, 1961.

\bibitem{Freidlin1}
M.I. Freidlin.
\newblock Quasi-deterministic approximation, metastability and stochastic
  resonance.
\newblock {\em Physica~D}, 137:333--352, 2000.

\bibitem{FreidlinWentzell_book}
M.I. Freidlin and A.D. Wentzell.
\newblock {\em Random perturbations of dynamical systems}, volume 260 of {\em
  Grundlehren der Mathematischen Wissenschaften [Fundamental Principles of
  Mathematical Sciences]}.
\newblock Springer-Verlag, New York, second edition, 1998.
\newblock Translated from the 1979 Russian original by Joseph Sz{\"u}cs.

\bibitem{Freund}
D.D. Freund.
\newblock A note on {Kaplun} limits and double asymptotics.
\newblock {\em Proc. Amer. Math. Soc.}, 35(2):464--470, 1972.

\bibitem{geritz2012mechanistic}
S.~Geritz and M.~Gyllenberg.
\newblock A mechanistic derivation of the {DeAngelis--Beddington} functional
  response.
\newblock {\em J. Theor. Biol.}, 314:106--108, 2012.

\bibitem{Ginoux2012}
J.M. Ginoux and C.~Letellier.
\newblock {Van der Pol and the history of relaxation oscillations: Toward the
  emergence of a concept}.
\newblock {\em Chaos}, 22(2):023120, 2012.

\bibitem{Gracia_2006}
J.L. Gracia, E.~O'Riordan, and M.L. Pickett.
\newblock A parameter robust second order numerical method for a singularly
  perturbed two-parameter problem.
\newblock {\em Appl. Numer. Math.}, 56(7):962--980, 2006.

\bibitem{Gross2008}
Thilo Gross and Bernd Blasius.
\newblock {Adaptive coevolutionary networks: a review}.
\newblock {\em Journal of The Royal Society Interface}, 5(20):259--271, 2008.

\bibitem{Gross2006}
Thilo Gross, Carlos J.~Dommar D'Lima, and Bernd Blasius.
\newblock {Epidemic Dynamics on an Adaptive Network}.
\newblock {\em Physical Review Letters}, 96(20):208701, 2006.

\bibitem{GuckenheimerWechselbergerYoung}
J.~Guckenheimer, M.~Wechselberger, and L.-S. Young.
\newblock Chaotic attractors of relaxation oscillations.
\newblock {\em Nonlinearity}, 19:701--720, 2006.

\bibitem{SteinrueckSchneiderGrillhofer}
W.~Schneider H.~Steinr{\"u}ck and W.~Grillhofer.
\newblock A multiple scales analysis of the undular hydraulic jump in turbulent
  open channel flow.
\newblock {\em Fluid Dyn. Res.}, 33(1):41--55, 2003.

\bibitem{Haberman}
R.~Haberman.
\newblock Slowly varying jump and transition phenomena associated with
  algebraic bifurcation problems.
\newblock {\em SIAM J. Appl. Math.}, 37(1):69--106, 1979.

\bibitem{Haiduc1}
R.~Haiduc.
\newblock {Horseshoes in the forced van der Pol system}.
\newblock {\em Nonlinearity}, 22:213--237, 2009.

\bibitem{henneke2016fast}
F.~Henneke and B.Q. Tang.
\newblock Fast reaction limit of a volume--surface reaction--diffusion system
  towards a heat equation with dynamical boundary conditions.
\newblock {\em Asymptotic Analysis}, 98(4):325--339, 2016.

\bibitem{Herceg_2011}
D.~Herceg.
\newblock Fourth-order finite-difference method for boundary value problems
  with two small parameters.
\newblock {\em Appl. Math. Comput.}, 218(2):616--627, 2011.

\bibitem{hilhorst2007mathematical}
D.~Hilhorst, J.R. King, and M.~R{\"o}ger.
\newblock Mathematical analysis of a model describing the invasion of bacteria
  in burn wounds.
\newblock {\em Nonlinear Anal. Theory Methods Appl.}, 66(5):1118--1140, 2007.

\bibitem{hilhorst2009fast}
D.~Hilhorst, M.~Mimura, and H.~Ninomiya.
\newblock Fast reaction limit of competition-diffusion systems.
\newblock {\em Handbook of differential equations: evolutionary equations},
  5:105--168, 2009.

\bibitem{HitczenkoMedvedev_09}
P.~Hitczenko and G.S. Medvedev.
\newblock Bursting oscillations induced by small noise.
\newblock {\em SIAM J. Appl. Math.}, 69:1359--1392, 2009.

\bibitem{Holmes}
M.H. Holmes.
\newblock {\em Introduction to Perturbation Methods}.
\newblock Springer, 1995.

\bibitem{Holmes5}
P.~Holmes.
\newblock Poincar{\'e}, celestial mechanics, dynamical-systems theory and
  {``chaos''}.
\newblock {\em Phys. Rep.}, 193(3):137--163, 1990.

\bibitem{Hoppensteadt1997}
F.C. Hoppensteadt and E.M. Izhikevich.
\newblock {\em Weakly Connected Neural Networks}, volume 126 of {\em Applied
  Mathematical Sciences}.
\newblock Springer, New York, NY, 1997.

\bibitem{HorsthemkeLefever}
W.~Horsthemke and R.~Lefever.
\newblock {\em Noise-Induced Transitions}.
\newblock Springer, 2006.

\bibitem{HughesProctor2}
D.~W. Hughes and M.~R.~E. Proctor.
\newblock Chaos and the effect of noise in a model of three-wave model
  coupling.
\newblock {\em Phys. D}, 46(2):163--176, 1990.

\bibitem{HughesProctor}
D.~W. Hughes and M.~R.~E. Proctor.
\newblock A low-order model of the shear instability of convection: chaos and
  the effect of noise.
\newblock {\em Nonlinearity}, 3(1):127--153, 1990.

\bibitem{huisman1997formal}
G.~Huisman and R.J. De~Boer.
\newblock A formal derivation of the {Beddington} functional response.
\newblock {\em J. Theor. Biol.}, 185(3):389--400, 1997.

\bibitem{Hurth_Kuehn}
T.~Hurth and C.~Kuehn.
\newblock Random switching near bifurcations.
\newblock {\em Stoch. Dyn.}, 20(2):2050008, 28, 2020.

\bibitem{Huygens1888}
Christiaan Huygens.
\newblock {\em {Oeuvres compl{\`{e}}tes de Christiaan Huygens. Publi{\'{e}}es
  par la Soci{\'{e}}t{\'{e}} hollandaise des sciences.}}
\newblock M. Nijhoff, La Haye, 1888.

\bibitem{iida2006diffusion}
M.~Iida, M.~Mimura, and H.~Ninomiya.
\newblock Diffusion, cross-diffusion and competitive interaction.
\newblock {\em J. Math. Biol.}, 53(4):617--641, 2006.

\bibitem{ImkellerLederer99}
P.~Imkeller and C.~Lederer.
\newblock An explicit description of the {L}yapunov exponents of the noisy
  damped harmonic oscillator.
\newblock {\em Dynamics and Stability of Systems}, 14(4):385--405, 1999.

\bibitem{ImkellerLederer2001}
P.~Imkeller and C.~Lederer.
\newblock Some formulas for {L}yapunov exponents and rotation numbers in two
  dimensions and the stability of the harmonic oscillator and the inverted
  pendulum.
\newblock {\em Dynam. Syst.}, 16(1):29--61, 2001.

\bibitem{Iuorio_2019}
A.~Iuorio, N.~Popovi{\'{c}}, and P.~Szmolyan.
\newblock Singular perturbation analysis of a regularized {MEMS} model.
\newblock {\em {SIAM} J. Appl. Dyn. Syst.}, 18(2):661--708, jan 2019.

\bibitem{Izhikevich2000}
E.M. Izhikevich.
\newblock Phase equations for relaxation oscillators.
\newblock {\em SIAM J. Appl. Math.}, 60(5):1789--1804, 2000.

\bibitem{izuhara2008reaction}
H.~Izuhara and M.~Mimura.
\newblock Reaction--diffusion system approximation to the cross-diffusion
  competition system.
\newblock {\em Hiroshima Math. J.}, 38(2):315--347, 2008.

\bibitem{DeJagerFuru}
E.M.~De Jager and J.~Furu.
\newblock {\em The Theory of Singular Perturbations}.
\newblock North-Holland, 1996.

\bibitem{Jansons_Lythe98}
Kalvis~M. Jansons and G.~D. Lythe.
\newblock Stochastic calculus: application to dynamic bifurcations and
  threshold crossings.
\newblock {\em J. Statist. Phys.}, 90(1--2):227--251, 1998.

\bibitem{Jardonetal}
H.~Jardon-Kojakhmetov, C.~Kuehn, M.~Sensi, and A.~Pugliese.
\newblock A geometric analysis of the {SIR}, {SIRS} and {SIRWS} epidemiological
  models.
\newblock {\em Nonl. Anal.: Real World Appl.}, 58:103220, 2021.

\bibitem{Jardon-Kojakhmetov2020}
Hildeberto Jard{\'{o}}n-Kojakhmetov and Christian Kuehn.
\newblock {On Fast-Slow Consensus Networks with a Dynamic Weight}.
\newblock {\em Journal of Nonlinear Science}, 30(6):2737--2786, 2020.

\bibitem{Jeffrey}
M.R. Jeffrey.
\newblock {\em Hidden Dynamics: The Mathematics of Switches, Decisions and
  Other Discontinuous Behaviour}.
\newblock Springer, 2018.

\bibitem{JelbartKristiansenWechselberger}
S.~Jelbart, K.U. Kristiansen, and M.~Wechselberger.
\newblock Singularly perturbed boundary-equilibrium bifurcations.
\newblock {\em Nonlinearity}, 34:7371--7414, 2021.

\bibitem{Jones}
C.K.R.T. Jones.
\newblock Geometric singular perturbation theory.
\newblock In {\em Dynamical Systems (Montecatini Terme, 1994)}, volume 1609 of
  {\em Lect. Notes Math.}, pages 44--118. Springer, 1995.

\bibitem{KabanovPergamenshchikov_2003}
Yuri Kabanov and Sergei Pergamenshchikov.
\newblock {\em Two-scale stochastic systems}, volume~49 of {\em Applications of
  Mathematics (New York)}.
\newblock Springer-Verlag, Berlin, 2003.
\newblock Asymptotic analysis and control, Stochastic Modelling and Applied
  Probability.

\bibitem{Kadalbajoo_2008}
M.K. Kadalbajoo and A.S. Yadaw.
\newblock {B-Spline} collocation method for a two-parameter singularly
  perturbed convection--diffusion boundary value problems.
\newblock {\em Appl. Math. Comput.}, 201(1-2):504--513, 2008.

\bibitem{Kaper}
T.J. Kaper.
\newblock An introduction to geometric methods and dynamical systems theory for
  singular perturbation problems. analyzing multiscale phenomena using singular
  perturbation methods.
\newblock In J.~Cronin and R.E. O'Malley, editors, {\em Analyzing Multiscale
  Phenomena Using Singular Perturbation Methods}, pages 85--131. Springer,
  1999.

\bibitem{Kasatkin2019}
Dmitry~V. Kasatkin, Vladimir~V. Klinshov, and Vladimir~I. Nekorkin.
\newblock {Itinerant chimeras in an adaptive network of pulse-coupled
  oscillators}.
\newblock {\em Physical Review E}, 99(2):022203, 2019.

\bibitem{KevorkianCole}
J.~Kevorkian and J.D. Cole.
\newblock {\em Multiple Scale and Singular Perturbation Methods}.
\newblock Springer, 1996.

\bibitem{Khasminskii}
R.Z Khas'minskii.
\newblock Necessary and sufficient conditions for the asymptotic stability of
  linear stochastic systems.
\newblock {\em Theory Probab. Appl.}, 12(1):144--147, 1967.

\bibitem{KooiPoggialeAugerKooijman}
B.W. Kooi, J.C. Poggiale, P.~Auger, and S.A.L.M. Kooijman.
\newblock Aggregation methods in food chains with nutrient recycling.
\newblock {\em Ecol. Model.}, 157(1):69--86, 2002.

\bibitem{Kopell1995}
N.~Kopell and D.~Somers.
\newblock Anti-phase solutions in relaxation oscillators coupled through
  excitatory interactions.
\newblock {\em J. Math. Biol.}, 33(3):261--280, 1995.

\bibitem{Kosiuk_2011}
I.~Kosiuk and P.~Szmolyan.
\newblock Scaling in singular perturbation problems: Blowing up a relaxation
  oscillator.
\newblock {\em {SIAM} J. Appl. Dyn. Syst.}, 10(4):1307--1343, 2011.

\bibitem{KosmidisPakdaman_03}
Efstratios~K. Kosmidis and K.~Pakdaman.
\newblock An analysis of the reliability phenomenon in the
  {F}itz{H}ugh--{N}agumo model.
\newblock {\em J. Comput.\ Neuroscience}, 14:5--22, 2003.

\bibitem{KrupaPopovicKopell}
M.~Krupa, N.~Popovic, and N.~Kopell.
\newblock Mixed-mode oscillations in three time-scale systems: A prototypical
  example.
\newblock {\em SIAM J. Appl. Dyn. Syst.}, 7(2):361--420, 2008.

\bibitem{KrupaPopovicKopellRotstein}
M.~Krupa, N.~Popovic, N.~Kopell, and H.G. Rotstein.
\newblock Mixed-mode oscillations in a three time-scale model for the
  dopaminergic neuron.
\newblock {\em Chaos}, 18:015106, 2008.

\bibitem{KS01}
M.~Krupa and P.~Szmolyan.
\newblock Extending geometric singular perturbation theory to nonhyperbolic
  points---fold and canard points in two dimensions.
\newblock {\em SIAM J. Math. Anal.}, 33(2):286--314, 2001.

\bibitem{KruSzm4}
M.~Krupa and P.~Szmolyan.
\newblock Extending slow manifolds near transcritical and pitchfork
  singularities.
\newblock {\em Nonlinearity}, 14:1473--1491, 2001.

\bibitem{KuehnUM}
C.~Kuehn.
\newblock Normal hyperbolicity and unbounded critical manifolds.
\newblock {\em Nonlinearity}, 27(6):1351--1366, 2014.

\bibitem{KuehnCurse}
C.~Kuehn.
\newblock The curse of instability.
\newblock {\em Complexity}, 20(6):9--14, 2015.

\bibitem{KuehnBook}
C.~Kuehn.
\newblock {\em Multiple Time Scale Dynamics}.
\newblock Springer, 2015.

\bibitem{kuehn2020numerical}
C.~Kuehn and C.~Soresina.
\newblock Numerical continuation for a fast-reaction system and its
  cross-diffusion limit.
\newblock {\em SN Partial Differ. Equ. Appl.}, 1:7, 2020.

\bibitem{KuehnSzmolyan}
C.~Kuehn and P.~Szmolyan.
\newblock Multiscale geometry of the {Olsen} model and non-classical relaxation
  oscillations.
\newblock {\em J. Nonlinear Sci.}, 25(3):583--629, 2015.

\bibitem{Kuske99}
R.~Kuske.
\newblock Probability densities for noisy delay bifurcations.
\newblock {\em J. Statist. Phys.}, 96(3--4):797--816, 1999.

\bibitem{Lawley}
S.D. Lawley, J.C. Mattingly, and M.C. Reed.
\newblock Sensitivity to switching rates in stochastically switched {ODE}s.
\newblock {\em Commun. Math. Sci.}, 12(7):1343--1352, 2014.

\bibitem{LawleyMattinglyReed2015}
Sean~D. Lawley, Jonathan~C. Mattingly, and Michael~C. Reed.
\newblock Stochastic switching in infinite dimensions with applications to
  random parabolic {PDE}.
\newblock {\em SIAM J. Math. Anal.}, 47(4):3035--3063, 2015.

\bibitem{lehtinen2019cyclic}
.~O Lehtinen and S.A.H. Geritz.
\newblock Cyclic prey evolution with cannibalistic predators.
\newblock {\em J. Theor. Biol.}, 479:1--13, 2019.

\bibitem{Leon2019a}
I.~Le{\'{o}}n and D.~Paz{\'{o}}.
\newblock {Phase reduction beyond the first order: The case of the mean-field
  complex Ginzburg--Landau equation}.
\newblock {\em Phys. Rev. E}, 100(1):012211, 2019.

\bibitem{Cui}
Dan Li, Shengqiang Liu, and Jing'an Cui.
\newblock Threshold dynamics and ergodicity of an sirs epidemic model with
  markovian switching.
\newblock {\em Journal of Differential Equations}, 263(12):8873 -- 8915, 2017.

\bibitem{Li}
Dan Li and Hui Wan.
\newblock Coexistence and exclusion of competitive {K}olmogorov systems with
  semi-{M}arkovian switching.
\newblock {\em Discrete and Continuous Dynamical Systems}, 41(9):4145--4183,
  2021.

\bibitem{LinSheaYoung09}
K.K. Lin, E.~Shea-Brown, and L.-S. Young.
\newblock Reliability of coupled oscillators.
\newblock {\em J. Nonlinear Sci.}, 19(5):497--545, 2009.

\bibitem{LinYoung08}
K.K. Lin and L.-S. Young.
\newblock Shear-induced chaos.
\newblock {\em Nonlinearity}, 21:899--922, 2008.

\bibitem{LinYoung10}
K.K. Lin and L.-S. Young.
\newblock Dynamics of periodically kicked oscillators.
\newblock {\em J. Fixed Point Theory Appl.}, 7(2):291--312, 2010.

\bibitem{Lindneretal_04}
B.~Lindner, J.~Garcia-Ojalvo, A.~Neiman, and L.~Schimansky-Geier.
\newblock Effects of noise in excitable systems.
\newblock {\em Physics Reports}, 392:321--424, 2004.

\bibitem{Lindner_Schimansky_1999}
Benjamin Lindner and Lutz Schimansky-Geier.
\newblock Analytical approach to the stochastic {F}itz{H}ugh-{N}agumo system
  and coherence resonance.
\newblock {\em Physical Review E}, 60(6):7270--7276, 1999.

\bibitem{Li14}
A.E. Lindsay, J.~Lega, and K.B. Glasner.
\newblock Regularized model of post-touchdown configurations in electrostatic
  {MEMS}: {E}quilibrium analysis.
\newblock {\em Physica D}, 280:95--108, 2014.

\bibitem{Loecherbach}
E.~L\"{o}cherbach.
\newblock Absolute continuity of the invariant measure in piecewise
  deterministic {M}arkov processes having degenerate jumps.
\newblock {\em Stochastic Process. Appl.}, 128(6):1797--1829, 2018.

\bibitem{LongoNovoObaya}
I.P. Longo, S.~Novo, and R.~Obaya.
\newblock Weak topologies for {Carath\'eodory} differential equations:
  continuous dependence, exponential dichotomy and attractors.
\newblock {\em J. Dyn. Diff. Eq.}, 31(3):1617--1651, 2019.

\bibitem{EngelGkogkasKuehn}
M.~Engel M.~Gkogkas and C.~Kuehn.
\newblock Homogenization of fully-coupled chaotic fast-slow systems via
  intermediate stochastic regularization.
\newblock {\em arXiv:2003.11297}, pages 1--29, 2020.

\bibitem{DeMaesschalckDumortier3}
P.~De Maesschalck and F.~Dumortier.
\newblock Slow-fast {Bogdanov-Takens} bifurcations.
\newblock {\em J. Diff. Eq.}, 250:1000--1025, 2011.

\bibitem{DeMaesschalckDumortierRoussarie}
P.~De Maesschalck, F.~Dumortier, and R.~Roussarie.
\newblock {\em Canard Cycles}.
\newblock Springer, 2021.

\bibitem{DeMaesschalckWechselberger}
P.~De Maesschalck and M.~Wechselberger.
\newblock Neural excitability and singular bifurcations.
\newblock {\em J. Math. Neurosci.}, 5(1):16, 2015.

\bibitem{Phu}
Florent Malrieu and Tran~Hoa Phu.
\newblock Lotka-{V}olterra with randomly fluctuating environments: a full
  description.
\newblock {\em Available at https://arxiv.org/abs/1607.04395}, 2018.

\bibitem{Malrieu_Zitt}
Florent Malrieu and Pierre-André Zitt.
\newblock On the persistence regime for lotka-volterra in randomly fluctuating
  environments.
\newblock {\em ALEA, Lat. Am. J. Probab. Math. Stat.}, 14:733--749, 2017.

\bibitem{Markram2011}
Henry Markram, Wulfram Gerstner, and Per~Jesper Sj{\"{o}}str{\"{o}}m.
\newblock {A history of spike-timing-dependent plasticity}.
\newblock {\em Frontiers in Synaptic Neuroscience}, 3:4, 2011.

\bibitem{Menon}
G.~Menon.
\newblock Gradient systems with wiggly energies and related averaging problems.
\newblock {\em Arch. Rat. Mech. Anal.}, 162(3):193--246, 2002.

\bibitem{metz2014dynamics}
J.A. Metz and O.~Diekmann.
\newblock {\em The dynamics of physiologically structured populations},
  volume~68.
\newblock Springer, 2014.

\bibitem{Meyer}
R.E. Meyer.
\newblock On the approximation of double limits by single limits and the
  {Kaplun} extension theorem.
\newblock {\em J. Inst. Math. Appl.}, 3:245--249, 1967.

\bibitem{MiaoPopovicSzmolyan}
Z.~Miao, N.~Popovi{\'c}, and P.~Szmolyan.
\newblock Oscillations in a {cAMP} signalling model for cell aggregation - a
  geometric analysis.
\newblock {\em J. Math. Anal. Appl.}, 483(1):123577, 2020.

\bibitem{MielkeTruskinovsky}
A.~Mielke and L.~Truskinovsky.
\newblock From discrete visco-elasticity to continuum rate-independent
  plasticity: rigorous results.
\newblock {\em Arch. Rat. Mech. Anal.}, 203(2):577--619, 2012.

\bibitem{MisRoz}
E.F. Mishchenko and N.Kh. Rozov.
\newblock {\em Differential Equations with Small Parameters and Relaxation
  Oscillations (translated from Russian)}.
\newblock Plenum Press, 1980.

\bibitem{MuratovVandeneijnden2007}
C.B. Muratov and E.~Vanden-Eijnden.
\newblock Noise-induced mixed-mode oscillations in a relaxation oscillator near
  the onset of a limit cycle.
\newblock {\em Chaos}, 18:015111, 2008.

\bibitem{MusokeKrauskopfOsinga}
E.~Musoke, B.~Krauskopf, and H.M. Osinga.
\newblock A surface of heteroclinic connections between two saddle slow
  manifolds in the {Olsen} model.
\newblock {\em Int. J. Bif. Chaos}, 30(16):2030048, 2020.

\bibitem{Nagumo1962}
J.~Nagumo, S.~Arimoto, and S.~Yoshizawa.
\newblock An active pulse transmission line simulating nerve axon.
\newblock {\em Proceedings of the IRE}, 50(10):2061--2070, 1962.

\bibitem{Nakao2015}
H.~Nakao.
\newblock {Phase reduction approach to synchronisation of nonlinear
  oscillators}.
\newblock {\em Contemp. Phy.}, 57(2):188--214, 2016.

\bibitem{Nayfeh1}
A.H. Nayfeh.
\newblock {\em Perturbation Methods}.
\newblock Wiley, 2004.

\bibitem{Olsen}
L.F. Olsen.
\newblock An enzyme reaction with a strange attractor.
\newblock {\em Phys. Lett. A}, 94(9):454--457, 1983.

\bibitem{OMalley25}
R.E. O'Malley.
\newblock Two-parameter singular perturbation problems for second-order
  equations.
\newblock {\em J. Math. Mach.}, 16(10):1143--1164, 1967.

\bibitem{OMalley11}
R.E. O'Malley.
\newblock On initial value problems for nonlinear systems of differential
  equations with two small parameters.
\newblock {\em Arch. Rat. Mech. Anal.}, 40(3):209--222, 1971.

\bibitem{OMalley20}
R.E. O'Malley.
\newblock Singular perturbation theory: a viscous flow out of {G\"{o}ttingen}.
\newblock {\em Ann. Rev. Fluid Mech.}, 42:1--17, 2010.

\bibitem{OMalley24}
R.E. O'Malley.
\newblock {\em Historical Developments in Singular Perturbations}.
\newblock Springer, 2014.

\bibitem{OMalley_1967}
R.E. O'Malley~Jr.
\newblock Singular perturbations of boundary value problems for linear ordinary
  differential equations involving two parameters.
\newblock {\em J. Math. Anal. Appl.}, 19(2):291--308, 1967.

\bibitem{OMalley_1974}
R.E. O'Malley~Jr.
\newblock {\em Introduction to singular perturbations}.
\newblock Academic Press, 1974.

\bibitem{OttStenlund10}
W.~Ott and M.~Stenlund.
\newblock From limit cycles to strange attractors.
\newblock {\em Commun. Math. Phys.}, 296(1):215--249, 2010.

\bibitem{PavliotisStuart}
G.A. Pavliotis and A.M. Stuart.
\newblock {\em Multiscale Methods: Averaging and Homogenization}.
\newblock Springer, 2008.

\bibitem{Pietras2019}
B.~Pietras and A.~Daffertshofer.
\newblock Network dynamics of coupled oscillators and phase reduction
  techniques.
\newblock {\em Phys. Rep.}, 819:1--105, 2019.

\bibitem{Pikovsky2003}
Arkady Pikovsky, Michael Rosenblum, and J{\"{u}}rgen Kurths.
\newblock {\em {Synchronization: A Universal Concept in Nonlinear Sciences}}.
\newblock Cambridge University Press, 2003.

\bibitem{Poincare}
H.~Poincar{\'e}.
\newblock M{\'e}moires et observations. {Sur} certaines solutions
  particuli{\`e}res du probl{\'e}me des trois corps.
\newblock {\em Bulletin Astronomique}, I(1):65--74, 1884.

\bibitem{Popovic_2004}
N.~Popovi{\'c} and P.~Szmolyan.
\newblock A geometric analysis of the {Lagerstrom} model problem.
\newblock {\em J. Differential Equations}, 199(2):290--325, 2004.

\bibitem{Prandtl}
L.~Prandtl.
\newblock {\"{U}ber Fl\"{u}ssigkeiten bei sehr kleiner Reibung}.
\newblock In {\em Verh. III - International Math. Kongress}, pages 484--491.
  Teubner, 1905.

\bibitem{Roos_2003}
H.G. Roos and Z.~Uzelac.
\newblock The {SDFEM} for a convection-diffusion problem with two small
  parameters.
\newblock {\em Comput. Methods Appl. Math.}, 3(3):443--458, 2003.

\bibitem{SandersVerhulstMurdock}
J.A. Sanders, F.~Verhulst, and J.~Murdock.
\newblock {\em Averaging Methods in Nonlinear Dynamical Systems}.
\newblock Springer, 2007.

\bibitem{Schenk96}
K.~R. Schenk-Hopp\'{e}.
\newblock Bifurcation scenarios of the noisy {D}uffing-van der {P}ol
  oscillator.
\newblock {\em Nonlinear Dynam.}, 11(3):255--274, 1996.

\bibitem{Seliger2002}
Philip Seliger, Stephen~C. Young, and Lev~S. Tsimring.
\newblock {Plasticity and learning in a network of coupled phase oscillators}.
\newblock {\em Physical Review E}, 65(4):041906, 2002.

\bibitem{shigesada1979spatial}
N.~Shigesada, K.~Kawasaki, and E.~Teramoto.
\newblock Spatial segregation of interacting species.
\newblock {\em J. Theor. Biol.}, 79(1):83--99, 1979.

\bibitem{SimpsonKuske_2011}
D.~W.~J. Simpson and R.~Kuske.
\newblock Mixed-mode oscillations in a stochastic, piecewise-linear system.
\newblock {\em Physica D}, 240:1189--1198, 2011.

\bibitem{Somers1993}
D.~Somers and N.~Kopell.
\newblock Rapid synchronization through fast threshold modulation.
\newblock {\em Biol. Cybern.}, 68(5):393--407, 1993.

\bibitem{soresina2021hopf}
C.~Soresina.
\newblock Hopf bifurcations in the {SKT} model and where to find them.
\newblock 2021.
\newblock Under review.

\bibitem{Stocks_Manella_McClintock_89}
N.~G. Stocks, R.~Manella, and P.~V.~E. McClintock.
\newblock Influence of random fluctuations on delayed bifurcations: {T}he case
  of additive white noise.
\newblock {\em Phys.\ Rev.~A}, 40:5361--5369, 1989.

\bibitem{Strickler}
Edouard Strickler.
\newblock Randomly switched vector fields sharing a zero on a common invariant
  face.
\newblock {\em Stoch. Dyn.}, 21(2):2150007, 20, 2021.

\bibitem{Strogatz2004}
S.H. Strogatz.
\newblock {\em {Sync: The Emerging Science of Spontaneous Order}}.
\newblock Penguin, 2004.

\bibitem{Swift_Hohenberg_Ahlers91}
J.~B. Swift, P.~C. Hohenberg, and Guenter Ahlers.
\newblock Stochastic {Landau} equation with time-dependent drift.
\newblock {\em Phys.\ Rev.~A}, 43:6572--6580, 1991.

\bibitem{TekaTabakBertram}
W.~Teka, J.~Tabak, and R.~Bertram.
\newblock The relationship between two fast/slow analysis techniques for
  bursting oscillations.
\newblock {\em Chaos}, 22:043117, 2012.

\bibitem{Tihonov}
A.N. Tihonov.
\newblock Systems of differential equations containing small parameters in the
  derivatives.
\newblock {\em Mat.\ Sbornik N. S.}, 31:575--586, 1952.

\bibitem{Tuckwell}
Henry~C. Tuckwell.
\newblock {\em Stochastic Processes in the Neurosciences}.
\newblock SIAM, Philadelphia, PA, 1989.

\bibitem{Valarmathi_2003}
S.~Valarmathi and N.~Ramanujam.
\newblock Computational methods for solving two-parameter singularly perturbed
  boundary value problems for second-order ordinary differential equations.
\newblock {\em Appl. Math. Comput.}, 136(2-3):415--441, 2003.

\bibitem{vanderPol1}
B.~van~der Pol.
\newblock On relaxation oscillations.
\newblock {\em Philos. Mag.}, 7:978--992, 1926.

\bibitem{VanDyke}
M.~van Dyke.
\newblock {\em Perturbation Methods in Fluid Mechanics}.
\newblock Academic Press, 1964.

\bibitem{Verhulst}
F.~Verhulst.
\newblock {\em Methods and Applications of Singular Perturbations: Boundary
  Layers and Multiple Timescale Dynamics}.
\newblock Springer, 2005.

\bibitem{Vulanovic_2001}
R.~Vulanovi{\'c}.
\newblock A higher-order scheme for quasilinear boundary value problems with
  two small parameters.
\newblock {\em Computing}, 67(4):287--303, 2001.

\bibitem{WangYoung02}
Q.~Wang and L.-S. Young.
\newblock From invariant curves to strange attractors.
\newblock {\em Commun. Math. Phys.}, 225(2):275--304, 2002.

\bibitem{WangYoung03}
Q.~Wang and L.-S. Young.
\newblock Strange attractors in periodically-kicked limit cycles and {H}opf
  bifurcations.
\newblock {\em Commun. Math. Phys.}, 240(3):509--529, 2003.

\bibitem{Wasow}
W.~Wasow.
\newblock {\em Asymptotic Expansions for Ordinary Differential Equations}.
\newblock Dover, 2002.

\bibitem{Wechselberger4}
M.~Wechselberger.
\newblock {\em Geometric Singular Perturbation Theory beyond the Standard
  Form}.
\newblock Springer, 2020.

\bibitem{WedgwoodLin13}
K.~C.~A. Wedgwood, K.~K. Lin, R.~Thul, and S.~Coombes.
\newblock Phase-amplitude descriptions of neural oscillator models.
\newblock {\em J. Math. Neurosci.}, 3, 2013.

\bibitem{Wiezcorek09}
S.~Wieczorek.
\newblock Stochastic bifurcation in noise-driven lasers and {H}opf oscillators.
\newblock {\em Phys. Rev. E}, 79:1--10, 2009.

\bibitem{Yin}
G.~George Yin and Chao Zhu.
\newblock {\em Hybrid switching diffusions}, volume~63 of {\em Stochastic
  Modelling and Applied Probability}.
\newblock Springer, New York, 2010.
\newblock Properties and applications.

\bibitem{Young08}
L.-S. Young.
\newblock Chaotic phenomena in three settings: large, noisy and out of
  equilibrium.
\newblock {\em Nonlinearity}, 21:245--252, 2008.

\end{thebibliography}

\end{document}